\documentclass[11pt]{article}
\usepackage[utf8]{inputenc}
\usepackage[T1]{fontenc}
\usepackage{amssymb,amsmath,amsthm,mathrsfs,indentfirst}
\usepackage{xcolor}
\numberwithin{equation}{section}
\newtheorem{Theorem}{Theorem}[section]
\newtheorem{Lemma}[Theorem]{Lemma}
\newtheorem{Corollary}[Theorem]{Corollary}
\newtheorem{Proposition}[Theorem]{Proposition}
\theoremstyle{definition}
\newtheorem{Remark}[Theorem]{Remark}
\newtheorem{Definition}[Theorem]{Definition}
\newtheorem{Example}[Theorem]{Example}
\newcommand{\F}{\mathbf{F}}
\newcommand{\U}{\mathbf{U}}
\newcommand{\V}{\mathbf{V}}
\newcommand{\A}{\mathbf{A}}
\newcommand{\Le}{\mathbb{L}_e}
\newcommand{\HH}{\mathbb{H}}
\newcommand{\Ae}{\mathbb{A}}
\newcommand{\Ee}{\mathbb{E}}
\newcommand{\Ff}{\mathbb{F}}
\DeclareMathOperator{\dive}{div}
\DeclareMathOperator{\Hess}{Hess}
\DeclareMathOperator{\dist}{dist}
\DeclareMathOperator{\tr}{tr}
\begin{document}
\title{\Large From global to robust exponential attractors for nonlinear elastic waves with variable Lam\'e coefficients}
\author{\normalsize Lito E. Bocanegra-Rodr\'iguez\\  {\footnotesize  Facultad de Ciencias Matem\'aticas,} \\ {\footnotesize Universidad Nacional Mayor de San Marcos,} \\ {\footnotesize Av. Venezuela S/N, 15081 Lima, Peru}  \\ \\
\normalsize  Rodrigo Nunes Monteiro \\ {\footnotesize Departamento de Matem\'atica,}
\\ {\footnotesize Universidade Estadual de Londrina,} \\ {\footnotesize 86057-970 Londrina, PR, Brazil}\\ \\
\normalsize  Paulo N. Seminario-Huertas\thanks{Corresponding author.
	\\ {\it E-mail addresses}:  lito.bocanegra@unmsm.edu.pe (L.E. Bocanegra-Rodr\'iguez, ORCID: 0000-0001-6726-8667), monteirorn@uel.br (R. Nunes Monteiro, ORCID: 0000-0002-6915-3895), paulo.seminario.huertas@upm.es (P.N. Seminario-Huertas, ORCID: 0000-0002-9467-1299)} \\ {\footnotesize Departamento de Matem\'atica Aplicada a las Tecnolog\'ias de la Informaci\'on y las Comunicaciones,} \\ {\footnotesize Universidad Polit\'ecnica de Madrid,} \\ {\footnotesize Campus de Montegancedo, 28660 Boadilla del Monte, Madrid, Spain} }
\date{}
\maketitle
\begin{abstract}
\noindent We study the long-time dynamics of nonlinear elastic waves in heterogeneous isotropic media with spatially dependent Lam\'e	coefficients and density. For nonlinearities of critical growth, we establish global well-posedness and prove that the associated dynamical system possesses a compact global attractor of finite fractal dimension, bounded in the strong phase space. In the subcritical regime and for a constant density, we also establish the existence of an exponential attractor. A key ingredient in the analysis is a modified elastic operator whose energy form remains coercive under the pointwise separation of the pressure and shear-wave speeds. This framework allows us to investigate the singular regime in which the compressional modulus degenerates and the two wave speeds coalesce. We prove that the corresponding exponential attractors form a robust family and depend H\"older continuously on the degeneration parameter with respect to the symmetric Hausdorff distance. In particular, the long-time dynamics remains stable as the elastic structure approaches this degenerate regime.
\end{abstract}
\section{Introduction}
\subsection{The mathematical model}
Elastic wave propagation in isotropic media is commonly described by the equations of linear
elasticity, where the unknown $u=u(x,t)=(u_1,u_2,u_3)^{\top}$ denotes the displacement vector of the
material. Throughout this paper, $\Omega\subset\mathbb{R}^3$ is a bounded domain with smooth boundary
$\partial\Omega$, representing the reference configuration of an elastic body.
The internal stresses generated by the deformation are characterized by the Cauchy stress tensor
\begin{equation}\label{eq:stress}
\sigma_e(u)=\lambda(x)\dive(u)I+\mu(x)\left(Ju+Ju^{\top}\right),
\end{equation}
or, equivalently,
\begin{equation*}
(\sigma_e(u))_{ij}
=\lambda(x)\dive(u)\delta_{ij}+\mu(x)\left(\frac{\partial u_i}{\partial x_j}
+\frac{\partial u_j}{\partial x_i}\right),
\end{equation*}
where $\lambda$ and $\mu$ denote the Lam\'e coefficients, $I$ is the identity matrix and $Ju$ stands
for the Jacobian matrix of $u$. The corresponding elastic restoring force is given by the divergence
of the stress tensor, namely $\Delta_eu=\nabla\cdot\sigma_e(u)$, which can be written in the form
\begin{equation}\label{eq:lame-expanded}
\Delta_e(u)=(\lambda+\mu)\nabla(\dive u)+\mu\Delta u+(\nabla\lambda)\dive u
+\left(Ju+Ju^{\top}\right)\nabla\mu.
\end{equation}
In this work, we investigate the nonlinear damped elastic wave equation
\begin{equation} \label{eq:main}
\begin{cases}
\rho(x)\, u_{tt}- \Delta_{e} (u) + \alpha  u_t + f(u) =  b(x) &\,\,\text{in}\,\, \Omega\times\mathbb{R}^{+}, \\
u=0 &\,\,\text{on}\,\,\partial \Omega\times\mathbb{R}^{+}, \\
u(x,0)=u_0(x), \, \qquad  u_t(x,0)=u_1(x)&\,\,x \in \Omega.
\end{cases}
\end{equation}

Here $\rho=\rho(x)>0$ is the space-dependent mass density, $\alpha>0$ is the frictional damping coefficient, $f(u)$ is the nonlinear force and $b=b(x)$ is a prescribed external load. Thus, \eqref{eq:main} describes the motion of an isotropic elastic body clamped along its boundary. Such models arise in continuum mechanics, structural dynamics and the theory of elastic vibrations. We refer to the classical treatises of Lam\'e \cite{Lame}, Love \cite{Love} and Kolsky \cite{Kolsky} for the derivation of \eqref{eq:stress}-\eqref{eq:lame-expanded} from the constitutive law of isotropic elasticity and to Pujol \cite{Pujol} for the seismological setting. The structural assumptions on the coefficients and the nonlinear force are specified below.

\subsection{Assumptions}
Throughout this paper, we adopt the following standing assumptions.
\begin{itemize}
\item[(H1)] \underline{Assumption on the external force.} The external force $b\in [L^2(\Omega)]^3$ does not depend on time.
\medskip

\item[(H2)] \underline{Assumptions on the elastic coefficients.} The damping coefficient satisfies
$\alpha>0$. Moreover, $\rho\in C^1(\overline{\Omega})$, $\lambda\in C^1(\overline{\Omega})$ and
$\mu\in C^2(\overline{\Omega})$ satisfy
\begin{align*}
\rho(x)>0, \,\, \mu(x) > 0 \,\,\,\,\mbox{and}\,\,\,\, \lambda(x)+ \mu(x) > 0, \,\,\,\forall x \in \overline{\Omega}.
\end{align*}
For later use, we introduce the following constants associated with the
coefficients: $\displaystyle \overline{\rho}=\max_{\overline\Omega}\rho$,
$\displaystyle \underline{\rho}=\min_{\overline\Omega}\rho$, $\displaystyle\overline{\mu}=\max_{\overline\Omega}\mu$,
$\displaystyle \underline{\mu}=\min_{\overline\Omega}\mu > 0$,
$\displaystyle \mu_1=\max_{x\in\overline{\Omega}}\frac{\big|\nabla\mu(x)\big|^{2}}{\lambda(x)+\mu(x)}$
and $\displaystyle \mu_2=\frac12\max_{x\in\overline\Omega,\;1\leqslant i,j\leqslant 3}
\left|\frac{\partial^2\mu}{\partial x_i\partial x_j}(x)\right|$, as well as
$\displaystyle \epsilon_0=\min_{x\in\overline{\Omega}}\big\{\lambda(x)+\mu(x)\big\}>0$.

\medskip
\item[(H3)] \underline{Assumptions on the nonlinear term.} The nonlinear term $f=(f_1,f_2,f_3)$ is
assumed to admit the decomposition
\begin{equation*}
f_i(u_1,u_2,u_3)=g_i(u_1,u_2,u_3) + h_i(u_i),
\end{equation*}
where $g_i \in C^1(\mathbb{R}^3)$, $h_i \in C^2(\mathbb{R})$ and
$f_i(0)=g_i(0)=h_i(0)=0$ for $i=1,2,3$. Moreover, we assume that there exists a potential function
$G \in C^2(\mathbb{R}^3)$, such that, $\nabla G=(g_1, g_2, g_3):= g$.
\medskip
\item[(H4)] \underline{Dissipativity assumptions.} There exist positive constants $M$ and $m_f$ such that
\begin{align}
f(u)\cdot u-G(u)-\sum_{i=1}^{3}\int_0^{u_i}h_i(s)\,ds
&\geqslant -M|u|^2-m_f,\,\,\,\forall\,u\in\mathbb R^3,\label{eq:diss1}\\
G(u)+\sum_{i=1}^{3}\int_0^{u_i}h_i(s)\,ds
&\geqslant-M|u|^2-m_f,\,\,\,\ \forall\,u\in\mathbb R^3.\label{eq:diss2}
\end{align}
Moreover,
\begin{equation}\label{eq:smallness}
2M+6\mu_2+\mu_1<\underline{\mu}\lambda_1,
\end{equation}
where $\lambda_1$ denotes the first eigenvalue of the Dirichlet Laplacian $-\Delta$ on $\Omega$.
\medskip
\item[(H5)] \underline{Growth assumptions.} There exist constants $1\leqslant p<3$ and $C_g>0$ such that
\begin{equation}\label{eq:growth-g}
|\nabla g_i(u)|\leqslant C_g \Big[1+|u_1|^{p-1}+|u_2|^{p-1}+|u_3|^{p-1}\Big],\,\,\, \forall\,u\in\mathbb R^3,
\end{equation}
for $i=1,2,3$. Furthermore, there exists $C_h>0$ such that
\begin{equation}\label{eq:growth-h}
|h_i''(x)|\leqslant C_h\big[1+|x|\big],\,\,\,\forall\,x\in\mathbb R,
\end{equation}
for $i=1,2,3$.
\end{itemize}

\begin{Remark}[On the smallness and growth conditions \eqref{eq:smallness}]
	The dissipativity assumptions \eqref{eq:diss1}-\eqref{eq:diss2} provide a quadratic lower bound for the nonlinear potential, with $M$ measuring its possible negative contribution. In terms of the coefficients, condition \eqref{eq:smallness} reads
	$$
	2M+3\max_{x\in\overline\Omega,\,1\leqslant i,j\leqslant3}
	\left|\frac{\partial^2\mu}{\partial x_i\partial x_j}(x)\right|+\max_{x\in\overline\Omega}
	\frac{|\nabla\mu(x)|^2}{\lambda(x)+\mu(x)}
	<\Bigl(\min_{\overline\Omega}\mu\Bigr)\lambda_1.
	$$
	The term involving $M$ is intrinsic to the nonlinearity and is already present for constant coefficients. The other two terms measure the spatial variation of $\mu$ and correspond to $6\mu_2$ and $\mu_1$, respectively. Condition \eqref{eq:smallness} ensures that the negative part of the nonlinear potential and the lower-order terms generated by $\mu$ can be controlled by the elastic energy. Both coefficient-dependent terms vanish when $\mu$ is constant. The admissible degree of heterogeneity is quantified in Remark \ref{rem:window}.
\end{Remark}

\begin{Remark}\label{rem:nonlin}
The above assumptions separate the nonlinear term into two components with distinct mathematical
roles. The coupling part $g=\nabla G$ possesses a gradient structure, which allows the construction
of the natural energy functional and plays a fundamental role in the dissipativity analysis. On the
other hand, the term $h=(h_1,h_2,h_3)$ contains the critical nonlinear behavior. The growth
assumption \eqref{eq:growth-h} imposed on $h$ is compatible with the critical Sobolev embedding in
three dimensions. Indeed, integrating twice yields $|h_i(x)|\leqslant C\big[|x|+|x|^3\big]$ and the cubic growth is exactly the threshold at which $[H_0^1(\Omega)]^3\subset[L^6(\Omega)]^3$ still
controls $h(u)$ in $[L^2(\Omega)]^3$, whereas the assumptions
\eqref{eq:diss1}-\eqref{eq:smallness} guarantee that the total energy is bounded from below and
coercive.
\end{Remark}

\subsection{Physical and analytical aspects of heterogeneous elastic media}\label{sub:physical}

The present model differs from the classical homogeneous one in allowing the Lam\'e moduli $\lambda,\mu$ and the density $\rho$ to depend on the space variable. This setting describes heterogeneous media such as stratified geological structures, welded or diffusion-bonded joints, composites and functionally graded materials. The regularity assumptions
$\rho,\lambda\in C^1(\overline\Omega)$ and
$\mu\in C^2(\overline\Omega)$ in {\rm(H2)} restrict the analysis to continuously varying media. Discontinuous coefficients lead instead to transmission problems and are not considered here. As explained below, heterogeneity makes the wave speeds position-dependent and introduces additional terms that affect the coercivity of the elastic energy.

The first effect concerns wave propagation. Since $\mu$ is the shear modulus and $\lambda+2\mu$ is the longitudinal modulus, the local speeds of shear and pressure waves are
$$
c_S(x)=\sqrt{\frac{\mu(x)}{\rho(x)}}
\,\,\,\mbox{and}\,\,\,
c_P(x)=\sqrt{\frac{\lambda(x)+2\mu(x)}{\rho(x)}}.
$$
The conditions $\rho>0$ and $\mu>0$ in {\rm(H2)} ensure that the medium has positive mass density and resists shear. Furthermore,
\begin{equation}\label{eq:speeds}
	\lambda+\mu=\rho\bigl(c_P^2-c_S^2\bigr),
\end{equation}
so the condition $\lambda+\mu>0$ means that pressure waves travel strictly faster than shear waves, uniformly on $\overline\Omega$. Section \ref{sec:robust} examines the limiting regime in which the gap between the squared wave speeds is allowed to close.

The second effect concerns the coercivity of the elastic energy. For a displacement $u$, this energy is given by
\begin{equation}\label{eq:strain-energy}
	E(u)
	=
	\int_\Omega
	\Bigl[
	\frac{\lambda(x)}{2}\bigl(\dive u\bigr)^2
	+\mu(x)|\varepsilon(u)|^2
	\Bigr]dx,
\end{equation}
where
$\varepsilon(u)=\tfrac12(\nabla u+\nabla u^\top)$ is the linearized strain. It is the quadratic form associated with
$-\tfrac12\Delta_e$, since
$-\tfrac12\int_\Omega\Delta_e(u)\cdot u\,dx=E(u)$ for every
$u\in[H^2(\Omega)\cap H_0^1(\Omega)]^3$.

The energy estimates require $E(u)$ to control the displacement. More precisely, one seeks a constant $C>0$, such that
\begin{equation}\label{eq:coercivity}
	E(u)
	\geqslant
	C\|\nabla u\|_2^2,
	\,\,\, \forall\,u\in[H_0^1(\Omega)]^3.
\end{equation}
Mechanically, this means that every nonzero displacement compatible with the clamped boundary condition encounters elastic resistance. Analytically, it makes $E(u)$ equivalent to the squared norm of $[H_0^1(\Omega)]^3$. For variable $\mu$, however, \eqref{eq:coercivity} does not follow directly from the pointwise condition $\lambda+\mu>0$.

Two identities clarify this difficulty. Introducing the strain
$\varepsilon^D=\varepsilon-\tfrac13(\tr\varepsilon)I$, we have
\begin{equation}\label{eq:density-split}
\frac{\lambda}{2}\bigl(\tr\varepsilon\bigr)^2+\mu|\varepsilon|^2=
\mu|\varepsilon^D|^2
+\frac{K}{2}\bigl(\tr\varepsilon\bigr)^2,
\end{equation}
where $K=\lambda+\tfrac23\mu$ is the bulk modulus. This identity separates the shear and volumetric contributions to the energy density. On the other hand, integration by parts and the Cauchy-Schwarz inequality give
\begin{equation}\label{eq:div-grad}
	\|\dive u\|_2^2
	=
	\int_\Omega\sum_{i,j=1}^3
	\partial_j u_i\,\partial_i u_j\,dx
	\leqslant
	\|\nabla u\|_2^2,
	\,\,\, \forall\,u\in[H_0^1(\Omega)]^3.
\end{equation}
Since
$|\varepsilon(u)|^2
=\tfrac12|\nabla u|^2
+\tfrac12\sum_{i,j=1}^3\partial_j u_i\,\partial_i u_j$,
we also obtain
\begin{equation}\label{eq:korn-id}
	\int_\Omega|\varepsilon(u)|^2dx
	=
	\frac12\|\nabla u\|_2^2
	+\frac12\|\dive u\|_2^2,
	\,\,\, \forall\,u\in[H_0^1(\Omega)]^3.
\end{equation}
While \eqref{eq:density-split} describes the pointwise behavior of the energy density, \eqref{eq:korn-id} provides the global control induced by the clamped boundary condition. Their consequences for coercivity are described in the following remark.

\begin{Remark}[The three conditions on the Lam\'e moduli]
	\label{rem:speeds}
	The pointwise positivity of the energy density, the coercivity of the total energy and the hyperbolicity of the system lead to the strict hierarchy
	$$
	\boxed{
		\begin{gathered}
			3\lambda+2\mu>0\\[-1mm]
			\Updownarrow\\[-1mm]
			c_P/c_S>2/\sqrt3
	\end{gathered}}
	\,\,\,\Longrightarrow\,\,\,
	\boxed{
		\begin{gathered}
			\lambda+\mu>0\\[-1mm]
			\Updownarrow\\[-1mm]
			c_P>c_S
	\end{gathered}}
	\,\,\,\Longrightarrow\,\,\,
	\boxed{
		\begin{gathered}
			\lambda+2\mu>0\\[-1mm]
			\Updownarrow\\[-1mm]
			c_P>0
	\end{gathered}}.
	$$
	
	The first condition is the classical requirement that the bulk modulus
	$K=\lambda+\tfrac23\mu$ be positive. If
	$K\geqslant\underline K>0$ and
	$\mu\geqslant\underline\mu>0$ on $\overline\Omega$, then
	\eqref{eq:density-split} and \eqref{eq:korn-id} give
	$$
	E(u)
	\geqslant
	\frac14\min\bigl\{2\underline\mu,3\underline K\bigr\}
	\|\nabla u\|_2^2,
	\,\,\, \forall\,u\in[H_0^1(\Omega)]^3.
	$$
	Thus, $E(u)$ is coercive without any restriction on the derivatives of the coefficients and the analysis can be carried out directly with $-\Delta_e$.
	
	Assumption {\rm(H2)} corresponds to the second condition. It includes the regime $K>0$ and also permits the additional window
	$$
	1<\frac{c_P}{c_S}\leqslant\frac{2}{\sqrt3},
	\,\,\,\mbox{equivalently}\,\,\,
	-\frac13\mu<K\leqslant0,
	$$
	which is excluded by the first condition. In this window, \eqref{eq:density-split} no longer provides pointwise coercivity. The volumetric contribution vanishes when $K=0$ and becomes negative wherever $K<0$. The pointwise argument used in the first regime is therefore no longer available. Nevertheless, the clamped boundary condition restricts the admissible strain fields through \eqref{eq:korn-id} and the total energy may remain coercive even when its density is indefinite.
	
	This mechanism is explicit when $\mu\equiv\mu_0$ is constant. In this case,
	$$
	E(u)
	=
	\frac12\int_\Omega(\lambda+\mu_0)(\dive u)^2dx
	+\frac{\mu_0}{2}\|\nabla u\|_2^2
	\geqslant
	\frac12
	\Bigl[
	\mu_0-\max_{\overline\Omega}(\lambda+\mu_0)^-
	\Bigr]
	\|\nabla u\|_2^2.
	$$
	Hence the third condition $\lambda+2\mu_0>0$ implies coercivity. For constant coefficients, this threshold is sharp. Indeed, the constant $1$ in \eqref{eq:div-grad} is optimal. To see this, take
	$u=\nabla\phi$ with $0\neq\phi\in C_c^\infty(\Omega)$. Then
	$u\in[H_0^1(\Omega)]^3$, $\dive u=\Delta\phi$ and
	$\nabla u=D^2\phi$. Two integrations by parts give
	$$
	\int_\Omega(\Delta\phi)^2dx
	=
	\int_\Omega|D^2\phi|^2dx
	$$
	and hence
	$\|\dive u\|_2=\|\nabla u\|_2$. Substitution into the expression for the energy yields
	$$
	2E(u)
	=
	(\lambda+2\mu_0)\|\nabla u\|_2^2.
	$$
	Therefore, if $\lambda+2\mu_0=0$, there exists a nonzero displacement $u$ for which $E(u)=0$ and coercivity fails. If $\lambda+2\mu_0<0$, the energy is negative for such a displacement. Thus, $\lambda+2\mu_0>0$ is the sharp threshold for coercivity in the constant-coefficient case.
	
	When $\mu$ is variable, the third condition no longer gives the same estimate. Replacing $\mu(x)$ by its minimum yields
	$$
	E(u)
	\geqslant
	\frac12
	\Bigl[
	\underline\mu-
	\max_{\overline\Omega}(\lambda+\underline\mu)^-
	\Bigr]
	\|\nabla u\|_2^2,
	$$
	which requires the global condition
	$\lambda+2\underline\mu>0$. This is stronger than the pointwise condition
	$\lambda+2\mu>0$ and compares $\lambda(x)$ with the minimum of $\mu$ over the entire domain. Assumption {\rm(H2)}, in contrast, compares $\lambda(x)$ and $\mu(x)$ at the same point and includes configurations not covered by this estimate, such as that of Example \ref{rem:example}. Section \ref{sec:framework} treats this setting through an equivalent modified form that is coercive under {\rm(H2)} and the prescribed control on the derivatives of $\mu$.
	
	Finally, the equality $c_P=c_S$ does not mark the loss of hyperbolicity, which occurs when $c_P$ vanishes. It corresponds to the closing of the wave-speed gap and to the limiting case $\lambda+\mu=0$ of the strict inequality imposed in {\rm(H2)}. Section \ref{sec:robust} shows that the modified form remains coercive in this limit under a uniform compatibility condition between $\lambda+\mu$ and $\nabla\mu$.
\end{Remark}

The distinction between pointwise and global coercivity discussed above also has a physical interpretation. A negative bulk modulus is incompatible with the stability of a free homogeneous isotropic solid, which explains the classical requirement $3\lambda+2\mu>0$. Under suitable constraints, however, the total elastic energy may remain coercive even when its density is not pointwise positive. A related mechanism occurs in composites with negative-stiffness inclusions, where an unstable phase is stabilized by the surrounding matrix \cite{lakes,lakesdrugan}. In the present problem, the stabilizing constraint is provided by the clamped boundary condition and expressed quantitatively by \eqref{eq:korn-id}.

A further motivation for studying the closing of the wave-speed gap comes from the elastic anomalies observed during the transition between the $\alpha$ and $\beta$ crystalline phases of quartz. Across this transition, the linear incompressibilities decrease sharply while the shear moduli change only slightly \cite{ohno,peng}. In view of \eqref{eq:speeds}, this behavior suggests a decrease of
$\lambda+\mu=\rho(c_P^2-c_S^2)$ with nearly constant $\mu$, which is the regime considered in Section \ref{sec:robust} and which also motivated the homogeneous model studied in \cite{bocanegra}. We do not claim that quartz reaches the range
$c_P/c_S\leqslant2/\sqrt3$. The cited observations support only the softening of the compressional response at nearly constant shear modulus.

We now return to the analytical assumptions. The following examples show that the coefficient conditions and the smallness requirement \eqref{eq:smallness} are compatible. The first is a perturbation of a homogeneous medium in the classical regime, while the second belongs to the additional window permitted by {\rm(H2)}.

\begin{Example}\label{ex:pert}
Fix $\psi\in C^2(\overline\Omega)$ and let
$\lambda(x)=\lambda_0$ and
$\mu(x)=\mu_0[1+\eta\psi(x)]$, where $\lambda_0,\mu_0>0$. For sufficiently small $|\eta|$, one has
	$\mu>0$ and $\lambda+\mu>0$ on $\overline\Omega$. Moreover,
	$\nabla\mu=\mu_0\eta\nabla\psi$ and
	$\partial_{ij}^2\mu=\mu_0\eta\,\partial_{ij}^2\psi$, while the denominator in the definition of $\mu_1$ remains bounded away from zero. Consequently,
	$$
	\mu_1=O(\eta^2),
	\,\,\,
	\mu_2=O(|\eta|),
	\,\,\,
	\underline\mu\longrightarrow\mu_0
	\,\,\,\mbox{and}\,\,\,
	\epsilon_0\longrightarrow\lambda_0+\mu_0\,\,\,\,\, \mbox{as}\,\,\,\, \eta\to0.
	$$
	Hence {\rm(H2)} holds for sufficiently small $|\eta|$. If, in addition,
	$2M<\mu_0\lambda_1$, then \eqref{eq:smallness} also holds for all sufficiently small $|\eta|$.
\end{Example}

\begin{Example}\label{rem:example}
	The family in Example \ref{ex:pert} belongs to the regime
	$3\lambda+2\mu>0$, where the elastic energy is already coercive. We now give a family lying outside this regime. Let
	$\Omega=B_1(0)\subset\mathbb R^3$, so that $\lambda_1=\pi^2$ and fix
	$\mu_0>0$ and $\eta\in(0,1)$. Set
	\begin{equation}\label{eq:example}
		\mu(x)=\mu_0\bigl[1+\eta x_1\bigr],
		\,\,\,
		\lambda(x)=\mu_0\eta-\mu(x)
		=\mu_0\bigl[\eta-1-\eta x_1\bigr].
	\end{equation}
	Then $\lambda+\mu\equiv\mu_0\eta$ and hence
	$$
	\epsilon_0=\mu_0\eta,
	\,\,\,
	\underline\mu=\mu_0(1-\eta),
	\,\,\,
	\overline\mu=\mu_0(1+\eta),
	\,\,\,
	|\nabla\mu|=\mu_0\eta
	\,\,\,\mbox{and}\,\,\,
	\nabla^2\mu=0.
	$$
	In particular, {\rm(H2)} holds and both coefficients are smooth.
	
	Moreover,
	$3\lambda+2\mu=3\mu_0\eta-\mu$, whose maximum on
	$\overline\Omega$ is $\mu_0(4\eta-1)$. Therefore,
	$$
	3\lambda+2\mu<0
	\,\,\,\mbox{on } \overline\Omega
	\,\,\,\mbox{whenever } \eta<\frac14.
	$$
	Thus, the bulk modulus is negative throughout the medium. Equivalently,
	$$
	\max_{\overline\Omega}
	\left(\frac{c_P}{c_S}\right)^2
	=
	\frac{1}{1-\eta}
	<
	\frac43
	$$
	under the same condition.
	
	The sufficient condition
	$\lambda>-\underline\mu$ also fails, since
	$\min_{\overline\Omega}\lambda=-\mu_0$, whereas
	$-\underline\mu=-\mu_0(1-\eta)>-\mu_0$. On the other hand,
	$\mu_1=\mu_0\eta$ and $\mu_2=0$, so \eqref{eq:smallness} becomes
	$$
	2M+\mu_0\eta
	<
	\mu_0(1-\eta)\pi^2.
	$$
	Hence, if $2M<\mu_0\pi^2$, this condition holds whenever
	$$
	0<\eta<
	\min\left\{
	\frac14,
	\frac{\mu_0\pi^2-2M}{\mu_0(1+\pi^2)}
	\right\}.
	$$
	Thus, \eqref{eq:smallness} applies to media that are not covered by either
	$3\lambda+2\mu>0$ or $\lambda>-\underline\mu$. This is precisely the regime that motivates the modified operator.
\end{Example}

The physical interpretation of the parameter-dependent problem is discussed in Subsection \ref{sec:physics}, where the degeneration of the wave-speed gap is described in terms of the bulk and shear moduli.

\subsection{Main results and contributions}
This paper investigates the long-time dynamics of a semilinear Lam\'e system with spatially varying Lam\'e coefficients and density. Our main results are summarized below.
\begin{itemize}
	\item[(i)] We establish the global well-posedness of problem
	\eqref{eq:main} in its natural energy space, including nonlinearities
	with critical Sobolev growth. Consequently, the solutions generate a
	strongly continuous semigroup $\{S(t)\}_{t \geqslant 0}$.
	\item[(ii)] We prove that $\{S(t)\}_{t \geqslant 0}$ possesses a unique compact global
	attractor $\mathbf{A}$ with finite fractal dimension and additional regularity. The
	attractor is also characterized as the unstable manifold emanating from
	the set of stationary points. The precise statement is given in Theorem \ref{thm:attractor}.
	\item[(iii)] In the subcritical regime and under the constant-density assumption, we
	establish the existence of an exponential attractor
	$\mathbf{A}^{\mathrm{exp}}$ containing the global attractor $\mathbf{A}$,
	as stated in Theorem \ref{thm:exp}. In particular, when the Lam\'e
	coefficients are also constant, this result yields an exponential attractor for the  subcritical case of the	homogeneous semilinear Lam\'e system
	considered in \cite{bocanegra}. The roles of the constant-density
	assumption and the subcritical growth condition are discussed in Remarks
	\ref{rem:rho} and \ref{rem:subcritical}, respectively.

\item[(iv)] We study a family of heterogeneous media for which
$\lambda_\delta+\mu$ converges to a possibly degenerate limit
$\lambda_0+\mu\geqslant0$, allowing the pressure and shear-wave speeds to coincide on a transition zone. Theorem \ref{thm:robust} establishes the existence of exponential attractors
$\{\mathbf A_\delta^{\exp}\}_{\delta\in[0,1]}$ with uniformly bounded fractal dimension and H\"older continuous dependence on $\delta$ in the symmetric Hausdorff distance. Theorem \ref{thm:usc} further proves the upper semicontinuity of the corresponding global attractors at $\delta=0$.
\end{itemize}

The proofs rely on the following analytical ingredients developed for the variable-coefficient setting.
\begin{itemize}
\item[(a)] We introduce in Section \ref{sec:framework} the \emph{modified Lam\'e operator} $\Le$, whose associated bilinear form is coercive on $[H_0^1(\Omega)]^3$ under the pointwise condition $\lambda+\mu>0$ in {\rm(H2)}. For a variable shear modulus, this condition is strictly weaker than the standard assumptions
$3\lambda+2\mu>0$ and $\lambda>-\underline\mu$ and Example \ref{rem:example} exhibits coefficients outside both regimes. We prove that $\Le$ is positive and self-adjoint with compact inverse, establish the ellipticity of its principal part in Proposition \ref{prop:elliptic} and identify
$D(\Le)=[H^2(\Omega)\cap H_0^1(\Omega)]^3$ using the regularity theory of Ne\v{c}as \cite{necas}. We also formulate the problem in the weighted $L^2$ space induced by the variable density and control the lower-order terms generated by $\mu$ under \eqref{eq:smallness}.	

\item[(b)] We identify the fractional-power spaces associated with $\Le$ with the corresponding fractional domains of the Dirichlet	Laplacian. This identification, established in Corollary 	\ref{cor:fract}, transfers the Sobolev embeddings required in the analysis to the variable-coefficient setting.

\item[(c)] We develop a finite bootstrap procedure that progressively improves the regularity estimates for the trajectories and yields the smoothing property of the inhomogeneous component required for the construction of an exponential attractor. The regularity increments and the bootstrap procedure are established in Lemmas \ref{lem:bootstrap} and \ref{lem:comp}.

\item[(d)] We derive from \cite[Theorem 3.1]{FGMZ} an abstract robustness criterion for exponential attractors in fixed Banach spaces, stated as Theorem \ref{thm:robust-abstract} in the form needed here. We verify its hypotheses through uniform contraction, smoothing estimates and a quantitative comparison between the parameter-dependent semigroups.
\end{itemize}

\subsection{Related work}

The results described above lie at the intersection of several lines of research on nonlinear elastic waves and dissipative dynamical systems. We first recall the works on elastic-wave propagation that provide the analytical background for the present problem. In the whole space, the analysis has focused mainly on the existence and continuation of small solutions. John \cite{john} proved almost global existence for elastic waves arising from small initial disturbances, while Sideris \cite{Sideris} obtained global existence under a suitable null condition. In bounded domains, dispersion is no longer available and dissipation becomes the principal mechanism controlling the evolution. Yamamoto \cite{ya1988,yam1989} studied the propagation of singularities and exponential energy decay for the linear elastic system with Dirichlet boundary conditions, whereas Bellassoued \cite{bellassoued} established energy decay in the presence of a locally distributed time-dependent nonlinear damping. These works rely on the geometric control condition of Bardos, Lebeau and Rauch \cite{BLR}. In the present problem, the damping acts throughout the domain and no geometric control condition is required.

The passage from energy decay to the asymptotic description of the dynamics leads to the theory of global attractors developed in the classical works of Hale \cite{hale}, Ladyzhenskaya \cite{ladyzhenskaya} and Temam \cite{temam}. For damped scalar wave equations, Arrieta, Carvalho and Hale \cite{arrieta} established the existence of a global attractor in the presence of nonlinearities with critical Sobolev growth. Their argument decomposes the dynamics into a decaying component and a component with improved compactness. This method relies on the scalar structure of the equation and on the regularity theory of the Dirichlet Laplacian and therefore does not apply directly to the present elastic system with coupled components and spatially varying Lam\'e coefficients. An alternative is provided by the quasi-stability theory of Chueshov and Lasiecka \cite{Von,CL-memoirs}, which yields asymptotic compactness, finite fractal dimension and additional regularity of complete trajectories without requiring an intrinsic smoothing property of the semigroup. This framework has been applied to several scalar and coupled wave systems, including equations with degenerate memory studied by Cavalcanti, Fatori and Ma \cite{Ma} and Bresse systems considered by Ma and Monteiro \cite{singular}.

For semilinear elastic systems, Bocanegra-Rodr\'iguez, Jorge Silva, Ma and Seminario-Huertas \cite{bocanegra} established the existence and qualitative properties of a global attractor for a Lam\'e system with constant coefficients. They also considered the degeneration $\lambda+\mu\to0$ and proved the upper semicontinuity of the corresponding global attractors. To the best of our knowledge, this was the first work devoted specifically to the global attractor of a semilinear Lam\'e system. The present analysis starts from this model but allows the Lam\'e coefficients and the density to vary spatially. It also replaces the qualitative perturbation analysis of global attractors by a quantitative robustness result for exponential attractors, while retaining the upper semicontinuity of the global attractors.

Exponential attractors, introduced by Eden, Foias, Nicolaenko and Temam \cite{EFNT}, complement the theory of global attractors by providing a finite-dimensional attracting set together with a quantitative rate of attraction. Our construction uses the criterion of Efendiev, Miranville and Zelik \cite{EMZ}, see also \cite[\S7.9]{Von}, based on decomposing a fixed-time map into a contractive component and a smoothing component. We also use the transitivity of exponential attraction established by Fabrie, Galusinski, Miranville and Zelik \cite[Theorem~5.1]{FGMZ}, as well as their perturbation framework for robust families of exponential attractors. Particularly relevant is the work of Pata and Squassina \cite{PS}, where the decomposition method is combined with a finite regularity bootstrap for the strongly damped wave equation. A related procedure appears in \cite{singular}. We follow this general strategy, but the regularity gain must be obtained from the inhomogeneous component because the linear semigroup associated with the present frictionally damped system has no smoothing effect, as explained in Remark \ref{rem:PS}.

\medskip
\underline{Organization of the paper.}
Section \ref{sec:framework} develops the modified Lam\'e operator and the functional framework for the variable-coefficient problem. Section \ref{sec:wp} establishes global well-posedness and the basic energy estimates. In Section \ref{sec:attractor}, we prove the existence, finite fractal dimension and additional regularity of the global attractor. Section \ref{sec:exp} constructs an exponential attractor in the subcritical regime. Finally, Section \ref{sec:robust} introduces the parameter-dependent family and proves the robustness of the exponential attractors together with the upper semicontinuity of the global attractors.

\section{Functional Framework}\label{sec:framework}
\subsection{Auxiliary Spaces}
For $1\leqslant q<\infty$, let $(L^q(\Omega),\|\cdot\|_q)$ denote the usual Lebesgue space.
Throughout the paper, we equip $[L^q(\Omega)]^3$ with the norm $\|u\|_q^q=\sum_{i=1}^3\|u_i\|_{L^q(\Omega)}^q$, while $[L^\infty(\Omega)]^3$ carries $\|u\|_\infty=\max_i\|u_i\|_{L^\infty(\Omega)}$ and $[L^2(\Omega)]^3$ is endowed with the inner product $\langle u,v\rangle=\int_\Omega u\cdot v\,dx$.
Let $A=-\Delta:D(A)\subset L^2(\Omega)\to L^2(\Omega)$, where $D(A)=H^2(\Omega)\cap H_0^1(\Omega)$.
Since $A$ is a positive self-adjoint operator with compact inverse, its fractional powers $A^{s/2}$,
$s\in\mathbb R$, are well defined and $D(A^{s/2})$ is a Hilbert space with the inner product
$\langle A^{s/2}\cdot,A^{s/2}\cdot\rangle$ and the induced norm $\|A^{s/2}\cdot\|_2$. Furthermore, for $s>r$ the embedding $D(A^{s/2})\subset D(A^{r/2})$ is continuous and compact. In particular,
\begin{equation*}
\|A^{r/2} u\|_2 \leqslant \lambda_1^{(r-s)/2}\, \|A^{s/2} u\|_2, \,\,\, \forall \, u \in D(A^{s/2}).
\end{equation*}
Indeed, writing $ u=\sum_{k\in \mathbb{N}} u_k e_k$ in the orthonormal basis of eigenfunctions of $A$ and using
$\lambda_k\geqslant\lambda_1$ together with $r-s<0$,
$$
\|A^{r/2}u\|_2^2
=\sum_{k\in \mathbb{N}}\lambda_k^{\,r}|u_k|^2
=\sum_{k\in \mathbb{N}}\lambda_k^{\,r-s}\,\lambda_k^{\,s}|u_k|^2
\leqslant\lambda_1^{\,r-s}\|A^{s/2}u\|_2^2.
$$
Moreover, by interpolation, if $s>r>q$ then for every $\varepsilon>0$ there exists $C(\varepsilon)>0$ such that
\begin{equation}\label{eq:interp}
\|A^{r/2} u\|_2 \leqslant \varepsilon \|A^{s/2} u\|_2 + C(\varepsilon)\|A^{q/2} u\|_2, \,\,\, \forall \, u \in D(A^{s/2}).
\end{equation}
By the concrete characterization of the fractional domains of the Dirichlet Laplacian due to
Fujiwara \cite{fujiwara} (see also Lions-Magenes \cite[Ch. 1, Thms. 11.6 and 15.1]{LM} and, for a self-contained exposition, Hajduk-Robinson \cite{HR}), we have
\begin{equation}\label{eq:charact}
D(A^{s/2})=
\begin{cases}
H^s(\Omega), & 0\leqslant s<\tfrac12,\\
H_0^s(\Omega), & \tfrac12<s\leqslant1,\\
H^s(\Omega)\cap H_0^1(\Omega), & 1\leqslant s\leqslant2.
\end{cases}
\end{equation}
Combining \eqref{eq:charact} with the Sobolev embedding $H^{s}(\Omega)\subset L^{6/(3-2s)}(\Omega)$ valid in dimension three for $0\leqslant s<3/2$, the following continuous embedding holds
\begin{equation}\label{eq:sobolev}
D(A^{s/2}) \subset L^{q}(\Omega), \,\, 2 \leqslant  q \leqslant 6/(3-2s), \,\,
s \in \left[0, \tfrac{3}{2}\right),
\end{equation}
and $D(A^{3/4})\subset L^{q}(\Omega)$ for every
$q<\infty$. No special treatment is needed when $s=\tfrac12$.
Although $D(A^{1/4})=H^{1/2}_{00}(\Omega)$ is a proper subspace of
$H^{1/2}(\Omega)$, the continuous inclusions
$H^{1/2}_{00}(\Omega)\subset H^{1/2}(\Omega)\subset L^{3}(\Omega)$
still yield \eqref{eq:sobolev}.
\medskip
Similarly, let $\Ae= A\,I:[D(A)]^3\subset [L^2(\Omega)]^3\to[L^2(\Omega)]^3$ denote the vector
Dirichlet Laplacian, where $I$ is the $3\times3$ identity matrix. Since $A$ is a positive
self-adjoint operator with compact inverse, the same properties hold for $\Ae$. Therefore, its fractional powers $\Ae^{s/2}$ are well defined and for $u=(u_1,u_2,u_3)\in D(\Ae^{s/2})$, this space is endowed with the norm $\|u\|_{D(\Ae^{s/2})}^2
=\sum_{i=1}^3\|A^{s/2}u_i\|_{L^2(\Omega)}^2$.

\subsection{The coercive case: $3\lambda+2\mu>0$ and other elementary regimes}\label{sub:coercive}
Consider the symmetric bilinear form associated with the stress tensor \eqref{eq:stress},
\begin{equation}\label{eq:physform}
	a(u,v)
	=
	\int_\Omega
	\Bigl[
	\lambda\,\dive(u)\dive(v)
	+2\mu\,\varepsilon(u):\varepsilon(v)
	\Bigr]dx.
\end{equation}
Since $a(u,u)=2E(u)$, this form is twice the elastic energy defined in \eqref{eq:strain-energy}. Moreover, integration by parts gives
\begin{equation}\label{eq:aibp}
-\int_\Omega\Delta_e(u)\cdot v\,dx
=a(u,v), \,\,\,
u\in\bigl[H^2(\Omega)\cap H_0^1(\Omega)\bigr]^3,\,\,v\in[H_0^1(\Omega)]^3.
\end{equation}
If
$$
a(u,u)\geqslant C\|\nabla u\|_2^2,
\,\,\, \forall\,u\in[H_0^1(\Omega)]^3,
$$
then $a$ is a coercive closed form and its associated positive self-adjoint operator is $-\Delta_e$. No modification of the elastic operator is therefore required. We record below two elementary regimes in which this coercivity holds. Under the weaker condition $\lambda+\mu>0$ in {\rm(H2)}, however, it is not guaranteed when $\mu$ is variable, which motivates the construction in Subsection \ref{sub:modified}.

The first regime is determined by the positivity of the bulk modulus. In this case, the decomposition \eqref{eq:density-split} gives pointwise coercivity of the energy density, which becomes coercivity of the form $a$ through \eqref{eq:korn-id}.

\begin{Proposition}[Positive bulk modulus]\label{prop:unperturbed}
Assume {\rm(H2)} and suppose, in addition, that $3\lambda(x)+2\mu(x)>0$ for every	$x\in\overline\Omega$. Define $\underline K	=	\min_{\overline\Omega}	\left(\lambda+\tfrac23\mu\right)>0$ and $C_K=\frac12\min\bigl\{2\underline\mu,3\underline K\bigr\}>0$. Then the following properties hold.
\begin{enumerate}
\item[$1.$] The form $a$ is symmetric, bounded and coercive on $[H_0^1(\Omega)]^3$. More precisely,
$$a(u,u)\geqslant C_K\|\nabla u\|_2^2, \,\,\, \forall\,u\in[H_0^1(\Omega)]^3.$$
Thus, the coercivity constant depends only on the uniform lower bounds of the shear and bulk moduli and is independent of their spatial derivatives.
		
\item[$2.$] The operator $-\Delta_e$ associated with $a$ is strictly positive and self-adjoint, with
$$D(-\Delta_e)
=
[H^2(\Omega)\cap H_0^1(\Omega)]^3
\,\,\,\mbox{and}\,\,\,
D\bigl((-\Delta_e)^{1/2}\bigr)
=[H_0^1(\Omega)]^3.$$
\end{enumerate}
Consequently, $-\Delta_e$ has all the properties of the elastic operator required in Sections $\ref{sec:wp}$-$\ref{sec:exp}$ and may be used in place of the modified operator introduced in Subsection $\ref{sub:modified}$. The precise correspondence is given in Remark $\ref{rem:dictionary}$.
\end{Proposition}

\begin{proof}
The symmetry of $a$ follows immediately from \eqref{eq:physform}. To prove coercivity, we use the doubled form of \eqref{eq:density-split},
		$$
		\lambda\bigl(\tr\varepsilon\bigr)^2
		+2\mu|\varepsilon|^2
		=
		2\mu|\varepsilon^D|^2
		+K\bigl(\tr\varepsilon\bigr)^2.
		$$
		Since
		$|\varepsilon|^2=|\varepsilon^D|^2+\tfrac13(\tr\varepsilon)^2$, it follows that
		$$
		2\mu|\varepsilon^D|^2
		+K\bigl(\tr\varepsilon\bigr)^2
		\geqslant
		\min\bigl\{2\underline\mu,3\underline K\bigr\}
		|\varepsilon|^2.
		$$
		Taking $\varepsilon=\varepsilon(u)$, integrating over $\Omega$ and using \eqref{eq:korn-id}, we obtain
		$$
		a(u,u)
		\geqslant
		\min\bigl\{2\underline\mu,3\underline K\bigr\}
		\int_\Omega|\varepsilon(u)|^2dx
		\geqslant
		C_{K}\|\nabla u\|_2^2.
		$$
		Finally, the boundedness of $a$ follows from
		$|\varepsilon(u)|\leqslant|\nabla u|$, \eqref{eq:div-grad} and the boundedness of $\lambda$ and $\mu$ on $\overline\Omega$.

The symmetric, bounded and coercive form $a$ is densely defined on
$[L^2(\Omega)]^3$. Since its form norm is equivalent to the norm of
$[H_0^1(\Omega)]^3$, it is also closed. The first representation theorem
\cite[Ch.~VI, Thm.~2.1]{kato-book} therefore associates with $a$ a unique
m-sectorial operator $T$ on $[L^2(\Omega)]^3$ satisfying
$$
a(u,v)=\langle Tu,v\rangle,
\,\,\,
u\in D(T),
\,\,\,
v\in[H_0^1(\Omega)]^3.
$$
Since $a$ is symmetric, $T$ is self-adjoint by
\cite[Ch.~VI, Thm.~2.6]{kato-book}. Moreover, item $1$ and the Poincar\'e inequality give
$$
\langle Tu,u\rangle
=
a(u,u)
\geqslant
C_K\|\nabla u\|_2^2
\geqslant
C_K\lambda_1\|u\|_2^2,
\,\,\,
u\in D(T),
$$
so that $T$ is strictly positive. The second representation theorem
\cite[Ch. VI, Thm.~2.23]{kato-book} yields
$$
D(T^{1/2})=[H_0^1(\Omega)]^3.
$$

It remains to identify $T$ and its domain. If
$u\in[H^2(\Omega)\cap H_0^1(\Omega)]^3$, then
$\Delta_eu\in[L^2(\Omega)]^3$ by \eqref{eq:lame-expanded} and
\eqref{eq:aibp} gives
$a(u,v)=\langle-\Delta_eu,v\rangle$ for every
$v\in[H_0^1(\Omega)]^3$. Hence
$u\in D(T)$ and $Tu=-\Delta_eu$.

Conversely, if $u\in D(T)$ and $w=Tu$, then $u$ is a weak solution of
$-\Delta_eu=w$ with homogeneous Dirichlet boundary conditions. The principal symbol of $-\Delta_e$ is
$$
A(x,\xi)\eta
=
\mu(x)|\xi|^2\eta
+\bigl(\lambda(x)+\mu(x)\bigr)(\xi\cdot\eta)\xi,
$$
and {\rm(H2)} gives
$$
\eta\cdot A(x,\xi)\eta
=
\mu(x)|\xi|^2|\eta|^2
+\bigl(\lambda(x)+\mu(x)\bigr)(\xi\cdot\eta)^2
\geqslant
\underline\mu|\xi|^2|\eta|^2,
\,\,\,
\xi,\eta\in\mathbb R^3.
$$
Thus, $-\Delta_e$ satisfies the uniform Legendre-Hadamard condition. Its principal coefficients are Lipschitz continuous, its lower-order coefficients are bounded and $\partial\Omega$ is smooth. The $[W^{2,2}(\Omega)]^3$ regularity theorem for strongly elliptic Dirichlet systems \cite{ADN}, see also \cite[Thms. 4.16 and 4.18]{mclean}, gives
$u\in[H^2(\Omega)\cap H_0^1(\Omega)]^3$. Consequently,
$$
D(T)
=
[H^2(\Omega)\cap H_0^1(\Omega)]^3
\,\,\,\mbox{and}\,\,\,
T=-\Delta_e.
$$
\end{proof}

\begin{Remark}[Constant shear modulus]\label{rem:muconst}
	If $\mu\equiv\mu_0$ is constant, the positivity of the bulk modulus is no longer required. Indeed, \eqref{eq:korn-id} gives
	$$
	a(u,u)
	=
	\int_\Omega
	(\lambda+\mu_0)(\dive u)^2dx
	+\mu_0\|\nabla u\|_2^2.
	$$
	Consequently, if $\lambda+\mu_0\geqslant0$, then
	$$
	a(u,u)
	\geqslant
	\mu_0\|\nabla u\|_2^2,
	\,\,\, \forall\,u\in[H_0^1(\Omega)]^3.
	$$
	Thus, the strict inequality in {\rm(H2)} may be relaxed to
	$\lambda+\mu_0\geqslant0$. Item $1$ of Proposition \ref{prop:unperturbed} then holds with coercivity constant $\mu_0$, and item $2$ remains valid. Moreover,
	$\nabla\mu=0$ and $\nabla^2\mu=0$ imply
	$\mu_1=\mu_2=0$, so \eqref{eq:smallness} reduces to
	$2M<\mu_0\lambda_1$. The lower-order corrections introduced in Subsection \ref{sub:modified} also vanish and the modified operator coincides with $-\Delta_e$.
	
	More generally, suppose that
	$\nabla\mu=0$ and $\lambda+\mu=0$ on an open subset of $\Omega$. Then
	$\nabla\lambda=0$ there and \eqref{eq:lame-expanded} reduces to
	$$
	\Delta_eu=\mu\Delta u.
	$$
	Thus, the elastic operator locally becomes a multiple of the vector Laplacian and the two characteristic speeds coincide. This local configuration, together with the relaxation
	$\lambda+\mu\geqslant0$, appears in Section \ref{sec:robust}, although $\mu$ is not assumed to be constant throughout $\Omega$.
\end{Remark}

\begin{Remark}[Comparison of the original and modified forms]\label{rem:korn}
The symmetric form $a$ in \eqref{eq:physform} has the pointwise quadratic form
$$
Q_0(x,\xi) =\lambda(x)\left[\sum_{i=1}^3\xi_{ii}\right]^2
+\frac{\mu(x)}{2}\sum_{i,j=1}^3\bigl(\xi_{ij}+\xi_{ji}\bigr)^2.
$$
This quadratic form does not satisfy the Legendre condition. Indeed, if $\xi$ is a nonzero antisymmetric matrix, then $\sum_{i=1}^3\xi_{ii}=0$ and	$\xi_{ij}+\xi_{ji}=0$, so that  $Q_0(x,\xi)=0$. This pointwise degeneracy does not prevent the form $a$ from being coercive on 	$[H_0^1(\Omega)]^3$, since the infinitesimal rotations not detected by $Q_0$ are controlled globally by \eqref{eq:korn-id}. It also does not affect the Legendre-Hadamard condition established in the proof of Proposition \ref{prop:unperturbed}.
	
Two sufficient conditions yield coercivity of $a$ directly. The first is $3\lambda+2\mu>0$, treated in Proposition \ref{prop:unperturbed}, with coercivity constant $C_{K}$. The second is $\lambda>-\underline\mu$ on $\overline\Omega$. Indeed, \eqref{eq:korn-id} gives
$$
2\int_\Omega\mu|\varepsilon(u)|^2dx\geqslant
	\underline\mu
	\left[
	\|\nabla u\|_2^2+\|\dive u\|_2^2
	\right],
	$$
	and therefore
$$
	a(u,u)
	\geqslant
	\underline\mu\|\nabla u\|_2^2
	+\int_\Omega(\lambda+\underline\mu)(\dive u)^2dx
	\geqslant
	\underline\mu\|\nabla u\|_2^2.
$$
Both conditions are stronger than {\rm(H2)} when $\mu$ is genuinely variable and both fail for the coefficients in Example \ref{rem:example}.
	
The next subsection introduces the modified operator $\Le$, which differs from $-\Delta_e$ by lower-order terms involving $\nabla\mu$ and $\nabla^2\mu$. Its associated form is coercive under {\rm(H2)} with constant $\underline\mu$ and its principal part satisfies the full Legendre condition. This allows the regularity theory of Ne\v{c}as \cite[Ch. 5, \S5.3.3]{necas} to be applied directly. Condition \eqref{eq:smallness} is subsequently used to control the additional lower-order terms in the dissipativity, quasi-stability and smoothing estimates. Thus, the modified operator is introduced to recover coercivity under {\rm(H2)}, rather than to restore symmetry or elliptic regularity.
\end{Remark}

\subsection{The modified Lam\'e operator: the case $\lambda+\mu>0$}\label{sub:modified}
We now turn to the regime actually covered by {\rm(H2)}, in which, by Remark \ref{rem:korn},
	neither $3\lambda+2\mu>0$ nor $\lambda>-\underline\mu$ need hold and the form $a$ of
	\eqref{eq:physform} need not be coercive on $[H_0^1(\Omega)]^3$.
The analysis developed in this paper relies on a suitable modification
of the Lam\'e operator. The starting point for this construction is the
following symmetric bilinear form on $[H_0^1(\Omega)]^3$
\begin{align}\label{eq:form}
	\left\langle u, v\right\rangle_{e}
	\!=\! \int_{\Omega}\Bigl[\!\left(\lambda \!+\! \mu\right)\dive(u)\dive(v) \!+\!  \mu \nabla u \!\cdot\! \nabla v
	\!+\! \dive(v) u \!\cdot\! \nabla \mu \!+\! \dive(u) v\!\cdot\! \nabla \mu  \!+\! \mu_1  u \!\cdot\! v\Bigr]dx,
\end{align}
where $\nabla u\cdot\nabla v=\sum_{i,j=1}^3\partial_ju_i\,\partial_jv_i$.
\begin{Proposition}\label{prop:coercive}
	The bilinear form $\langle\cdot,\cdot\rangle_e$ defines an inner product on $[H_0^1(\Omega)]^3$ and the induced norm $\|\cdot\|_e$ is equivalent to the standard $[H_0^1(\Omega)]^3$ norm. More
	precisely,
	\begin{equation}\label{eq:coerc}
		\underline\mu\,\|\nabla u\|_2^2\ \leqslant\ \|u\|_e^2 \leqslant
		\Big[2\max_{\overline\Omega}(\lambda+\mu)+\overline\mu+\tfrac{2\mu_1}{\lambda_1}\Big]\,\|\nabla u\|_2^2,
		\,\,\,\forall\,u\in[H_0^1(\Omega)]^3 .
	\end{equation}
\end{Proposition}
\begin{proof}
	Bilinearity and symmetry are immediate from \eqref{eq:form}. By Young  inequality, for a.e.
	$x\in\Omega$,
	\begin{equation}\label{eq:young-form}
		\begin{aligned}
			2\big|\dive(u)\,u\cdot\nabla\mu\big|
			&\leqslant (\lambda+\mu)\, \dive^2(u)+ \frac{\big|u\cdot\nabla\mu\big|^2}{\lambda+\mu}\\
			&\leqslant (\lambda+\mu)\, \dive^2(u)+ \frac{|\nabla\mu|^2}{\lambda+\mu}\,|u|^2\\
			& \leqslant\ (\lambda+\mu)\, \dive^2(u)+ \mu_1|u|^2.
		\end{aligned}
	\end{equation}
	Therefore,
	\begin{align*}
		\|u\|_{e}^{2}
		&=\int_{\Omega}\Bigl[
		(\lambda+\mu)\dive^2(u)+\mu|\nabla u|^2
		+2\dive(u)\,u\cdot\nabla\mu+\mu_1|u|^2
		\Bigr]dx\\
		&\geqslant
		\int_{\Omega}\mu|\nabla u|^2dx
		+\int_{\Omega}
		\left(\mu_1-\frac{|\nabla\mu|^2}{\lambda+\mu}\right)|u|^2dx
		\geqslant
		\underline{\mu}\|\nabla u\|_2^2.
	\end{align*}
	Here the mixed term is absorbed by the term involving
	$(\lambda+\mu)\dive^2(u)$ through \eqref{eq:young-form}, while the
	remaining zero-order term is nonnegative by the definition of
	$\mu_1$. In particular, $\|u\|_e=0$ implies $u=0$ and hence
	$\langle\cdot,\cdot\rangle_e$ is an inner product.

	Applying \eqref{eq:young-form} with the opposite sign and then using
	\eqref{eq:div-grad} and Poincar\'e inequality, we obtain
	\begin{align*}
		\|u\|_{e}^{2}&\leqslant
		2\int_{\Omega}(\lambda+\mu)\dive^2(u)\,dx +\overline{\mu}\|\nabla u\|_2^2+2\mu_1\|u\|_2^2\\
		&\leqslant
		\left[
		2\max_{\overline\Omega}(\lambda+\mu)
		+\overline{\mu}
		+\frac{2\mu_1}{\lambda_1}
		\right]\|\nabla u\|_2^2.
	\end{align*}

	Thus $\|\cdot\|_e$ is equivalent to the standard norm of $[H_0^1(\Omega)]^3$.
\end{proof}
\bigskip
To identify the differential operator associated with \eqref{eq:form}, consider the \underline{multiplication operator} $\HH:[L^2(\Omega)]^3\to[L^2(\Omega)]^3$ defined by
$$\HH u=\Hess(\mu)\,u,
$$
where $\Hess(\mu)$ denotes the Hessian matrix of $\mu$. Since
$\mu\in C^2(\overline{\Omega})$, the operator $\HH$ is bounded. Moreover,
the symmetry of $\Hess(\mu)$ implies that $\HH$ is self-adjoint and satisfies
\begin{equation}\label{eq:H-sym}
	\langle \HH u,v\rangle
	=\int_\Omega v^\top\Hess(\mu)\,u\,dx
	=\langle u,\HH v\rangle\,\,\,\mbox{and}\,\,\,
	|\langle\HH u,u\rangle|\leqslant 6\mu_2\|u\|_2^2.
\end{equation}
The relation between the bilinear form \eqref{eq:form} and the original
Lam\'e operator is established by the following integration-by-parts
identity.
\begin{Lemma}[Integration by parts]\label{lem:ibp}
	For every $u\in[H^2(\Omega)\cap H_0^1(\Omega)]^3$, the identity
	\begin{equation}\label{eq:ibp}
		-\int_\Omega \Delta_e(u)\cdot v\,dx
		=\langle u,v\rangle_e+\langle\HH u,v\rangle-\mu_1\langle u,v\rangle
	\end{equation}
	holds for every $v\in[H_0^1(\Omega)]^3$.
\end{Lemma}
\begin{proof}
	Fix $u\in[H^2(\Omega)\cap H_0^1(\Omega)]^3$. Since $[C_c^\infty(\Omega)]^3$ is dense in
	$[H_0^1(\Omega)]^3$, it is enough to establish \eqref{eq:ibp} for $v\in[C_c^\infty(\Omega)]^3$. Indeed, for fixed $u$, both sides of \eqref{eq:ibp} define continuous linear functionals of $v$ on $[H_0^1(\Omega)]^3$. Since
	$(\Delta_eu)_i=\sum_j\partial_j(\sigma_e(u))_{ij}$, an integration by parts yields
	\begin{align*}
		\int_\Omega \Delta_e(u)\cdot v\,dx
		&=-\int_\Omega\lambda\,\dive(u)\dive(v)\,dx-\int_\Omega\mu\,\nabla u\cdot\nabla v\,dx
		-\underbrace{\sum_{i,j}\int_\Omega \mu\,\partial_iu_j\,\partial_jv_i\,dx}_{=: {\rm I}} .
	\end{align*}
	The term $I$ can be rewritten by integrating by parts with respect to
	$x_j$. The boundary term vanishes since $v$ has compact support and using $\partial_j(\mu\,\partial_iu_j)=(\partial_j\mu)\partial_iu_j+\mu\,\partial_i(\dive u)$, we obtain
	$$
	{\rm I}=\underbrace{-\sum_{i,j}\int_\Omega(\partial_j\mu)(\partial_iu_j)\,v_i\,dx}_{=:{\rm I}_1}
	\underbrace{-\sum_i\int_\Omega\mu\,\partial_i(\dive u)\,v_i\,dx}_{=:{\rm I}_2} .
	$$
	A further integration by parts with respect to $x_i$ gives
	$$
	{\rm I}_1
	=\sum_{i,j}\int_\Omega u_j\,\partial_i\big[(\partial_j\mu)v_i\big]dx
	=\langle \HH(u),v\rangle+\int_\Omega \dive(v)\,u\cdot\nabla\mu\,dx
	$$
	and
	$$
	{\rm I}_2=\sum_i\int_\Omega \dive(u)\,\partial_i(\mu v_i)\,dx
	=\int_\Omega\mu\,\dive(u)\dive(v)\,dx+\int_\Omega \dive(u)\,v\cdot\nabla\mu\,dx .
	$$
	Combining these identities gives
	\begin{align*}
		-\int_\Omega\Delta_e(u)\cdot v\,dx=
		&
		\int_\Omega\Big[(\lambda\!+\! \mu)\dive(u)\dive(v)+\mu\nabla u\cdot\nabla v\Big]\,dx\\
		&+	\int_\Omega	\Big[	\dive(u)\,v\cdot\nabla\mu \,+\, \dive(v)\,u\cdot\nabla\mu	\Big]\,dx
		+\langle\HH u,v\rangle.
	\end{align*}
	By the definition of $\langle\cdot,\cdot\rangle_e$, the right-hand side
	is equal to $\langle u,v\rangle_e
	+\langle\HH u,v\rangle
	-\mu_1\langle u,v\rangle.$
	This proves \eqref{eq:ibp}.
\end{proof}
We are now in a position to define the modified Lam\'e operator. We do it through its form, which is the order in which self-adjointness is actually available.
\begin{Definition}\label{def:Le}
	The \underline{modified Lam\'e operator} $\Le$ associated with the form $\langle\cdot,\cdot\rangle_e$ is defined by
	$$
	D(\Le)=
	\left\{
	u\in[H_0^1(\Omega)]^3 \left|
	\begin{array}{l}
		\text{there exists } w\in[L^2(\Omega)]^3 \text{ such that}\\[1mm]
		\langle u,v\rangle_e=\langle w,v\rangle
		\text{ for every }v\in[H_0^1(\Omega)]^3
	\end{array} \right.
	\right\}.
	$$
	For $u\in D(\Le)$, we set $\Le u=w$. Equivalently,
	\begin{equation}\label{weak-Le}
		\langle\Le u,v\rangle= \langle u,v\rangle_e,\,\,\,\forall\, u\in D(\Le),\, \forall \,		v\in[H_0^1(\Omega)]^3.
	\end{equation}
\end{Definition}
\begin{Proposition}[Self-adjointness and positivity] \label{prop:selfadj} Assume {\rm(H2)}. Then $\Le$ is a strictly positive self-adjoint operator on $[L^2(\Omega)]^3$ and
	$$\langle\Le u,u\rangle=\|u\|_e^2\geqslant\underline{\mu}\lambda_1\|u\|_2^2,\,\,\, u\in D(\Le). $$
	Moreover, $\Le^{-1}$ is compact
	\begin{equation}\label{eq:half-power}
		D\bigl(\Le^{1/2}\bigr)=[H_0^1(\Omega)]^3
		\,\,\,\mbox{and}\,\,\,\,
		\bigl\|\Le^{1/2}u\bigr\|_2=\|u\|_e,\,\,\, u\in[H_0^1(\Omega)]^3 .
	\end{equation}
\end{Proposition}
\begin{proof}
	By Proposition \ref{prop:coercive}, the symmetric form $\langle\cdot,\cdot\rangle_e$ is bounded and coercive on	$[H_0^1(\Omega)]^3$, which is dense in $[L^2(\Omega)]^3$. Since its 	form norm is equivalent to the standard norm of $[H_0^1(\Omega)]^3$, the form is closed. The first representation theorem \cite[Ch. VI, Thm. 2.1]{kato-book} therefore associates with it an m-sectorial operator, which is precisely the operator $\Le$ of Definition \ref{def:Le}. Its symmetry implies that $\Le$ is self-adjoint \cite[Ch.~VI, Thm.~2.6]{kato-book}.

	For $u\in D(\Le)$, \eqref{eq:coerc}, \eqref{weak-Le} and Poincar\'e	inequality give
	$$\langle\Le u,u\rangle
	=\|u\|_e^2
	\geqslant\underline{\mu}\|\nabla u\|_2^2
	\geqslant\underline{\mu}\lambda_1\|u\|_2^2.
	$$
	Thus $\Le$ is strictly positive and boundedly invertible. The second representation theorem \cite[Ch. VI, Thm. 2.23]{kato-book} yields \eqref{eq:half-power}. Finally, the Lax-Milgram theorem shows that $\Le^{-1}$ maps $[L^2(\Omega)]^3$ boundedly into $[H_0^1(\Omega)]^3$. Its compactness then follows from the compact embedding $[H_0^1(\Omega)]^3\subset[L^2(\Omega)]^3$.
\end{proof}
\begin{Proposition}[Strong ellipticity of $\Le$] \label{prop:elliptic}
	Assume {\rm (H2)}. The quadratic form associated with the principal
	part of $\Le$ is given by
	$$Q(x,\zeta)=\bigl[\lambda(x)+\mu(x)\bigr]\!
	\left[\sum_{i=1}^3\zeta_{ii}\right]^2\!\!
	+\mu(x)\sum_{i,j=1}^3\zeta_{ij}^2,\,\,\, \zeta\in\mathbb R^{3\times3}.
	$$
	Then $Q$ satisfies the uniform Legendre condition
	$$Q(x,\zeta) \geqslant \underline{\mu}|\zeta|^2, \,\,\, \forall\, x\in\overline{\Omega},\, \forall\, \zeta\in\mathbb R^{3\times3}.$$
	Consequently, the principal part of $\Le$ is uniformly strongly
	elliptic in the sense of Legendre-Hadamard. More precisely,
	$$ Q(x,\eta\otimes\xi)
	\geqslant
	\underline{\mu}\,|\eta|^2|\xi|^2,\,\,\, x\in\overline{\Omega},\,\eta,\xi\in\mathbb R^3.
	$$
\end{Proposition}
\begin{proof} Assumption {\rm(H2)} ensures that	$\lambda(x)+\mu(x)>0$ and $\mu(x)\geqslant\underline{\mu}>0$ for every $x\in\overline{\Omega}$. Therefore, for every $\zeta\in\mathbb R^{3\times3}$,
	$$
	Q(x,\zeta)
	\geqslant\underline{\mu}\sum_{i,j=1}^3\zeta_{ij}^2=\underline{\mu}|\zeta|^2.
	$$
	This proves the uniform Legendre condition. It remains to verify the strong ellipticity of $\Le$. The coefficients	of the principal part of the variational form associated with $\Le$ are given by $a_{jl}^{ik}(x)=\bigl(\lambda(x)+\mu(x)\bigr)\delta_{ij}\delta_{kl} +\mu(x)\delta_{ik}\delta_{jl}$. For $\xi,\eta\in\mathbb R^3$, set $\zeta=\eta\otimes\xi$ then $\zeta_{ij}=\eta_i\xi_j$, $\sum_{i=1}^3\zeta_{ii}=\eta\cdot\xi$ and $	|\zeta|^2=|\eta\otimes\xi|^2=|\eta|^2|\xi|^2$. Moreover, the Legendre estimate consequently gives
	\begin{align*}
		\sum_{i,j,k,l=1}^3a_{jl}^{ik}(x)\xi_j\xi_l\eta_i\eta_k
		=\bigl[\lambda(x)+\mu(x)\bigr](\eta\cdot\xi)^2+\mu(x)|\eta|^2|\xi|^2=Q(x,\eta\otimes\xi)\geqslant
		\underline{\mu}|\eta|^2|\xi|^2.
	\end{align*}
	Thus, the principal part of $\Le$ satisfies the uniform	Legendre-Hadamard condition. Since the terms $-\HH$ and $\mu_1I$ are of order zero, they do not affect the principal coefficients.
	Consequently, $\Le$ is uniformly strongly elliptic.
\end{proof}
\begin{Theorem}[Domain and differential representation of $\Le$]\label{thm:domain}
	Assume {\rm (H2)}. Then $D(\Le)	= \bigl[H^2(\Omega)\cap H_0^1(\Omega)\bigr]^3$ and $\Le$ admits the differential representation
	\begin{equation}\label{eq:Le-expression}
		\Le u=-\Delta_e u-\HH u+\mu_1u, \,\,\, u\in D(\Le).
	\end{equation}
	Moreover, the norms $\|\Le \,\cdot\|_2$ and $\|\cdot\|_{H^2}$ are equivalent on $D(\Le)$.
\end{Theorem}
\begin{proof} We first prove that $D(\Le)\subseteq [H^2(\Omega)\cap H_0^1(\Omega)]^3$. Let $u\in D(\Le)$ and set $w=\Le u\in[L^2(\Omega)]^3$. By the definition of the operator associated with the form, we have $\langle u,v\rangle_e= \langle w,v\rangle$ for all $v\in[H_0^1(\Omega)]^3$. Thus, $u$ is a weak solution of the Dirichlet problem associated 	with $\Le u=w$.	By Proposition \ref{prop:elliptic}, the principal part of $\Le$	satisfies the uniform Legendre condition and hence the uniform Legendre-Hadamard condition. Its principal coefficients are Lipschitz continuous on $\overline{\Omega}$, since	$\lambda\in C^1(\overline{\Omega})$ and $\mu\in C^2(\overline{\Omega})$ and the coefficients of its 	lower-order terms are bounded. Since $\partial\Omega$ is of class $C^\infty$, the $[W^{2,2}(\Omega)]^{3}$-regularity theorem for strongly elliptic Dirichlet systems \cite[Ch. 5, \S5.3.3, Lemma 3.3]{necas} (see also \cite{ADN} and \cite[Thms.4.16 and 4.18]{mclean})  yields $u\in[H^2(\Omega)\cap H_0^1(\Omega)]^3$ together with the estimate
	\begin{equation}\label{e1}
		\|u\|_{H^2} \leqslant C\|w\|_2 = C\|\Le u\|_2 .
	\end{equation}
	We now prove the reverse inclusion. Let	$u\in[H^2(\Omega)\cap H_0^1(\Omega)]^3$. By Lemma \ref{lem:ibp}, $\langle u,v\rangle_e=\langle-\Delta_eu-\HH u+\mu_1u,v\rangle$ for all $v\in[H_0^1(\Omega)]^3$. Moreover, the regularity of the coefficients implies that $	-\Delta_eu-\HH u+\mu_1u\in[L^2(\Omega)]^3$. It follows from the definition of the operator associated with the form that $u\in D(\Le)$ and $\Le u=-\Delta_eu-\HH u+\mu_1u$. Therefore, $D(\Le)=[H^2(\Omega)\cap H_0^1(\Omega)]^3$. Finally, the differential representation of $\Le$ and the boundedness of its coefficients give
	$$\|\Le u\|_2\leqslant C\|u\|_{H^2},\,\,u\in D(\Le).$$
	Together with \eqref{e1}, this proves the equivalence of the norms on $D(\Le)$.
\end{proof}
Since $\Le$ and $\Ae$ are positive self-adjoint operators, the Heinz
inequality yields the following comparison of their fractional powers.
\begin{Corollary}[Equivalence of fractional domains]\label{cor:fract}
	The operators $\Le$ and $\Ae$ have the same domain and there exist constants $m_1,m_2>0$ such that
	\begin{equation}\label{eq:operator-norm-equivalence}
		m_1\|\Ae u\|_2\leqslant \|\Le u\|_2 \leqslant m_2\|\Ae u\|_2,\,\,\,u\in D(\Ae)=D(\Le).
	\end{equation}
	Moreover, for every $r\in[0,1]$, their fractional domains coincide $D(\Le^r)=D(\Ae^r)$ and the corresponding fractional norms are equivalent. More precisely,
	\begin{equation}\label{eq:heinz}
		m_1^{\,r}\|\Ae^r u\|_2\leqslant	\|\Le^r u\|_2 \leqslant	m_2^{\,r}\|\Ae^r u\|_2,\,\,\,u\in D(\Ae^r)=D(\Le^r).
	\end{equation}
\end{Corollary}
\begin{proof}
	By Theorem \ref{thm:domain} and the corresponding classical result for  $\Ae$, we have $D(\Le)=D(\Ae)$.  Moreover, both $\|\Le\cdot\|_2$ and $\|\Ae\cdot\|_2$ are equivalent to the $[H^2(\Omega)]^{3}$-norm on this common domain, which proves \eqref{eq:operator-norm-equivalence}.

	We now consider the fractional powers. Both operators are positive and self-adjoint, by
	Proposition \ref{prop:selfadj} and by the classical properties of the Dirichlet Laplacian,
	respectively. Recall the Heinz inequality \cite{heinz} (see also \cite{kato61}). If $S$ and $T$ are positive self-adjoint operators on a Hilbert space satisfying
	$$ D(T)\subseteq D(S)\,\,\,\text{and}\,\,\, \|Su\|\leqslant\|Tu\|,\,\,\,u\in D(T)$$
	then, for every $r\in[0,1]$,
	$$ D(T^r)\subseteq D(S^r)\,\,\,\text{and}\,\,\, \|S^ru\|\leqslant\|T^ru\|,\,\,\, u\in D(T^r).$$
	From the upper estimate in \eqref{eq:operator-norm-equivalence}, we	may apply this result with
	$S=\Le$ and $T=m_2\Ae$. Since $(m_2\Ae)^r=m_2^r\Ae^r$, we obtain $D(\Ae^r)\subseteq D(\Le^r)$ 	and $\|\Le^ru\|_2\leqslant m_2^r\|\Ae^ru\|_2$ for all $u\in D(\Ae^r)$. Similarly, the lower estimate in \eqref{eq:operator-norm-equivalence} allows us to apply the Heinz inequality with $	S=m_1\Ae$ and $T=\Le$. This proves \eqref{eq:heinz}.
\end{proof}
\begin{Remark}\label{rem:fract}
	As a consequence of Corollary \ref{cor:fract}, the interpolation inequality \eqref{eq:interp} and the Sobolev embeddings \eqref{eq:sobolev} remain valid, possibly with different constants, when $\Ae$ is replaced by $\Le$. The elementary comparison of fractional powers, on the other hand, does not need Corollary \ref{cor:fract} and follows directly from the spectral theorem applied to $\Le$: writing $u=\sum_k u_k e_k$ in an orthonormal basis of eigenvectors of $\Le$, with eigenvalues $\nu_k\geqslant\nu_1=\min\sigma(\Le)\geqslant\underline\mu\lambda_1$ by Proposition \ref{prop:selfadj}, one has, for $0\leqslant r<s$ and for every $u\in D(\Le^{s/2})$
	\begin{align}\label{eq:power-Le}
		\|\Le^{r/2}u\|_2^2=\sum_k\nu_k^{\,r-s}\nu_k^{\,s}|u_k|^2\leqslant\nu_1^{\,r-s}\|\Le^{s/2}u\|_2^2\,\mbox{, that is,}\,\,\,
		\|\Le^{r/2}u\|_2 \leqslant \bigl[\underline\mu\lambda_1\bigr]^{(r-s)/2} \|\Le^{s/2}u\|_2.
	\end{align}
\end{Remark}
\begin{Remark}[A fractional embedding]\label{rem:key} The following estimate will be used repeatedly in Section \ref{sec:exp}. Let $r,s\geqslant0$ satisfy $r+s\leqslant1$ and set $ \sigma=1-r-s\in[0,1]$. Since $\Le$ is positive and self-adjoint, the composition rule for its fractional powers gives $\Le^{\sigma/2}\Le^{r+s}u=\Le^{(1+r+s)/2}u$. In particular, $u\in D\bigl(\Le^{(1+r+s)/2}\bigr)$ if and only if $\Le^{r+s}u\in D\bigl(\Le^{\sigma/2}\bigr)$. By Remark \ref{rem:fract}, the Sobolev embedding \eqref{eq:sobolev} remains valid for $\Le$. Hence, $D(\Le^{\sigma/2})\subset[L^{6/(3-2\sigma)}(\Omega)]^3=[L^{6/(1+2r+2s)}(\Omega)]^3$. Applying this embedding to $\Le^{r+s}u$, we obtain
	\begin{equation}\label{eq:key}
		\bigl\|\Le^{r+s}u\bigr\|_{\frac{6}{1+2r+2s}}
		\leqslant C\bigl\|\Le^{(1+r+s)/2}u\bigr\|_2,\,\,\,u\in D\bigl(\Le^{(1+r+s)/2}\bigr).
	\end{equation}
	This argument avoids an explicit identification of $D(\Le^\theta)$ with a Sobolev space of order $2\theta$ when $2\theta>1$, where the characterization of the fractional domain requires additional boundary conditions. Such an identification is nevertheless legitimate whenever $\theta\leqslant1$, since Corollary \ref{cor:fract} then applies and reduces it to \eqref{eq:charact}; this is used once, in Case $2$ of the proof of Lemma \ref{lem:bootstrap}.
\end{Remark}
\bigskip
Finally, we introduce on $[L^2(\Omega)]^3$ the weighted inner product
\begin{equation}\label{eq:rho-ip}
	\langle u,v\rangle_\rho=\int_\Omega \rho(x)\,u(x)\cdot v(x)\,dx\,\,\,\mbox{and}\,\,\,
	\|u\|_\rho^2=\int_\Omega\rho\,|u|^2dx,
\end{equation}
so that $\tfrac12\|u_t\|_\rho^2$ is the kinetic energy of the heterogeneous body. Since
$\underline\rho\leqslant\rho(x)\leqslant\overline\rho$ in $\Omega$, the norms $\|\cdot\|_\rho$ and $\|\cdot\|_2$
are equivalent. More precisely,
\begin{align}\label{eq:rho-equiv}
	\underline\rho \|u\|_2^2  \leqslant  \|u\|_\rho^2  \leqslant\ \overline\rho\,\|u\|_2^2.
\end{align}
\begin{Remark}[How much heterogeneity the framework tolerates]\label{rem:window}
		The constant $\mu_1$ of {\rm(H2)} is a maximum, so the smallness condition
			\eqref{eq:smallness} can be inverted: it is \emph{equivalent} to the pointwise requirement
		\begin{equation}\label{eq:window-sharp}
			\bigl|\nabla\mu(x)\bigr|^{2}\ <\ \bigl(\lambda(x)+\mu(x)\bigr)\,
			\Bigl[\underline\mu\lambda_1-2M-6\mu_2\Bigr]
			\ =\ \rho(x)\bigl(c_P^2(x)-c_S^2(x)\bigr)\Bigl[\underline\mu\lambda_1-2M-6\mu_2\Bigr],
		\end{equation}
		for every $x\in\overline\Omega$ and it is \emph{implied by} the cruder global requirement
		\begin{equation}\label{eq:window-global}
			\max_{\overline\Omega}\bigl|\nabla\mu\bigr|^{2}\ <\ \epsilon_0\Bigl[\underline\mu\lambda_1-2M-6\mu_2\Bigr],
			\qquad\mbox{that is,}\qquad \bigl|\nabla\mu\bigr|=O\bigl(\sqrt{\epsilon_0}\bigr),
		\end{equation}
		since $\epsilon_0=\min_{\overline\Omega}(\lambda+\mu)$. Written in this way the restriction
			acquires a meaning: by \eqref{eq:speeds} it bounds the variation of the shear modulus at a point by
			the local gap between the two characteristic speeds, so that \emph{the shear modulus must become
				locally homogeneous at least as fast as that gap closes} pointwise in
			\eqref{eq:window-sharp}, uniformly in \eqref{eq:window-global}. Which of the two is imposed is
			immaterial as long as the gap stays away from zero and decisive as soon as it is allowed to
			degenerate; that distinction is taken up in Subsection \ref{sec:physics}, where it is also
			quantified.
\end{Remark}
\begin{Remark}[The unperturbed case revisited]\label{rem:dictionary}
		Suppose that, besides {\rm(H2)}, the medium satisfies $3\lambda+2\mu>0$ on $\overline\Omega$. By
		Proposition \ref{prop:unperturbed} the elastic form $a$ of \eqref{eq:physform} is then coercive by
		itself and $-\Delta_e$ is the strictly positive self-adjoint operator associated with it, with
		$D(-\Delta_e)=[H^2(\Omega)\cap H_0^1(\Omega)]^3$ and $D\bigl((-\Delta_e)^{1/2}\bigr)=[H_0^1(\Omega)]^3$. Exactly the properties of $\Le$ established in Proposition \ref{prop:selfadj} and
		Theorem \ref{thm:domain}. Consequently every result of Sections
		{\rm\ref{sec:wp}-\ref{sec:exp}} remains valid after replacing $\Le$ by $-\Delta_e$,
		$\langle\cdot,\cdot\rangle_e$ by $a$ and the coercivity constant $\underline\mu$ of
		\eqref{eq:coerc} by $C_K$ and setting $\HH=0$ and $\mu_1=\mu_2=0$. Indeed, $\HH$ and $\mu_1I$ enter
		the analysis only through \eqref{eq:form}, \eqref{eq:main-abs}, \eqref{eq:energy},
		\eqref{eq:E-lin}, the constant $C_4$ of Lemma \ref{lem:diff-est} and \eqref{eq:Ers} and all the
		corresponding terms then vanish. In particular the smallness condition \eqref{eq:smallness} reduces
		to $2M<C_K\lambda_1$. The same dictionary applies, with $C_K$ replaced by $\mu_0$ and under
		{\rm(H2)} alone, when the shear modulus is constant, in which case $\Le=-\Delta_e$
		(Remark \ref{rem:muconst}).
\end{Remark}
\section{Well-posedness}\label{sec:wp}
In this section we establish the well-posedness of problem \eqref{eq:main}. The functional framework
introduced in the previous section allows us to reformulate the system as an abstract evolution
equation in a suitable Hilbert space.
\subsection{Abstract formulation}
The characterization of $\Le$ established in Theorem \ref{thm:domain} allows us to rewrite \eqref{eq:main} as
\begin{equation} \label{eq:main-abs}
\begin{cases}
\rho\, u_{tt}+\Le u+ \alpha  u_t + f(u) = \mu_1u -\HH u+ b(x) &\,\,\text{in}\,\, \Omega\times\mathbb{R}^{+}, \\
u=0 &\,\,\text{on}\,\,\partial \Omega\times\mathbb{R}^{+}, \\
u(x,0)=u_0(x)\,\,\,\mbox{and}\,\,\, u_t(x,0)=u_1(x)&\,\,x \in \Omega.
\end{cases}
\end{equation}
The natural energy and strong phase spaces associated with
\eqref{eq:main-abs} are
$$\mathscr{H}^0= \big(D(\Le^{1/2}), \|\Le^{1/2} \cdot\|_2 \big)\times \big([L^2(\Omega)]^3, \| \cdot\|_\rho \big)$$
and
$$\mathscr{H}^1=  \big(D(\Le), \|\Le \cdot\|_2 \big)\times \big(D(\Le^{1/2}), \| \Le^{1/2} \cdot \|_2 \big).$$
Introducing the state variable $\U=(u,v)$ with $v=u_t$, system \eqref{eq:main-abs} becomes the
abstract Cauchy problem
\begin{equation}\label{eq:cauchy}
\begin{cases}
\U_t=\Ee\U+\Ff(\U),\\
\U(0)=\U_0=(u_0,u_1).
\end{cases}
\end{equation}
Here $\Ee:D(\Ee)\subset\mathscr H^0\to\mathscr H^0$ and $\Ff:\mathscr H^0\to\mathscr H^0$ are
\begin{equation} \label{eq:EF}
\Ee \,\U=\begin{pmatrix} v \\ -\rho^{-1}\big[\Le(u) +\alpha v\big]\end{pmatrix},
\, D(\Ee)=\mathscr{H}^1
\,\,\,\mbox{and}\,\,\,
\Ff(\U)= \begin{pmatrix} 0  \\ \rho^{-1}\big[\mu_1u \!-\! \HH(u) \!+\! b\!-\! f(u)\big] \end{pmatrix}.
\end{equation}
\subsection{Properties of the operators $\Ee$ and $\Ff$}
\begin{Lemma}\label{lem:generator}
The operator $\Ee$ generates a $C_0$-semigroup of contractions $\{S_h(t)\}_{t\geqslant0}$ on
$\mathscr{H}^0$.
\end{Lemma}
\begin{proof}
We first observe that $\Ee$ is dissipative. Indeed, for every $\U=(u,v)\in D(\Ee)$, using
\eqref{eq:rho-ip} and \eqref{weak-Le},
\begin{align*}
\langle \Ee\U,\U\rangle_{\mathscr H^0}
&=\big\langle\Le^{1/2}v,\Le^{1/2}u\big\rangle
-\big\langle \rho^{-1}\big[\Le u+\alpha v\big],\,v\big\rangle_\rho\\
&=\ \langle\Le u,v\rangle-\langle\Le u,v\rangle-\alpha\|v\|_2^2\\
&= -\alpha\|v\|_2^2\ \leqslant\ 0,
\end{align*}
the cancellation being a consequence of the choice of the weight in \eqref{eq:rho-ip}.
To apply the Lumer-Phillips theorem, it remains to prove that $I-\Ee$ is surjective. Let
$(\tilde u,\tilde v)\in\mathscr H^0$ be arbitrary. Solving $(I-\Ee)\U=(\tilde u,\tilde v)$ is
equivalent to finding $(u,v)\in D(\Ee)$ with
$$
v=u-\tilde u\,\,\,\text{and}\,\,\, (\alpha+\rho)\,u+\Le u=\rho\,\tilde v+(\alpha+\rho)\,\tilde u.
$$
Indeed, the first component of $(I-\Ee)\U=(\tilde u,\tilde v)$ reads $u-v=\tilde u$ and
multiplying the second one by $\rho$ gives $(\rho+\alpha)v+\Le u=\rho\tilde v$. Substituting
$v=u-\tilde u$ yields the stated equation. The latter admits the variational formulation
\begin{equation}\label{eq:varform}
\big\langle(\alpha+\rho)u,\phi\big\rangle +\langle u,\phi\rangle_e
= \big\langle \rho\,\tilde v +(\alpha+\rho)\tilde u,\phi\big\rangle,\,\,\, \forall\,\phi\in[H_0^1(\Omega)]^3 .
\end{equation}
Since $\alpha+\rho\geqslant\alpha+\underline\rho>0$, the bilinear form on the left-hand side of
\eqref{eq:varform} is continuous and coercive on $[H_0^1(\Omega)]^3$ by
Proposition \ref{prop:coercive}, so the Lax-Milgram theorem guarantees the existence of a unique
solution $u\in D(\Le^{1/2})=[H_0^1(\Omega)]^3$. Moreover, rewriting \eqref{eq:varform} as
$\langle u,\phi\rangle_e=\langle w,\phi\rangle$ for every $\phi\in[H_0^1(\Omega)]^3$, with
$w:=\rho\tilde v+(\alpha+\rho)(\tilde u-u)\in[L^2(\Omega)]^3$, Definition \ref{def:Le} yields
$u\in D(\Le)$ at once; no elliptic regularity is needed here. Setting
$v=u-\tilde u\in D(\Le^{1/2})$ we obtain $\U=(u,v)\in D(\Ee)$ satisfying
$(I-\Ee)\U=(\tilde u,\tilde v)$. Hence $I-\Ee$ is surjective, $\Ee$ is maximal dissipative and  by the Lumer-Phillips theorem  it generates a $C_0$-semigroup of contractions on $\mathscr H^0$.
\end{proof}
\begin{Lemma}[Local Lipschitz continuity of $\Ff$] \label{lem:lipschitz} For every bounded subset $D\subset\mathscr H^0$, there exists a constant $C_D>0$, such that
\begin{align}\label{lip1}
\big\|\Ff(\U^{1})-\Ff(\U^{2})\big\|_{\mathscr{H}^0}
\leqslant C_D\|\U^{1}-\U^{2}\|_{\mathscr{H}^0},\,\,\, \forall\,\U^{1},\U^{2}\in D.
\end{align}
Moreover, there exists $C>0$ such that
\begin{align}\label{estF}\|\Ff(\U)\|_{\mathscr H^0}\leqslant C\big[1+\|\U\|_{\mathscr H^0}\big]^3,\,\,\,\forall \,\U\in\mathscr H^0.
\end{align}
\end{Lemma}
\begin{proof}
Let $\U^i=(u^i,v^i)\in\mathscr H^0$ $(i=1,2)$ and set $\omega=u^1-u^2$. Since the first component of $\Ff$ vanishes, we have $\|\Ff(\U^1)-\Ff(\U^2)\|_{\mathscr H^0}=
\left\|\rho^{-1}\bigl[\mu_1\omega-\HH\omega	-\bigl(f(u^1)-f(u^2)\bigr)	\bigr]		\right\|_\rho$. By \eqref{eq:rho-equiv}, $\|\rho^{-1}w\|_\rho \leqslant	\underline{\rho}^{-1/2}\|w\|_2$. It is therefore enough to estimate the terms on the right-hand side  in $[L^2(\Omega)]^3$.
For the linear terms, the boundedness of $\HH$ on $[L^2(\Omega)]^3$ gives
\begin{equation}\label{linear-F}
\|\mu_1\omega-\HH\omega\|_2\leqslant C\|\omega\|_2.
\end{equation}
We now estimate the critical part of $f$. By \eqref{eq:growth-h}, we find $|h_i'(s)|\leqslant C[1+|s|^2]$, $s\in\mathbb R$. Hence, the mean value theorem yields
 $$|h_i(a)-h_i(b)| \leqslant C\bigl(1+|a|^2+|b|^2\bigr)|a-b|,\,\,\, a,b\in\mathbb R.$$
Applying H\"older inequality and using $\frac12=\frac{2}{6}+\frac16$, we obtain
\begin{equation}\label{h-lipschitz}
\|h(u^1)-h(u^2)\|_2 \leqslant C\bigl[1+\|u^1\|_6^2+\|u^2\|_6^2\bigr]\|\omega\|_6.
\end{equation}
For the subcritical part, \eqref{eq:growth-g} and the mean value theorem give
$|g_i(a)-g_i(b)|\leqslant
C\bigl[1+|a|^{p-1}+|b|^{p-1}\bigr]|a-b|$.
Since $1\leqslant p<3$ and $\frac12=\frac{p-1}{6}+\frac{4-p}{6}$, H\"older's inequality implies
\begin{equation}\label{g-lipschitz}
\|g(u^1)-g(u^2)\|_2
\leqslant
C\bigl[1+\|u^1\|_6^{p-1}+\|u^2\|_6^{p-1}\bigr]\|\omega\|_{\frac{6}{4-p}}.
\end{equation}
Since $2\leqslant\frac{6}{4-p}<6$, the Sobolev embedding $[H_0^1(\Omega)]^3\subset[L^{6}(\Omega)]^3$ controls both $\|\omega\|_6$ and $\|\omega\|_{6/(4-p)}$. Moreover, by \eqref{eq:half-power} and \eqref{eq:coerc}, $\|\omega\|_{H_0^1}\leqslant C\|\U^1-\U^2\|_{\mathscr H^0}$.
If $\U^1,\U^2$ belong to a bounded subset $D\subset\mathscr H^0$, then $\|u^1\|_6$ and $\|u^2\|_6$ are bounded by a constant depending only on $D$. Combining \eqref{linear-F}, \eqref{h-lipschitz} and \eqref{g-lipschitz}  we conclude \eqref{lip1}. Finally, the same estimates, together with $b\in[L^2(\Omega)]^3$ and $f(0)=0$, give
$$
\|\Ff(\U)\|_{\mathscr H^0}
\leqslant
C\left[1+\|u\|_{H_0^1}+\|u\|_{H_0^1}^{p}+\|u\|_{H_0^1}^{3}\right] \\
\leqslant	C\bigl[1+\|\U\|_{\mathscr H^0}\bigr]^3,
$$	where we used $p<3$.
This completes the proof.
\end{proof}
\subsection{Local well-posedness}
\begin{Theorem}[Local well-posedness]\label{thm:local}
Assume that {\rm(H1)-(H5)} hold. Then, for every initial datum $\U_0\in\mathscr{H}^0$, there exist a
maximal time $T_{\max}\in(0,\infty]$ and a unique mild solution
$\U\in C([0,T_{\max});\mathscr{H}^0)$ of \eqref{eq:cauchy} satisfying
\[
\U(t) =S_h(t)\U_0 + \int_0^t S_h(t-s)\Ff(\U(s))\,ds,\,\,\, 0\leqslant t<T_{\max} .
\]
Furthermore, either $T_{\max}=+\infty$ or
$\limsup_{t\to T_{\max}^{-}}\|\U(t)\|_{\mathscr{H}^0}=+\infty$. Finally, if $\U_0\in\mathscr{H}^1$
then $\U\in C([0,T_{\max});\mathscr{H}^1)\cap C^1([0,T_{\max});\mathscr{H}^0)$ and $\U$ is a classical
solution of \eqref{eq:cauchy}.
\end{Theorem}
\begin{proof}
By Lemma \ref{lem:generator}, $\Ee$ generates the $C_0$-semigroup of contractions
$\{S_h(t)\}_{t\geqslant 0}$ on $\mathscr{H}^0$; by Lemma \ref{lem:lipschitz}, $\Ff$ is locally Lipschitz
continuous on $\mathscr{H}^0$. The conclusion follows from the classical theory of semilinear
evolution equations, see \cite[Ch.~6, Thms. 1.4 and 1.5]{pazy}.
\end{proof}
\begin{Lemma}[Local continuous dependence]\label{lem:locdep} Let $\U^1$ and $\U^2$ be the maximal mild solutions of \eqref{eq:cauchy} with initial data $\U_0^1,\U_0^2\in\mathscr H^0$ and let	$T<\min\{T_{\max}^1,T_{\max}^2\}$. If $\|\U^i(t)\|_{\mathscr H^0}\leqslant R$ $(i=1,2)$ for $t\in[0,T]$, then there exists $C_R>0$, depending only on $R$, such that
$$\|\U^1(t)-\U^2(t)\|_{\mathscr H^0}\leqslant e^{C_Rt}\|\U_0^1-\U_0^2\|_{\mathscr H^0},\,\,\, t\in[0,T].$$
\end{Lemma}
\begin{proof}
Let $C_R$ be the Lipschitz constant of $\Ff$ on the ball of radius $R$ in $\mathscr H^0$, given by Lemma \ref{lem:lipschitz}. Subtracting the variation-of-constants formulas for $\U^1$ and $\U^2$ and using the contractivity of $\{S_h(t)\}_{t\geqslant0}$, we obtain
$$\|\U^1(t)-\U^2(t)\|_{\mathscr H^0}
\leqslant
\|\U_0^1-\U_0^2\|_{\mathscr H^0}
+C_R\int_0^t
\|\U^1(s)-\U^2(s)\|_{\mathscr H^0}\,ds.
$$
Gronwall inequality completes the proof.
\end{proof}
\subsection{Energy identity and a priori estimates}
Let $\U(t)=(u(t),u_t(t))$ be the maximal local solution provided by Theorem \ref{thm:local}. The
total energy functional associated with \eqref{eq:main-abs} is
\begin{align}\label{eq:energy}
\mathscr{E}_{\U}(t)=\frac{1}{2}\|\U(t)\|_{\mathscr{H}^{0}}^2
\!+\! \int_{\Omega}G(u)\,dx
\!+\! \sum_{i=1}^{3}\int_{\Omega}\int_{0}^{u_i}\!\! h_i(s)\,ds dx
\!+\! \frac{1}{2}\left\langle \HH(u),u\right\rangle
\!-\! \left\langle b,u\right\rangle \!-\! \frac{\mu_1}{2}\|u\|_2^2,
\end{align}
where
$$\frac12\|\U\|^2_{\mathscr H^0}=\frac12\|u\|_e^2+\frac12\int_\Omega\rho|u_t|^2dx.$$
\begin{Proposition}\label{prop:energy}
Let $\U(t)=(u(t),u_t(t))\in C([0,T_{\max});\mathscr H^0)$ be the maximal local solution of
\eqref{eq:main-abs}. Then the following statements hold.
\begin{itemize}
\item[$1.$] The energy is non-increasing and satisfies
\begin{equation}\label{eq:energy-id}
\mathscr{E}_{\U}(t)+\alpha\int_s^t\|u_t(\tau)\|_2^2\, d\tau=\mathscr{E}_{\U}(s),
\,\,\, 0\leqslant s\leqslant t<T_{\max}.
\end{equation}
\item[$2.$] There exist positive constants $C_0$, $C_1$ and $C_2$ such that
\begin{equation}\label{eq:energy-equiv}
C_1\|\U(t)\|_{\mathscr H^0}^{2}-C_0\ \leqslant\ \mathscr{E}_{\U}(t)
 \leqslant C_2\bigl[1+\|\U(t)\|_{\mathscr H^0}^{4}\bigr],\,\,\, 0\leqslant t<T_{\max}.
\end{equation}
\end{itemize}
\end{Proposition}
\begin{proof}
We first consider initial data $\U_0\in\mathscr H^1$, so that the corresponding solution is
sufficiently regular. Multiplying \eqref{eq:main-abs} by $u_t$ and integrating over $\Omega$ and using that $\langle\HH(u),u_t\rangle=\tfrac12\tfrac{d}{dt}\langle\HH(u),u\rangle$ by the
symmetry \eqref{eq:H-sym} and that $b$ does not depend on $t$, we obtain
$\mathscr{E}_{\U}^{\prime}(t)+\alpha\|u_t(t)\|_2^2=0$. Integration over $(s,t)$ yields
\eqref{eq:energy-id}. For arbitrary $\U_0\in\mathscr H^0$ the result follows by a standard density
argument, using the local continuous dependence of Lemma \ref{lem:locdep} and the continuity
of every term of \eqref{eq:energy} on $\mathscr H^0$.
Next, we prove the lower and upper estimates in \eqref{eq:energy-equiv}. Denote the nonlinear and lower-order part of the energy by
\begin{align*}
I(u)= \int_\Omega G(u)\,dx + \sum_{i=1}^3 \int_\Omega\int_0^{u_i}h_i(s)\,ds\,dx
+ \frac{1}{2}\langle \HH(u),u\rangle - \langle b,u\rangle - \frac{\mu_1}{2}\|u\|_2^2 .
\end{align*}
By \eqref{eq:diss2},
$\int_\Omega G(u)+\sum_i\int_\Omega\int_0^{u_i}h_i  \geqslant-M\|u\|_2^2-m_f|\Omega|$ and by
\eqref{eq:H-sym}, $\tfrac12\langle\HH(u),u\rangle\geqslant-3\mu_2\|u\|_2^2$. Using Young and Poincar\'e inequalities together with \eqref{eq:coerc}, we find
\[
|\langle b,u\rangle|\leqslant\frac{1}{4\epsilon}\|b\|_2^2+\epsilon\|u\|_2^2
\leqslant\frac{1}{4\epsilon}\|b\|_2^2+\frac{\epsilon}{\lambda_1\underline\mu}\|u\|_e^2 .
\]
Since $\|u\|_2^2\leqslant \lambda_1^{-1}\|\nabla u\|_2^2\leqslant (\lambda_1\underline\mu)^{-1}\|u\|_e^2$, adding
the three estimates gives
\begin{align*}
I(u)\geqslant -\left[\frac{2M+6\mu_2+\mu_1}{2\lambda_1\underline\mu}
+\frac{\epsilon}{\lambda_1\underline\mu}\right]\|u\|_e^2 -\frac{1}{4\epsilon}\|b\|_2^2 -m_f|\Omega|.
\end{align*}
Therefore
\begin{align*}
\mathscr{E}_{\U}(t)\geqslant\left[\frac{1}{2}-\frac{2M+6\mu_2+\mu_1}{2\lambda_1\underline\mu}
-\frac{\epsilon}{\lambda_1\underline\mu}\right]\|u(t)\|_e^2+\frac{1}{2}\|u_t(t)\|_\rho^2
-m_f|\Omega|-\frac{1}{4\epsilon}\|b\|_2^2 .
\end{align*}
By \eqref{eq:smallness} we may choose $\epsilon>0$ so small that
$C_\epsilon=\tfrac{1}{2}-\tfrac{2M+6\mu_2+\mu_1}{2\lambda_1\underline\mu}-\tfrac{\epsilon}{\lambda_1\underline\mu}>0$.
Setting $C_1=\min\{C_\epsilon,\tfrac12\}>0$ and $C_0=m_f|\Omega|+\tfrac{1}{4\epsilon}\|b\|_2^2$ we obtain the lower bound.
From \eqref{eq:growth-g} and $g(0)=0$, we get
$|g(u)|\leqslant C\big[|u|+|u|^{p}\big]$. Hence, $|G(u)|\leqslant C[|u|^2+|u|^{p+1}]$. Using \eqref{eq:growth-h} together with $h_i(0)=0$ and $h_i\in C^2(\mathbb R)$, we find $|h_i(x)|\leqslant C\big[|x|+|x|^{3}\big]$ and $\big|\int_0^{x}h_i(s)\,ds\big|\leqslant C[x^2+x^4]$. Therefore
\begin{align*}
I(u)\leqslant C\left[1+\|u\|_{p+1}^{p+1}+\sum_{i=1}^3\|u_i\|_4^4\right] + \frac{1}{2}\|b\|_2^2 + C\|u\|_2^2 .
\end{align*}
Using the Sobolev embeddings $[H_0^1(\Omega)]^3\subset[L^{p+1}(\Omega)]^3$ (recall $p+1<4$) and
$[H_0^1(\Omega)]^3\subset[L^4(\Omega)]^3$, together with \eqref{eq:coerc}, we infer
$\mathscr{E}_{\U}(t) \leqslant C_2\big[1+\|\U(t)\|_{\mathscr H^0}^4\big]$.
\end{proof}
\subsection{Global well-posedness}
\begin{Theorem}[Global well-posedness]\label{thm:global}
Assume that {\rm(H1)-(H5)} hold. Then the following assertions are valid.
\begin{itemize}
\item[$1.$] For every $\U_0\in\mathscr{H}^0$, problem \eqref{eq:cauchy} admits a unique global mild
solution $\U\in C([0,\infty);\mathscr{H}^0)$, given by the variation-of-constants formula.
\item[$2.$] If $\U_0\in\mathscr{H}^1=D(\Ee)$, then the corresponding solution is a global strong
solution and $\U\in C([0,\infty);\mathscr{H}^1)\cap C^1([0,\infty);\mathscr{H}^0)$. In particular
$\U_t(t)=\Ee\,\U(t)+\Ff(\U(t))$ holds in $\mathscr H^0$ for every $t\geqslant0$.
\item[$3.$] For every bounded set $D\subset\mathscr H^0$, there exists
$C_D>0$, such that, any two solutions with initial data $\U_0^1,\U_0^2\in D$ satisfy
\begin{equation}\label{eq:cont-dep}
\|\U^1(t)-\U^2(t)\|_{\mathscr{H}^0}^2\leqslant e^{C_D t}\|\U_0^1-\U_0^2\|_{\mathscr{H}^0}^2,\,\,\, t\geqslant0.
\end{equation}
\item[$4.$] The family $S(t):\mathscr H^0\to\mathscr H^0$ defined by $S(t)\U_0=\U(t)$ is a strongly
continuous nonlinear semigroup on $\mathscr H^0$.
\end{itemize}
\end{Theorem}
\begin{proof}
Let $\U_0\in\mathscr{H}^0$. By Theorem \ref{thm:local}, \eqref{eq:cauchy} admits a unique maximal mild solution. By Proposition \ref{prop:energy} we have that $\sup_{0\leqslant t<T_{\max}}\|\U(t)\|_{\mathscr H^0}<\infty$ and the blow-up alternative excludes
$T_{\max}<\infty$. This proves $1$ and $2$ follows from the regularity statement of
Theorem \ref{thm:local} together with the same global bound.
To prove item $3$, set $r_D=\sup_{\U_0\in D}\|\U_0\|_{\mathscr H^0}$. The energy equivalence
\eqref{eq:energy-equiv} and the energy identity \eqref{eq:energy-id} yield a constant $R_D>0$, depending only on $r_D$, such that, every solution issuing from $D$ satisfies $
\|\U(t)\|_{\mathscr H^0}\leqslant R_D$ for every $t\geqslant0$. Consequently, Lemma \ref{lem:locdep} applies to any pair of such solutions on every interval $[0,T]$, with a Lipschitz constant $C_{R_D}$ independent of $T$. Since $T>0$ is arbitrary,
$$
\|\U^1(t)-\U^2(t)\|_{\mathscr H^0}
\leqslant
e^{C_{R_D}t}
\|\U_0^1-\U_0^2\|_{\mathscr H^0},
\,\,\,\forall t\geqslant0.
$$
Squaring this estimate and setting $C_D=2C_{R_D}$ gives \eqref{eq:cont-dep}.
Finally, existence, uniqueness and continuous dependence on the initial data imply that
$S(t)\U_0=\U(t)$ is well defined and strongly continuous, uniqueness gives the semigroup property and clearly $S(0)=I$. This proves $4$ and completes the proof.
\end{proof}
The preceding proof also provides a bound that is uniform with respect
to the initial data in bounded sets. We record its explicit dependence
on the energy constants for later use in the parameter-dependent
analysis.
\begin{Corollary}[Bounded forward orbits]\label{cor:bounded} Let $D\subset\mathscr H^0$ be bounded and set	$r_D=\sup_{\U\in {D}}\|\U\|_{\mathscr H^0}$. Then $$\|S(t)\U_0\|_{\mathscr H^0}\leqslant R_D,\,\,\,\U_0\in D,\,\,t\geqslant0,
$$
where $R_D=\left[\frac{C_2\bigl(1+{r_D}^4\bigr)+C_0}{C_1}\right]^{\frac{1}{2}}$. In particular, $R_D$ depends on $D$ only through $r_D$.
\end{Corollary}
\begin{proof}
The energy identity and Proposition \ref{prop:energy} give
$$C_1\|S(t)\U_0\|_{\mathscr H^0}^2-C_0
\leqslant
\mathscr E_{\U}(t)
\leqslant
\mathscr E_{\U}(0)
\leqslant
C_2\bigl[1+\|\U_0\|_{\mathscr H^0}^4\bigr].
$$
Since $\|\U_0\|_{\mathscr H^0}\leqslant r_D$ for $\U_0\in {D}$, the conclusion follows.
\end{proof}
\section{Global attractor}\label{sec:attractor}
In this section we investigate the long-time behavior of the nonlinear semigroup
$(\mathscr{H}^{0},S(t))$ generated by problem \eqref{eq:main-abs}. For the reader's convenience, we
first recall the main concepts and results from the theory of infinite-dimensional dynamical systems
that will be used throughout. The presentation follows the framework developed in \cite{Von}.
\subsection{Preliminaries}
\begin{Definition} Consider a dynamical system $(\mathscr{H},\{S(t)\}_{t\geqslant 0})$.
\begin{itemize}
\item[$1.$] A compact set $\A\subset\mathscr H$ is called a \underline{global attractor} if
$S(t)\A=\A$ for every $t\geqslant0$ and $\A$ attracts every bounded subset of $\mathscr H$, that is,
$\lim_{t\to\infty}\dist_{\mathscr H}(S(t)D,\A)=0$ for every bounded $D\subset\mathscr H$, where
$\dist_{\mathscr H}$ denotes the Hausdorff semidistance.
\item[$2.$] A closed set $\mathbf{B}\subset\mathscr H$ is called an \underline{absorbing set} if, for every
bounded $D\subset\mathscr H$, there exists $t_D\geqslant0$ such that $S(t)D\subset \mathbf{B}$ for all
$t\geqslant t_D$. The system is \underline{dissipative} if it possesses a bounded absorbing set.
\item[$3.$] $(\mathscr{H},S(t))$ is called a \underline{gradient} dynamical system if there exists a continuous strict Lyapunov functional $\Psi:\mathscr{H}\to\mathbb{R}$, that is, the function
$t\mapsto\Psi(S(t)\U_0)$ is non-increasing for every $\U_0\in\mathscr{H}$ and whenever
$t\mapsto\Psi(S(t)\U_0)$ is constant on $[0,\infty)$ one necessarily has $\U_0\in\mathscr{N}$, where
$\mathscr{N}$ denotes the set of stationary points.
\item[$4.$] Let $X$ and $Y$ be reflexive Banach spaces such that the embedding $X\subset Y$ is
compact and let $\mathscr H=X\times Y$. Assume that $S(t)\U_0=(u(t),u_t(t))$ with
$u\in C([0,\infty);X)\cap C^1([0,\infty);Y)$. The semigroup is said to be \underline{quasi-stable} on
a bounded set $D\subset\mathscr H$ if there exist a compact seminorm $\eta_X$ on $X$, locally bounded
nonnegative functions $a(t),c(t)$ and a function $\mathsf{b}\in L^1(\mathbb R^+)$ with
$\lim_{t\to\infty}\mathsf{b}(t)=0$, such that, for every $\U_0^1,\U_0^2\in D$,
\[
\|S(t)\U_0^1-S(t)\U_0^2\|_{\mathscr H}^2 \leqslant a(t)\|\U_0^1-\U_0^2\|_{\mathscr H}^2
\]
and
\[
\|S(t)\U_0^1-S(t)\U_0^2\|_{\mathscr H}^2 \leqslant \mathsf{b}(t)\|\U_0^1-\U_0^2\|_{\mathscr H}^2
+ c(t) \sup_{0\leqslant s\leqslant t} \eta_X\big(u^1(s)-u^2(s)\big)^2 .
\]
\end{itemize}
\end{Definition}
The following theorem summarizes the abstract results that will be used throughout this section. It
combines Corollary 7.5.7, Proposition 7.9.4, Theorem 7.9.6 and Theorem 7.9.8 in \cite{Von}.
\begin{Theorem}\label{thm:abstract}
Let $(\mathscr H,S(t))$ be a gradient dynamical system. Assume that its Lyapunov functional $\Psi$ is
bounded from above on bounded subsets of $\mathscr H$, that every sublevel set
$\Psi_R=\{\U\in\mathscr H:\Psi(\U)\leqslant R\}$ is bounded, that the set $\mathscr N$ of stationary
points is bounded and that $(\mathscr H,S(t))$ is quasi-stable on every bounded forward invariant
subset of $\mathscr H$. Then:
\begin{enumerate}
\item[$1.$] $(\mathscr H,S(t))$ possesses a compact global attractor $\A=W^u(\mathscr N)$.
\item[$2.$] If $(\mathscr H,S(t))$ is quasi-stable on $\A$, then $\A$ has finite fractal dimension.
\item[$3.$] If the function $c(t)$ in the quasi-stability estimate is bounded, then every complete
trajectory $(u(t),u_t(t))$ contained in $\A$ satisfies
$u_t\in L^\infty(\mathbb R;X)\cap C(\mathbb R;Y)$ and $u_{tt}\in L^\infty(\mathbb R;Y)$, with a bound
$\|u_t(t)\|_X^2+\|u_{tt}(t)\|_Y^2\leqslant R^2$ independent of the trajectory.
\end{enumerate}
\end{Theorem}
\subsection{Main result}
We are now in a position to establish the existence and main properties of the global attractor associated with \eqref{eq:main-abs}.
\begin{Theorem}\label{thm:attractor}
Assume {\rm(H1)}-{\rm(H5)}. Then the dynamical system $(\mathscr H^0,S(t))$ generated by \eqref{eq:main-abs} possesses a	unique compact global attractor $\mathbf A$. Moreover, the following properties hold.
\begin{enumerate}
\item[$1.$] The global attractor coincides with the unstable manifold emanating from the set $\mathscr N$ of stationary points, namely $\mathbf A=W^u(\mathscr N)$.
\item[$2.$] The attractor $\mathbf A$ has finite fractal dimension in $\mathscr H^0$.
\item[$3.$] The attractor $\mathbf A$ is bounded in the strong phase space $\mathscr H^1$.
\end{enumerate}
\end{Theorem}
\subsection{Proof of Theorem \ref{thm:attractor}}
\subsubsection{Lyapunov functional and gradient property}
\begin{Lemma}\label{lem:lyap} Assume {\rm(H1)}-{\rm(H5)}. The energy $\mathscr E_{\U}$ satisfies the following properties.
\begin{enumerate}
\item[$1.$] $\mathscr E_{\U}$ is a strict Lyapunov functional for $(\mathscr H^0,S(t))$.
\item[$2.$] $\mathscr E_{\U}$ is bounded from above on every bounded subset of $\mathscr H^0$ and each sublevel set $\{\U\in\mathscr H^0 \,|\, \mathscr E_{\U}(\U)\leqslant R\}$ is bounded in $\mathscr H^0$.
\item[$3.$] The set $\mathscr N$ of stationary points is bounded in $\mathscr H^0$.
\end{enumerate}
\end{Lemma}
\begin{proof}
Let $\U_0=(u_0,u_1)\in\mathscr H^0$. By the energy identity \eqref{eq:energy-id}, the map
$t\mapsto\Psi(S(t)\U_0)\equiv\mathscr E_{\U}(t)$ is non-increasing. If it is constant on
$[0,\infty)$, then the dissipation term in \eqref{eq:energy-id} vanishes identically, hence
$u_t(t)=0$ for all $t\geqslant 0$ and $S(t)\U_0=\U_0$ for every $t\geqslant0$. Therefore
$\U_0\in\mathscr N$ and $\Psi$ is a strict Lyapunov functional. The continuity of $\Psi$ on
$\mathscr H^0$ follows from the growth bounds established in the proof of Proposition \ref{prop:energy} together with the Sobolev embeddings.
The upper estimate in \eqref{eq:energy-equiv}, shows that $\Psi$ is bounded from above on bounded sets. If $\Psi(\U)\leqslant R$ then $\|\U\|^2_{\mathscr H^0}\leqslant (R+C_0)/C_1$, so every sublevel set is bounded.
It remains to prove that $\mathscr N$ is bounded. Let $\U=(u,0)\in\mathscr N$, then $u$ solves
\begin{equation}\label{eq:stationary}
\Le u+f(u)=\mu_1u-\HH(u)+b(x)\ \text{ in }\Omega\,\,\,\mbox{and}\,\,\, u=0\ \text{ on }\partial\Omega.
\end{equation}
Testing \eqref{eq:stationary} with $u$ and using \eqref{weak-Le} gives
$\|u\|_e^2+\langle f(u),u\rangle=\mu_1\|u\|_2^2-\langle\HH(u),u\rangle+\langle b,u\rangle$. Adding \eqref{eq:diss1} and \eqref{eq:diss2} we obtain the pointwise bound $f(u)\cdot u\geqslant-2M|u|^2-2m_f$, whence $\langle f(u),u\rangle\geqslant-2M\|u\|_2^2-2m_f|\Omega|$. Together with \eqref{eq:H-sym}, Young's inequality and \eqref{eq:coerc} this yields, for every $\varepsilon>0$
\[
\left[1-\frac{2M+6\mu_2+\mu_1}{\lambda_1\underline\mu}-\frac{\varepsilon}{4\lambda_1\underline\mu}\right]\|u\|_e^2
\leqslant 2m_f|\Omega|+\frac{1}{\varepsilon}\|b\|_2^2 .
\]
By \eqref{eq:smallness} the constant
$\beta_0:=1-\bigl(2M+6\mu_2+\mu_1\bigr)/(\lambda_1\underline\mu)$ is positive. Choosing
$\varepsilon=2\lambda_1\underline\mu\,\beta_0$ the bracket equals $\beta_0/2$ and we obtain the
explicit bound
\begin{equation}\label{eq:bound-N}
\sup_{(u,0)\in\mathscr N}\|u\|_e^2 \leqslant
\frac{2}{\beta_0}\Bigl[2m_f|\Omega|+\frac{\|b\|_2^2}{2\lambda_1\underline\mu\,\beta_0}\Bigr]:=C_{\mathscr N}^2 .
\end{equation}
Hence $\mathscr N$ is bounded in $\mathscr H^0$, with a radius depending only on the structural
constants.
\end{proof}
\subsubsection{Quasi-stability property}
\begin{Theorem}[Quasi-stability]\label{thm:quasi}
Assume {\rm(H1)-(H5)}. The dynamical system $(\mathscr H^0,S(t))$ is quasi-stable on every bounded subset $D\subset\mathscr H^0$. More precisely, for any two mild solutions
$\U^1(t)=(u(t),u_t(t))$ and $\U^2(t)=(v(t),v_t(t))$ with initial data
$\U^{1}_{0},\U^{2}_{0}\in D$, there exist $C_D>0$ and a nonnegative function $\mathsf{b}\in L^1(\mathbb R^+)$ with
$\lim_{t\to\infty}\mathsf{b}(t)=0$, such that
$$
\|\U^1(t)-\U^2(t)\|_{\mathscr H^0}^2
\leqslant \mathsf{b}(t)\|\U^{1}_{0}-\U^{2}_{0}\|_{\mathscr H^0}^2
+C_D\sup_{0\leqslant s\leqslant t}\|u(s)-v(s)\|_{\overline p}^2,$$
where $\overline p=\max\Big\{4,\frac{6}{4-p}\Big\}<6$.
\end{Theorem}
Throughout this subsection $D$ denotes a bounded subset of $\mathscr H^0$. By Corollary \ref{cor:bounded} the corresponding mild solutions remain uniformly bounded in
$\mathscr H^0$ and $C_D$ denotes a generic positive constant depending only on $D$. Setting
$\U(t)=\U^1(t)-\U^2(t)$ and $\omega(t)=u(t)-v(t)$, the function $\omega$ solves
\begin{equation} \label{eq:diff}
\begin{cases}
\rho\,\omega_{tt} + \Le(\omega) + \alpha \omega_t + f(u) - f(v) = \mu_1\omega - \HH(\omega) & \text{in } \Omega \times \mathbb{R}^{+}, \\
\omega = 0 & \text{on } \partial \Omega \times \mathbb{R}^{+}, \\
\omega(0) = u_0 - v_0, \,\, \omega_t(0) = u_1 - v_1 & \text{in } \Omega .
\end{cases}\end{equation}
The linear energy associated with \eqref{eq:diff} is given by
\begin{equation} \label{eq:E-lin}
E_{\U}(t) = \frac{1}{2}\|\omega\|_e^2 + \frac{1}{2}\|\omega_t\|_\rho^2 + \frac{1}{2} \langle \HH(\omega),\omega \rangle - \frac{\mu_1}{2}\|\omega\|_2^2.
\end{equation}
We also introduce the auxiliary functionals
\begin{equation} \label{eq:chi}
	\chi(t) = \langle \omega,  \omega_t \rangle_\rho \,\,\,\text{and}\,\,\,\Upsilon(t) = \gamma  E_\U(t) + \chi(t),
\end{equation}
where $\gamma$ is a positive constant to be chosen later.
\begin{Lemma}\label{lem:equiv}
There exist positive constants $\sigma_1$, $\sigma_2$ and $\sigma_3$ such that
\begin{equation}\label{eq:E-equiv}
\sigma_1\|\U\|_{\mathscr H^0}^2\leqslant E_{\U}(t) \leqslant \sigma_2\|\U\|_{\mathscr H^0}^2
\end{equation}
and
\begin{equation}\label{eq:Ups-equiv0}
\big[\gamma-\sigma_3\big]E_{\U}(t)\leqslant \Upsilon(t)\leqslant \big[\gamma+\sigma_3\big]E_{\U}(t).
\end{equation}
\end{Lemma}
\begin{proof}
By \eqref{eq:H-sym} and Poincar\'e inequality combined with \eqref{eq:coerc},
\[
\left[\frac12-\frac{\mu_1+6\mu_2}{2\underline{\mu}\lambda_1}\right]\|\omega\|_e^2+\frac12\|\omega_t\|_\rho^2
\leqslant E_{\U}(t)\leqslant\left[\frac12+\frac{3\mu_2}{\underline{\mu}\lambda_1}\right]\|\omega\|_e^2
+\frac12\|\omega_t\|_\rho^2 ,
\]
and, in view of \eqref{eq:smallness},
$\tfrac12-\tfrac{\mu_1+6\mu_2}{2\underline{\mu}\lambda_1}>\tfrac{M}{\underline{\mu}\lambda_1}>0 $. This gives \eqref{eq:E-equiv}. Next, by the Cauchy-Schwarz and Young inequalities together with \eqref{eq:rho-equiv} and \eqref{eq:coerc}, we find
$$|\chi(t)|=\big|\langle\omega,\omega_t\rangle_\rho\big|\leqslant\|\omega\|_\rho\|\omega_t\|_\rho
\leqslant\frac{\overline\rho}{2}\Big[\frac{1}{\underline\mu\lambda_1}\|\omega\|_e^2+\|\omega_t\|_2^2\Big]
\leqslant \sigma_3E_{\U}(t).$$
Recalling that $\Upsilon(t)=\gamma E_{\U}(t)+\chi(t)$, we obtain \eqref{eq:Ups-equiv0}.
\end{proof}
\begin{Lemma}\label{lem:diff-est}
Let $D\subset\mathscr H^0$ be bounded. Then there exist positive constants $C_4$, $C_5$ and $C_D$
such that
\begin{equation}\label{eq:E'}
E_{\U}'(t)\leqslant -\frac{\alpha}{2}\|\omega_t\|_2^2-\langle\omega_t,h(u)-h(v)\rangle
+C_D\|\omega\|_{\overline p}^2
\end{equation}
and
\begin{equation}\label{eq:chi'}
\chi'(t)\leqslant -\frac{C_4}{2\sigma_2}E_{\U}(t)-\frac{C_4\underline{\mu}}{2}\|\nabla\omega\|_2^2
+C_5\|\omega_t\|_2^2+C_D\|\omega\|_4^2 .
\end{equation}
\end{Lemma}
\begin{proof}
Differentiating $E_{\U}$ and using \eqref{eq:diff}, we find
$E_{\U}'(t) = -\alpha\|\omega_t\|_2^2 -\langle\omega_t,g(u)-g(v)\rangle-\langle\omega_t,h(u)-h(v)\rangle$.
By \eqref{eq:growth-g} and the mean value theorem,
$|g_i(u)-g_i(v)|\leqslant C_g[1+\sum_j(|u_j|^{p-1}+|v_j|^{p-1})]|\omega|$, so H\"older  inequality with
$\tfrac{p-1}{6}+\tfrac{4-p}{6}+\tfrac12=1$ give
$$|\langle\omega_t,g(u)-g(v)\rangle|\leqslant C_D\|\omega\|_{6/(4-p)}\|\omega_t\|_2.$$ Young  inequality and $\overline p\geqslant \tfrac{6}{4-p}$ then yield \eqref{eq:E'}.
For $\chi$, differentiating \eqref{eq:chi} and using \eqref{eq:diff},
\begin{align*}
\chi'(t)= \|\omega_t\|_\rho^2 -\|\omega\|_e^2-\alpha\langle\omega,\omega_t\rangle
-\langle f(u)-f(v),\omega\rangle+\mu_1\|\omega\|_2^2-\langle\HH(\omega),\omega\rangle .
\end{align*}
By \eqref{eq:H-sym} and \eqref{eq:coerc},
$\mu_1\|\omega\|_2^2-\langle\HH(\omega),\omega\rangle\leqslant\tfrac{\mu_1+6\mu_2}{\underline\mu\lambda_1}\|\omega\|_e^2$,
and Young's inequality applied to $\alpha\langle\omega,\omega_t\rangle$ gives
\[
\chi'(t)\leqslant-C_4\|\omega\|_e^2+C\|\omega_t\|_2^2-\langle f(u)-f(v),\omega\rangle,
\]
where $C_4=\frac12\left(1-\frac{\mu_1+6\mu_2}{\underline{\mu}\lambda_1}\right)>0$.  We now estimate the nonlinear term. Since $p<3$, H\"older's inequality with
$\tfrac{p-1}{6}+2\cdot\tfrac{7-p}{12}=1$ and $\tfrac{12}{7-p}\leqslant 4$ give
\begin{align}\label{eq:estg}
|\langle g(u)-g(v),\omega\rangle|\leqslant C_D\|\omega\|_4^2 ,
\end{align}
while from \eqref{eq:growth-h}, $|h_j(u_j)-h_j(v_j)|\leqslant C[1+|u_j|^2+|v_j|^2]|\omega_j|$ and H\"older's inequality with $\tfrac26+\tfrac23=1$ give
$|\langle h_j(u_j)-h_j(v_j),\omega_j\rangle|\leqslant C\|\omega_j\|_2^2+C[\|u_j\|_6^2+\|v_j\|_6^2]\|\omega_j\|_3^2$. Using $L^4(\Omega)\subset L^3(\Omega)\subset L^2(\Omega)$,
\begin{align}\label{eq:esth}
\Big|\sum_{j=1}^3\langle h_j(u_j)-h_j(v_j),\omega_j\rangle\Big|\leqslant C_D\|\omega\|_4^2 .
\end{align}
Combining \eqref{eq:estg}-\eqref{eq:esth}, we obtain
$|\langle f(u)-f(v),\omega\rangle|\leqslant C_D\|\omega\|_4^2$.
Finally we recover the energy. Splitting the term $-C_4\|\omega\|_e^2$ into two halves and using
$E_\U\leqslant \sigma_2\|\U\|^2_{\mathscr H^0}=\sigma_2\big(\|\omega\|_e^2+\|\omega_t\|_\rho^2\big)$ from
\eqref{eq:E-equiv},
\[
-\frac{C_4}{2}\|\omega\|_e^2=-\frac{C_4}{2}\|\U\|^2_{\mathscr H^0}+\frac{C_4}{2}\|\omega_t\|_\rho^2
\leqslant-\frac{C_4}{2\sigma_2}E_{\U}(t)+\frac{C_4\overline\rho}{2}\|\omega_t\|_2^2 ,
\]
the surplus being absorbed into $C_5\|\omega_t\|_2^2$, whereas for the remaining half
$-\tfrac{C_4}{2}\|\omega\|_e^2\leqslant-\tfrac{C_4\underline\mu}{2}\|\nabla\omega\|_2^2$ by
\eqref{eq:coerc}. This gives \eqref{eq:chi'}.
\end{proof}
A crucial ingredient in the proof of Theorem \ref{thm:quasi} is the control of the terms involving the critical component $h$ of the nonlinearity. To this end we use the following auxiliary result,
established in \cite[Lemma 2.1]{Ma} for the critical case.
\begin{Lemma}[{\cite[Lemma 2.1]{Ma}}] \label{lem:Ma}
There exists a positive constant $C=C_D$, such that, for every $\eta>0$,
\begin{align*}
&\int_0^t e^{\eta s}\big\langle h(v(s))-h(u(s)),\omega_t(s)\big\rangle ds\\
& \qquad \qquad \qquad
\leqslant C e^{\eta t}\sup_{0<s<t}\|\omega\|_4^2
+ C\sum_{j=1}^3\int_0^t\big(\|\partial_tu_j(s)\|_2+\|\partial_tv_j(s)\|_2\big)e^{\eta s}\|\nabla\omega_j(s)\|_2^2\,ds.
\end{align*}
\end{Lemma}
With the above lemmas established, we can now proceed to finalize the proof of the Theorem \ref{thm:quasi}.
\begin{proof}[Proof of Theorem \ref{thm:quasi}]
Differentiating $\Upsilon$ and combining \eqref{eq:E'} and \eqref{eq:chi'},
\begin{align}\begin{split}\label{eq:Ups'}
\Upsilon'(t)&\leqslant-\frac{C_4}{2\sigma_2}E_{\U}(t)-\frac{C_4\underline{\mu}}{2}\|\nabla\omega\|_2^2
-\left[\frac{\gamma\alpha}{2}-C_5\right]\|\omega_t\|_2^2\\
&\quad-\gamma\langle\omega_t,h(u)-h(v)\rangle+(1+\gamma)C_D\|\omega\|_{\overline p}^2 .
\end{split}\end{align}
We now fix $\gamma>\max\{\sigma_3,\,2C_5/\alpha\}$, so that $\Upsilon$ is equivalent to $E_{\U}$ by
\eqref{eq:Ups-equiv0} and the coefficient of $\|\omega_t\|_2^2$ in \eqref{eq:Ups'} is nonpositive.
Using the upper estimate in \eqref{eq:Ups-equiv0},
\begin{equation} \label{eq:Ups-diff}
\Upsilon'(t)\leqslant-\eta\,\Upsilon(t)-\frac{C_4\underline{\mu}}{2}\|\nabla\omega\|_2^2
-\gamma\langle\omega_t,h(u)-h(v)\rangle+C_D\|\omega\|_{\overline p}^2,
\,\,\,\eta=\frac{C_4}{2\sigma_2(\gamma+\sigma_3)}>0 .
\end{equation}
Notice that $\eta$ is independent of the bounded set $D$. Multiplying \eqref{eq:Ups-diff} by
$e^{\eta t}$ and integrating over $(0,t)$ we obtain
\begin{align}\begin{split}\label{eq:Ups-int}
e^{\eta t}\Upsilon(t)&\leqslant \Upsilon(0)-\frac{C_4\underline{\mu}}{2}\int_0^te^{\eta s}\|\nabla\omega(s)\|_2^2\,ds
+C_D\int_0^te^{\eta s}\|\omega(s)\|_{\overline p}^2\,ds + \gamma\,J(t),
\end{split}\end{align}
where $J(t)=\int_0^te^{\eta s}\langle\omega_t(s),h(v(s))-h(u(s))\rangle\,ds$.
Applying Lemma \ref{lem:Ma} and, for every $\varepsilon>0$, Young's inequality in the form
\[
C\big(\|u_t\|_2+\|v_t\|_2\big)\|\nabla\omega\|_2^2
\leqslant\varepsilon\|\nabla\omega\|_2^2+\frac{C^2}{2\varepsilon}\big(\|u_t\|_2^2+\|v_t\|_2^2\big)\|\nabla\omega\|_2^2,
\]
together with \eqref{eq:E-equiv} and \eqref{eq:Ups-equiv0} give
\begin{align}\begin{split}\label{eq:J-est}
J(t)&\leqslant Ce^{\eta t}\sup_{0\leqslant s\leqslant t}\|\omega(s)\|_4^2
+\varepsilon\int_0^te^{\eta s}\|\nabla\omega(s)\|_2^2\,ds
+\frac{C^2}{2\varepsilon}\int_0^t\!\! \big[\|u_t(s)\|_2^2+\|v_t(s)\|_2^2\big]e^{\eta s}\Upsilon(s)\,ds .
\end{split}\end{align}
Choosing $\varepsilon=\tfrac{C_4\underline{\mu}}{4\gamma}$, the gradient term in \eqref{eq:J-est} is absorbed by the negative gradient term in \eqref{eq:Ups-int}. Therefore
\begin{align*}
e^{\eta t}\Upsilon(t)\leqslant \Upsilon(0) +C_D e^{\eta t}\sup_{0\leqslant s\leqslant t}\|\omega(s)\|_{\overline p}^2
+C\int_0^t\big[\|u_t(s)\|_2^2+\|v_t(s)\|_2^2\big]e^{\eta s}\Upsilon(s)\,ds .
\end{align*}
By Gronwall's inequality,
\begin{align}\label{eq:quasi1}
e^{\eta t}\Upsilon(t)\leqslant\left[\Upsilon(0)+C_D e^{\eta t}\sup_{0\leqslant s\leqslant t}\|\omega(s)\|_{\overline p}^2\right]
\exp\left\{C\int_0^t\big[\|u_t(s)\|_2^2+\|v_t(s)\|_2^2\big]ds\right\} .
\end{align}
Finally, by the energy identity \eqref{eq:energy-id} applied to each trajectory,
\[
\int_0^t\big[\|u_t(s)\|_2^2+\|v_t(s)\|_2^2\big]ds
=\frac{1}{\alpha}\Big[\mathscr{E}_{\U^{1}}(0)-\mathscr{E}_{\U^{1}}(t)+\mathscr{E}_{\U^{2}}(0)-\mathscr{E}_{\U^{2}}(t)\Big]
\leqslant C_D,\,\,\,\forall\,t\geqslant0.
\]
Since the initial data belong to the bounded set $D$, the corresponding
trajectories are uniformly bounded in $\mathscr H^0$. Moreover,
\eqref{eq:energy-equiv} provides a uniform lower bound for the energy
along these trajectories. Therefore, \eqref{eq:quasi1} yields
\begin{align}\label{quasi}
\Upsilon(t)
\leqslant C_D e^{-\eta t}\Upsilon(0)+C_D \sup_{0\leqslant s\leqslant t} \|\omega(s)\|_{\overline p}^2.
\end{align}
Combining \eqref{quasi} with \eqref{eq:E-equiv} and
\eqref{eq:Ups-equiv0}, we obtain the desired quasi-stability estimate
with $\mathsf{b}(t)=C_D e^{-\eta t}$, which belongs to $L^1(\mathbb R^+)$ and
$\eta_X(\omega)=\|\omega\|_{\overline p}$. This is a compact seminorm
on $X=D(\Le^{1/2})=[H_0^1(\Omega)]^3$, since $\overline p<6$ and the
embedding $[H_0^1(\Omega)]^3\subset
[L^{\overline p}(\Omega)]^3$ is compact. Finally, the first
quasi-stability estimate is precisely \eqref{eq:cont-dep}.
\end{proof}
\begin{Remark}\label{rem:homog}
Consider the homogeneous problem obtained from \eqref{eq:cauchy} by setting $\Ff(\U)=0$. Its energy is $\widetilde E(t)=\tfrac12\|\omega\|_e^2+\tfrac12\|\omega_t\|_\rho^2$ and a direct computation yields
$$
\widetilde E'(t)=-\alpha\|\omega_t\|_2^2\,\,\,\mbox{and}\,\,\,
\widetilde\chi'(t)=
\|\omega_t\|_\rho^2-\|\omega\|_e^2
-\alpha\langle\omega,\omega_t\rangle
\leqslant
-\tfrac12\|\omega\|_e^2+C\|\omega_t\|_2^2,
$$
where $\widetilde\chi(t)=\langle\omega,\omega_t\rangle_\rho$ and $\omega$ now denotes the solution of the homogeneous problem.
Defining $\widetilde\Upsilon(t)=\gamma\widetilde E(t)+\widetilde\chi(t)$ and choosing $\gamma>0$ sufficiently large, we obtain  $\widetilde\Upsilon'(t)\leqslant-c\widetilde\Upsilon(t)$ for some $c>0$. Therefore, Gronwall's inequality and the equivalence of $\widetilde\Upsilon$ with the squared norm of $\mathscr H^0$ imply that there exist constants $C,\eta>0$, independent of the initial	data, such that
\begin{equation}\label{eq:lin-decay}
\|S_h(t)\U_0\|_{\mathscr H^0}
\leqslant Ce^{-\eta t}\|\U_0\|_{\mathscr H^0},	\,\,\, t\geqslant 0.
\end{equation}
Thus, the homogeneous semigroup $\{S_h(t)\}_{t\geqslant0}$ is uniformly exponentially stable on $\mathscr H^0$.
\end{Remark}
\subsubsection{Completion of the proof of Theorem \ref{thm:attractor}}
\begin{proof} By Lemma \ref{lem:lyap}, $(\mathscr H^0,S(t))$ is a gradient dynamical system whose Lyapunov functional is bounded from above on bounded subsets of $\mathscr H^0$ and has bounded sublevel sets. Moreover, the set $\mathscr N$ of stationary points is bounded. By Theorem \ref{thm:quasi}, the system is quasi-stable on every bounded forward invariant subset of $\mathscr H^0$. Therefore, all the hypotheses of Theorem~\ref{thm:abstract} are satisfied.
Theorem \ref{thm:abstract}$(1)$ yields the existence of the compact	global attractor $\mathbf A=W^u(\mathscr N)$. Since the quasi-stability estimate holds in particular on $\mathbf A$, Theorem \ref{thm:abstract}$(2)$ implies	that $\mathbf A$ has finite fractal dimension in $\mathscr H^0$.
It remains to prove that $\mathbf A$ is bounded in $\mathscr H^1$. Let $(u(t),u_t(t))$, $t\in\mathbb R$, be a complete trajectory in $\mathbf A$. Since the function $c(t)$ in the quasi-stability estimate is constant, Theorem~\ref{thm:abstract}$(3)$ gives
$u_t\in L^\infty(\mathbb R;[H_0^1(\Omega)]^3)\cap C(\mathbb R;[L^2(\Omega)]^3)$ and $u_{tt}\in L^\infty(\mathbb R;[L^2(\Omega)]^3)$,
with bounds uniform over all complete trajectories in $\mathbf A$. Since $\mathbf A$ is bounded in $\mathscr H^0$, the function $u(t)$ is uniformly bounded in $[H_0^1(\Omega)]^3$ and hence in $[L^6(\Omega)]^3$. By \eqref{eq:growth-g}-\eqref{eq:growth-h}, we have	$|f(u)|\leqslant C[|u|+|u|^p+|u|^3]$ with $p<3$. Consequently, $f(u(t))$ is uniformly bounded in $[L^2(\Omega)]^3$. Writing
$$
\Le u(t)=	-\rho u_{tt}(t)-\alpha u_t(t)-f(u(t))+\mu_1u(t)-\HH u(t)+b,
$$
we see that every term on the right-hand side is uniformly bounded in $[L^2(\Omega)]^3$. Therefore, 	$u\in L^\infty(\mathbb R;D(\Le))$. Moreover, the estimates for $u_t$ give $u_t\in L^\infty(\mathbb R;D(\Le^{1/2}))$. Hence, $(u,u_t)\in L^\infty(\mathbb R;\mathscr H^1)$,	with a bound uniform over all complete trajectories in $\mathbf A$. Since every point of $\mathbf A$ lies on a complete trajectory, the	attractor $\mathbf A$ is bounded in $\mathscr H^1$.
\end{proof}
\subsection{Quantitative dissipativity}
The existence of the global attractor was established above through
the gradient structure and asymptotic compactness, without requiring
an explicit absorbing set. For later use, particularly in the
parameter-dependent analysis, we now derive quantitative dissipativity
estimates.
\begin{Corollary}[Explicit bound for the attractor]\label{cor:attr-bound}
Assume {\rm(H1)-(H5)}. Then
\begin{equation}\label{eq:CL752}
\sup_{\U\in\mathbf A}\mathscr E_{\U}
\leqslant \sup_{\U\in\mathscr N}\mathscr E_{\U} \leqslant C_2\Bigl[1+C_{\mathscr N}^{4}\Bigr],
\end{equation}
where $C_{\mathscr N}$ is given by \eqref{eq:bound-N}. Consequently,
\begin{equation}\label{eq:radiusA}
\|\U\|_{\mathscr H^0} \leqslant
\left[\frac{C_2\bigl[1+C_{\mathscr N}^4\bigr]+C_0}{C_1}\right]^{\frac{1}{2}}=: C_{\mathbf A},\,\,\, \U\in\mathbf A.
\end{equation}
In particular, every ball in $\mathscr H^0$ centered at the origin with radius larger than $C_{\mathbf A}$ is absorbing. The bound	$C_{\mathbf A}$ depends only on $\alpha$, $\underline\rho$, $\overline\rho$, $\underline\mu$, $\lambda_1$, $M$, $m_f$, $\mu_1$, $\mu_2$, $|\Omega|$ and $\|b\|_2$.
\end{Corollary}
\begin{proof}
Let $\U\in\mathbf A$. By Theorem \ref{thm:attractor}, $\mathbf A=W^u(\mathscr N)$ and hence
$\U=\Gamma(0)$ for some complete trajectory ${\Gamma}\subset\mathbf A$ satisfying
$\dist_{\mathscr H^0}({\Gamma}(t),\mathscr N)\to0$ as $t\to-\infty$. Since the energy is non-increasing along trajectories, \eqref{eq:energy-id} gives $\mathscr E_{{\Gamma}(0)}\leqslant\mathscr E_{{\Gamma}(t)}$ for every $t\leqslant0$.

Choose $\V_t\in\mathscr N$ such that $\|{\Gamma}(t)-\V_t\|_{\mathscr H^0}\to0$ as $t\to-\infty$. Since $\mathbf A$ and $\mathscr N$ are bounded and the energy is uniformly continuous on bounded subsets of $\mathscr H^0$, the growth estimates used in Proposition \ref{prop:energy} yield $\mathscr E_{{\Gamma}(t)}=\mathscr E_{\V_t}+o(1)$ as $t\to-\infty$. Therefore, $\mathscr E_{\U}\leqslant\sup_{\V\in\mathscr N}\mathscr E_{\V}$, which proves the first inequality in \eqref{eq:CL752}. The second follows from the upper bound in \eqref{eq:energy-equiv} and \eqref{eq:bound-N}, since every stationary point has the form $\V=(u,0)$ and satisfies $\|\V\|_{\mathscr H^0}=\|u\|_e\leqslant C_{\mathscr N}$.

The lower bound in \eqref{eq:energy-equiv} now gives \eqref{eq:radiusA}. Finally, if $r>C_{\mathbf A}$, the ball $B_{\mathscr H^0}(0,r)$ contains a neighbourhood of $\mathbf A$. Since $\mathbf A$ attracts every bounded subset of $\mathscr H^0$, this ball is absorbing.
\end{proof}
Although Corollary \ref{cor:attr-bound} shows that every ball whose
radius is larger than $C_{\mathbf A}$ is absorbing, it does not
provide any control of the corresponding entrance time. Such control
will be needed in the parameter-dependent analysis. We therefore
derive a second absorbing ball directly from the energy estimates,
together with an explicit entrance time.
\begin{Proposition}[Absorbing ball]\label{prop:absorbing}
Assume {\rm(H1)-(H5)}. Then $(\mathscr H^0,S(t))$ is dissipative. More precisely, there exists $R_0>0$, depending only on $\alpha$, $\underline\rho$, $\overline\rho$, $\underline\mu$, $\lambda_1$, $M$, $m_f$, $\mu_1$, $\mu_2$, $|\Omega|$ and $\|b\|_2$, such that the closed ball
$$
\mathbf B_0=\bigl\{\U\in\mathscr H^0\,|\, \|\U\|_{\mathscr H^0}\leqslant R_0\bigr\}$$
absorbs every bounded subset of $\mathscr H^0$. Moreover, for every bounded set ${D}\subset\mathscr H^0$, the entrance time ${t_D}$ satisfying $S(t)D\subseteq\mathbf B_0$ for all $t\geqslant t_D$ can be chosen as a nondecreasing function of $r(D)=\sup_{\U\in D}\|\U\|_{\mathscr H^0}$.
\end{Proposition}
\begin{proof} To derive a dissipative estimate, we introduce the perturbed energy functional
\begin{equation}\label{eq:Upsilon}
\Upsilon(t)=\mathscr E_{\U}(t)+\gamma\langle u(t),u_t(t)\rangle_\rho,\,\,\,\gamma>0.
\end{equation}

\underline{Step 1. Two-sided bounds.} By \eqref{eq:coerc}, \eqref{eq:rho-equiv} and Poincar\'e inequality,	$|\langle u,u_t\rangle_\rho|\leqslant C\|\U\|_{\mathscr H^0}^2$, where $C>0$ is independent of $\U$. Since $\gamma$ will be fixed below, we may assume at this stage that $\gamma\leqslant\min\{1,C_1/(2C)\}$. Then \eqref{eq:energy-equiv} yields
\begin{equation}\label{eq:Upsilon-equiv}
\frac{C_1}{2}\|\U\|_{\mathscr H^0}^2-C_0\leqslant\Upsilon\leqslant C_2'\bigl[1+\|\U\|_{\mathscr H^0}^4\bigr],
\end{equation}
where $C_2'=C_2+C$.

\underline{Step 2. Strict decay.}
Differentiating \eqref{eq:Upsilon} and using
\eqref{eq:main-abs} and \eqref{eq:energy-id}, we obtain
$$
\Upsilon'(t)
=
-\alpha\|u_t\|_2^2
+\gamma\Bigl[
\|u_t\|_\rho^2-\|u\|_e^2
-\alpha\langle u,u_t\rangle
-\langle f(u),u\rangle
+\mu_1\|u\|_2^2-\langle\HH u,u\rangle
+\langle b,u\rangle
\Bigr].
$$
Adding \eqref{eq:diss1} and \eqref{eq:diss2} gives
$\langle f(u),u\rangle
\geqslant-2M\|u\|_2^2-2m_f|\Omega|$, while
\eqref{eq:H-sym} yields
$\mu_1\|u\|_2^2-\langle\HH u,u\rangle
\leqslant[\mu_1+6\mu_2]\|u\|_2^2$. Moreover,
$$
[2M+\mu_1+6\mu_2]\|u\|_2^2
\leqslant
[1-3C_3]\|u\|_e^2,
$$
where $C_3=
\frac13\left[
1-\frac{2M+\mu_1+6\mu_2}
{\underline\mu\lambda_1}
\right]>0 $ by \eqref{eq:smallness}. Young  and Poincar\'e  inequalities also
give
$\alpha|\langle u,u_t\rangle|
\leqslant C_3\|u\|_e^2+C\|u_t\|_2^2$ and
$|\langle b,u\rangle|
\leqslant C_3\|u\|_e^2+C\|b\|_2^2$. Consequently,
$$
\Upsilon'(t)
\leqslant
-\alpha\|u_t\|_2^2
+\gamma\Bigl[
C'_4\|u_t\|_2^2-C_3\|u\|_e^2+C_5
\Bigr],
$$
where $C'_4,C_5>0$ depend only on the structural constants. We now fix
$\gamma>0$ so that $\gamma\leqslant \min\left[1,\frac{C_1}{2C},\frac{\alpha}{2C'_4}\right]$. This choice is compatible with Step 1 and gives
$$
\Upsilon'(t)\leqslant-\frac{\alpha}{2}\|u_t\|_2^2-\gamma C_3\|u\|_e^2+\gamma C_5.
$$
Since
$\|u_t\|_\rho^2\leqslant\overline\rho\|u_t\|_2^2$, it follows that
\begin{equation}\label{eq:Upsilon-decay}
\Upsilon'(t)\leqslant-C_6\|\U(t)\|_{\mathscr H^0}^2+C_7,
\end{equation}
where $C_6=\min\left[\frac{\alpha}{2\overline\rho},\gamma C_3\right]>0$ and $C_7=\gamma C_5$.

\underline{Step 3. Absorption.}
Fix $\gamma$ as above and set
$$
\Upsilon_\infty
=C_2'\left[1+\left(\frac{2C_7}{C_6}\right)^2\right].
$$
If $\Upsilon(t)>\Upsilon_\infty$, the upper bound in
\eqref{eq:Upsilon-equiv} implies
$\|\U(t)\|_{\mathscr H^0}^2>2C_7/C_6$. Hence,
\eqref{eq:Upsilon-decay} gives
$\Upsilon'(t)\leqslant-C_7<0$. Consequently, the sublevel set
$\{\U \in \mathscr{H}^{0}\,|\, \Upsilon(\U)\leqslant\Upsilon_\infty\}$ is positively invariant and
every trajectory starting outside it enters it within a time not exceeding $[\Upsilon(0)-\Upsilon_\infty]/C_7$.
By the lower bound in \eqref{eq:Upsilon-equiv}, this sublevel set is contained in $\mathbf B_0$, where $R_0^2=\frac{2[\Upsilon_\infty+C_0]}{C_1}$. Now let $D\subset\mathscr H^0$ be bounded and set
$r(D)=\sup_{\U\in D}\|\U\|_{\mathscr H^0}$. The upper bound in
\eqref{eq:Upsilon-equiv} gives
$\Upsilon(0)\leqslant C_2'[1+r(D)^4]$ for every initial datum in
$D$. Therefore, an entrance time for $D$ is
$$
t_D=
\max\left\{0,\frac{C_2'[1+r(D)^4]-\Upsilon_\infty}{C_7}\right\}.
$$
This time depends only on $r(D)$ and is nondecreasing as a function of $r(D)$.
\end{proof}
\section{Exponential attractor}\label{sec:exp}
In this section we strengthen the asymptotic description obtained above by proving the existence of an exponential attractor for the semigroup generated by \eqref{eq:main-abs}. We begin by recalling the precise definition.
\begin{Definition}\label{exp}
Let $H$ be a Banach space and $\{S(t)\}_{t\geqslant0}$ a semigroup on $H$. A compact set $\A^{\exp}\subset H$ is called an \underline{exponential attractor} for $(H,S(t))$ if it has finite fractal dimension, is positively invariant and attracts bounded subsets of $H$ exponentially. More precisely, $S(t)\A^{\exp}\subseteq\A^{\exp}$ for every $t\geqslant0$ and there exist $\gamma>0$ and a nondecreasing function $C:\mathbb R^+\to\mathbb R^+$ such that
$$\dist_H\bigl(S(t)D,\A^{\exp}\bigr)\leqslant C(r(D))e^{-\gamma t}
$$ for every bounded set $D\subset H$ and every $t\geqslant0$, where $r({D})=\sup_{\U\in D}\|\U\|_H$.
\end{Definition}
Our proof is based on the following abstract criterion, established in \cite{EFNT,EMZ}
(see also \cite[\S7.9]{Von}).
\begin{Theorem} \label{thm:exp-crit}
Let $H_0$ and $H_1$ be two Hilbert spaces such that $H_1$ is compactly embedded into $H_0$. Let
$\{S(t)\}_{t\geqslant0}$ be a $C_0$-semigroup on $H_{0}$ and suppose that there exist a bounded
positively invariant set $\mathbf{D}\subset H_0$ and a time $t_*>0$ satisfying
\begin{itemize}
\item[$1.$] The map $[0,t_*]\times \mathbf{D} \ni (t,\U_0) \mapsto S(t)\U_0 \in \mathbf{D}$ is Lipschitz with constant $L_*>0$.
\item[$2.$] The map $S(t_*)|_{\mathbf{D}}: \mathbf{D} \to {H_0}$ admits a decomposition $S(t_*)|_{{\mathbf{D}}}=S_0+S_1$ with $S_0: \mathbf{D}\to H_0$ and $S_1: \mathbf{D}\to H_1$ such that, for all $\U^1,\U^2\in \mathbf{D}$,
\begin{equation*}
\left\|S_0 \U^1\!-\! S_0 \U^2 \right\|_{H_0}\leqslant \frac{1}{8}\left\| \U^1 \!-\! \U^2\right\|_{H_0}
\,\text{and}\,\,\,
\left\|S_1 \U^1 \!-\! S_1 \U^2 \right\|_{H_1}\leqslant C_*\left\| \U^1 \!-\! \U^2\right\|_{H_0},
\end{equation*}
where $C_*>0$.
\end{itemize}
Then there exists a positively invariant compact set $\mathbf{K}\subset \mathbf{D}$ such that
$$
\dist_{H_0}(S(t)\mathbf{D},\mathbf{K})\leqslant C e^{-\nu t},\,\,\, \forall\, t\geqslant 0,
$$
where $C=2L_*\sup_{\U\in{\mathbf{D}}}\|\U\|_{H_0}e^{\nu}$ and $\nu=\tfrac{\log2}{t_*}$. Moreover,
$\dim_F(\mathbf{K})\leqslant 1+\tfrac{\log N^*}{\log2}$, where $N^*$ is the minimum number of balls
of radius $\tfrac{1}{8{C_*}}$ of $H_0$ necessary to cover the unit ball of $H_1$.
\end{Theorem}
The proof also makes use of the transitivity of exponential attraction, which we state for
convenience (see \cite[Theorem 5.1]{FGMZ}).
\begin{Lemma}\label{lem:trans}
Let $\mathbf{K}_1,\mathbf{K}_2,\mathbf{K}_3$ be subsets of $H$ such that
$\dist_{H}(S(t)\mathbf{K}_1,\mathbf{K}_2)\leqslant L_1 e^{-\vartheta_1 t}$ and
$\dist_{H}(S(t)\mathbf{K}_2,\mathbf{K}_3)\leqslant L_2 e^{-\vartheta_2 t}$
for some $\vartheta_1,\vartheta_2>0$ and $L_1,L_2\geqslant0$. Assume also that
$\|S(t)\U^1 - S(t)\U^2\|_H \leqslant L_0 e^{\vartheta_0 t}\|\U^1 - \U^2\|_H$ for all
$\U^1,\U^2\in\bigcup_{t\geqslant0}S(t)\mathbf{K}_j$ $(j=1,2,3)$ and some $\vartheta_0,L_0\geqslant0$. Then
$\dist_{H}(S(t)\mathbf{K}_1,\mathbf{K}_3)\leqslant Le^{-\vartheta t}$ for all $t\geqslant0$, where
$\vartheta=\tfrac{\vartheta_1\vartheta_2}{\vartheta_0+\vartheta_1+\vartheta_2}$ and $L=L_0L_1+L_2$.
\end{Lemma}
\subsection{Main result}
For the construction of the exponential attractor, we impose the following additional assumptions.
\begin{itemize}
\item[{\rm(H6)}] The critical component vanishes, $h\equiv0$, the density is constant and normalized so that $\rho\equiv1$ and the growth exponent satisfies $2\leqslant p<3$.
\item[{\rm(H7)}] The functions $g$ and $\mu$ satisfy $g_i\in C^2(\mathbb R^3)$, $i=1,2,3$ and $\mu\in C^3(\overline{\Omega})$. Moreover, there exists $M_g>0$ such that
$$\left|\frac{\partial^2g_i}{\partial u_k\partial u_l}(u)\right|
\leqslant M_g\Bigl[	1+|u_1|^{p-2}+|u_2|^{p-2}+|u_3|^{p-2}
	\Bigr], \,\,\,\forall\, u\in\mathbb R^3, \,i,k,l\in\{1,2,3\}.$$
\end{itemize}
Under {\rm(H6)}-{\rm(H7)}, the nonlinear source is subcritical and
has the additional regularity required for the smoothing and
compactness estimates developed below.
\begin{Remark}[On the normalization $\rho\equiv1$]\label{rem:rho} The normalization $\rho\equiv1$ is \textit{not} merely notational. The estimates developed in this section are obtained by testing the	equation with powers of $\Le$ and rely on the fact that $\Le$ is 	positive and self-adjoint in $[L^2(\Omega)]^3$. If $\rho$ is  variable, the corresponding operator is $\rho^{-1}\Le$, which is not self-adjoint in $[L^2(\Omega)]^3$ and its fractional powers cannot be directly identified with those of $\Le$. Nevertheless, the restriction $\rho\equiv1$ can be removed by working in the weighted space $L_\rho^2=([L^2(\Omega)]^3,\langle\cdot,\cdot\rangle_\rho)$. Indeed, define $\Le^\rho=\rho^{-1}\Le$ with $D(\Le^\rho)=D(\Le)$. Then $\Le^\rho$ is positive and self-adjoint in $L_\rho^2$, since $\langle\Le^\rho u,v\rangle_\rho =\langle\Le u,v\rangle=\langle u,v\rangle_e$ for $u\in D(\Le)$ and $v\in D(\Le^{1/2})$. Moreover, $L^2_\rho$ and $[L^2(\Omega)]^3$ coincide as sets and carry equivalent norms by \eqref{eq:rho-equiv}, while	$D(\Le^\rho)=D(\Le)$ with equivalent graph norms. The interpolation characterization of the fractional domains of positive self-adjoint 	operators therefore yields $D((\Le^\rho)^r)=D(\Le^r)$ with equivalent norms for every	$r\in[0,1]$. Consequently, the arguments below extend to a variable  density after replacing $\Le$ by $\Le^\rho$ and the standard inner product by $\langle\cdot,\cdot\rangle_\rho$. We assume 	$\rho\equiv1$ only to avoid carrying this weighted formulation throughout the section.
\end{Remark}
\medskip
We are now in a position to state the main result of this section.
\begin{Theorem}\label{thm:exp} Assume {\rm(H1)}-{\rm(H7)}. Then the dynamical system 	$(\mathscr H^0,S(t))$ generated by \eqref{eq:main-abs} possesses an	exponential attractor $\A^{\exp}$.  Moreover, $\mathbf A\subseteq\A^{\exp}$.
\end{Theorem}
\begin{Remark}[Comparison with the strongly damped wave equation]	\label{rem:PS}
The decomposition and finite bootstrap used here parallel the strategy of Pata and Squassina \cite{PS}, but the source of the regularity gain is different. For the strongly	damped wave equation, the linear semigroup is analytic and therefore regularizes instantaneously. In the present problem, the damping is of	order zero and the exponentially stable linear semigroup has no smoothing effect. The regularity gain comes instead from the	inhomogeneous component under the subcritical assumptions	{\rm(H6)}-{\rm(H7)}. Consequently, each bootstrap increment is limited by $(3-p)/2$ and degenerates as $p\to3^{-}$.
\end{Remark}
\subsection{Proof of Theorem \ref{thm:exp}}
\subsubsection{Smoothing and regularity estimates}
Throughout this section, for  $r\in[0,1]$ we work with the scale of Hilbert spaces
\begin{equation}\label{Hr-scale}
\mathscr H^r = 	D\bigl(\Le^{(r+1)/2}\bigr)	\times D\bigl(\Le^{r/2}\bigr),
\end{equation}
endowed with the norm
$\|(u,v)\|_{\mathscr H^r}^2
=\|\Le^{(r+1)/2}u\|_2^2
+\|\Le^{r/2}v\|_2^2$.
In particular, the endpoints $\mathscr H^0$ and $\mathscr H^1$
coincide with the energy and strong phase spaces introduced above.
\medskip
The regularity argument is based on separating the exponentially decaying linear component from the component that carries the smoothing effect. For any initial datum $\U_0=(u_0,u_1)\in\mathscr H^0$, we decompose the corresponding solution $\U(t)=(u(t),u_t(t))$ as $u(t)=u^h(t)+u^c(t)$, where
\begin{eqnarray} \label{eq:hom}
\left\lbrace\begin{array}{ll}
u^h_{tt} +\Le(u^h)+ \alpha  u^h_t = 0 &\,\,\text{in}\,\, \Omega \times\mathbb{R}^{+}, \\
u^h=0 &\,\,\text{on}\,\,\partial \Omega\times\mathbb{R}^{+}, \\
u^h(x,0)=u_0(x), \ \  u^h_t(x,0)=u_1(x)&\,\,x \in \Omega,
\end{array}\right.
\end{eqnarray}
and
\begin{eqnarray} \label{eq:comp}
\left\lbrace\begin{array}{ll}
u^c_{tt} +\Le(u^c)+ \alpha  u^c_t + f(u)= \mu_1u -\HH(u)+ b(x) &\,\,\text{in}\,\, \Omega\times\mathbb{R}^{+}, \\
u^c=0 &\,\,\text{on}\,\,\partial \Omega\times\mathbb{R}^{+}, \\
u^c(x,0)=0, \ \  u^c_t(x,0)=0&\,\,x \in \Omega .
\end{array}\right.
\end{eqnarray}
Accordingly, we set $\U^h=(u^h,u_t^h)$ and $\U^c=(u^c,u_t^c)$, so that $\U=\U^h+\U^c$.
 The homogeneous component is governed by the linear semigroup, namely $\U^h(t)=S_h(t)\U_0$. Its exponential decay extends to every space of the scale \eqref{Hr-scale}. Indeed, for
 $\U_0\in\mathscr H^r$, applying $\Le^{r/2}$ to \eqref{eq:hom} and then using \eqref{eq:lin-decay} yields
 \begin{equation}\label{eq:lin-decay-r}
 \|\U^h(t)\|_{\mathscr H^r}
\leqslant	Ce^{-\eta t}\|\U_0\|_{\mathscr H^r},\,\,\, r\in[0,1],\, t\geqslant0.
 \end{equation}
 Thus, all the additional regularity needed below must be obtained from
 the inhomogeneous component $\U^c$.
\begin{Lemma}[Uniform estimates along the Hilbert scale]\label{lem:bootstrap} Assume {\rm(H1)}-{\rm(H7)}. Let $r\in[0,1]$ be fixed and suppose that $\|\U_0\|_{\mathscr{H}^r}\leqslant R_r$. Then there exists a constant $C_r=C_r(R_r)\geqslant0$ such that
$$\|\U(t)\|_{\mathscr{H}^r}\leqslant C_r,\,\,\,\forall \, t\geqslant 0.$$
\end{Lemma}
\begin{proof}
\underline{Case 1: $0\leqslant r \leqslant \tfrac12$.} The result for $r=0$ is Corollary \ref{cor:bounded}. We
proceed by a bootstrap argument on the regularity index. Assume that the estimate has been
established at level $r$ and let the increment $s>0$ satisfy
\begin{equation}\label{eq:step}
0<s\leqslant s(r):=\min\left\{1-r,\ 1-\tfrac{1}{2}(1-2r)(p-1),\ \tfrac{3-p}{2}\right\} .
\end{equation}
Observe that, by {\rm(H6)}, $\theta:=\tfrac{3-p}{2}\in(0,\tfrac12]$ and that $s(r)=\theta$ for
every $r\in[0,\tfrac12]$. Indeed $1-r\geqslant\tfrac12\geqslant\theta$ and, since
$0\leqslant 1-2r\leqslant1$, one has $1-\tfrac12(1-2r)(p-1)\geqslant1-\tfrac{p-1}{2}=\theta$.
In particular, since any positive increment not exceeding $s(r)$ is admissible, every level
$\bar r\in(0,\tfrac12]$ is reached from $r=0$ in at most $\lceil\bar r/\theta\rceil$ steps of equal
size $\bar r/\lceil\bar r/\theta\rceil\leqslant\theta$.
Let $\U_0\in\mathscr H^{r+s}$ with $\|\U_0\|_{\mathscr H^{r+s}}\leqslant R_{r+s}$. Since
$\mathscr H^{r+s}\subset\mathscr H^r$, the bootstrap hypothesis provides a uniform bound in
$\mathscr H^r$, which will be used repeatedly below. Our goal is to show that
\begin{align}\label{eq:B1}
\|\U(t)\|_{\mathscr H^{r+s}}\leqslant C_{r+s},\,\,\,\forall\,t\geqslant0 .
\end{align}
To this end we introduce the higher-order energy
\begin{equation}\label{eq:Ers}
E_{r+s}(t) = \frac{1}{2}\|\U(t)\|_{\mathscr{H}^{r+s}}^2 \!+  \big\langle f(u), \Le^{r+s} u \big\rangle
\!+ \big\langle \HH(u), \Le^{r+s} u \big\rangle \!-\! \big\langle b(x), \Le^{r+s} u \big\rangle
\!-\! \frac{\mu_1}{2}\big\|\Le^{\frac{r+s}{2}}u\big\|_2^2.
\end{equation}
We further introduce the corrector $\chi_{r+s}(t)=\big\langle\Le^{(r+s)/2}u,\Le^{(r+s)/2}u_t\big\rangle$ and the perturbed energy $\Upsilon_{r+s}(t)=E_{r+s}(t)+\gamma \chi_{r+s}(t)$, where
$\gamma>0$ will be chosen later.
\underline{Step 1: equivalence.} By H\"older inequality with conjugate exponents
$\frac{6}{5-2r-2s}$ and $\frac{6}{1+2r+2s}$,
\[
\big|\big\langle f(u),\Le^{r+s}u\big\rangle\big|
\leqslant
\|f(u)\|_{\frac{6}{5-2r-2s}}
\big\|\Le^{r+s}u\big\|_{\frac{6}{1+2r+2s}}.
\]
Since $h\equiv0$, condition \eqref{eq:growth-g} and the mean value
theorem yield
$\|f(u)\|_{6/(5-2r-2s)}
\leqslant C\bigl[1+\|u\|_{6p/(5-2r-2s)}^p\bigr]$. For $0\leqslant r<\frac12$, Remark \ref{rem:fract} gives $D(\Le^{(1+r)/2})\!\subset\! [L^{6/(1-2r)}(\Omega)]^3$ and the restriction $s\leqslant\theta$ in \eqref{eq:step} implies $\frac{6p}{5-2r-2s}\leqslant\frac{6}{1-2r}$, since
$5-2r-2s-p(1-2r)\geqslant2(1-r+rp)>0$. For $r=\frac12$, Remark \ref{rem:fract} gives
$D(\Le^{3/4})\subset[L^q(\Omega)]^3$ for every $q<\infty$ and the restriction is automatic.
Consequently, $\|u\|_{6p/(5-2r-2s)} \leqslant C\|\Le^{(1+r)/2}u\|_2$. On the other hand, Remark \ref{rem:key} gives $\|\Le^{r+s}u\|_{6/(1+2r+2s)}\leqslant C\|\Le^{(1+r+s)/2}u\|_2$.
Combining these estimates and using the bootstrap hypothesis and Young  inequality, we obtain
\begin{align*}
\big|\big\langle f(u),\Le^{r+s}u\big\rangle\big|
\leqslant C\Bigl[1+\big\|\Le^{(1+r)/2}u\big\|_2^p\Bigr] \big\|\Le^{(1+r+s)/2}u\big\|_2
\leqslant \frac{M}{8\lambda_1\underline{\mu}}\|\U(t)\|_{\mathscr H^{r+s}}^2+C_r.
\end{align*}
\noindent Similarly,
$$
\big|\big\langle\HH u,\Le^{r+s}u\big\rangle\big|
\leqslant
\frac{M}{8\lambda_1\underline{\mu}}\|\U(t)\|_{\mathscr H^{r+s}}^2+C_r
 \,\,\,\mbox{and}\,\,\,
\big|\big\langle b,\Le^{r+s}u\big\rangle\big|
\leqslant\frac{M}{8\lambda_1\underline{\mu}}\|\U(t)\|_{\mathscr H^{r+s}}^2+C.
 $$
 Moreover, the last term in \eqref{eq:Ers} satisfies
 $\frac{\mu_1}{2}\|\Le^{(r+s)/2}u\|_2^2
 \leqslant
 \frac{\mu_1}{2\lambda_1\underline{\mu}}
 \|\U(t)\|_{\mathscr H^{r+s}}^2$
 by \eqref{eq:power-Le}. Consequently, there exist constants $C'_1,C'_2>0$, where $C'_2$ may depend on $R_r$, such that
 \begin{equation}\label{eq:Ups-equiv}
 [C'_1-\gamma]\|\U(t)\|_{\mathscr H^{r+s}}^2-C_r
 \leqslant
 \Upsilon_{r+s}(t)
 \leqslant
 [C'_2+\gamma]\|\U(t)\|_{\mathscr H^{r+s}}^2+C_r.
 \end{equation}
 Here, $C'_1=\frac12-\frac{3M+4\mu_1}{8\lambda_1\underline{\mu}}>0$, the positivity being a
 consequence of \eqref{eq:smallness}, since $3M+4\mu_1\leqslant4(2M+\mu_1)<4\lambda_1\underline{\mu}$.

\underline{Step 2: differential inequality.} Testing \eqref{eq:main-abs} by $\Le^{r+s}u_t$ and using the identities
$$\langle f(u),\Le^{r+s}u_t\rangle
=\frac{d}{dt}\langle f(u),\Le^{r+s}u\rangle
-\langle Jf(u)u_t,\Le^{r+s}u\rangle,$$
$$ \langle\HH u,\Le^{r+s}u_t\rangle
=\frac{d}{dt}\langle\HH u,\Le^{r+s}u\rangle
-\langle\HH u_t,\Le^{r+s}u\rangle \,\,\,\mbox{and}\,\,\,\mu_1\langle u,\Le^{r+s}u_t\rangle
=\frac{\mu_1}{2}\frac{d}{dt}
\|\Le^{(r+s)/2}u\|_2^2,$$ together with the fact that $b$ is independent of time, we obtain
\begin{equation*}
E_{r+s}'(t)
=-\alpha\|\Le^{(r+s)/2}u_t\|_2^2+\langle Jf(u)u_t,\Le^{r+s}u\rangle+\langle\HH u_t,\Le^{r+s}u\rangle.
\end{equation*}
The term involving $\mu_1$ is incorporated into
\eqref{eq:Ers} with a favorable sign. Next, using \eqref{eq:growth-g}, H\"older's inequality with $\frac{2-2s}{6}+\frac{3-2r}{6}
+\frac{1+2r+2s}{6}=1$, embeddings $D(\Le^{(1+r)/2})\subset [L^{6(p-1)/(2-2s)}(\Omega)]^3$ and $D(\Le^{r/2})
\subset[L^{6/(3-2r)}(\Omega)]^3$, together with \eqref{eq:key}, we obtain
$$
\begin{aligned}
	\big|\langle Jf(u)u_t,\Le^{r+s}u\rangle\big|
	&\leqslant
	C\Bigl[1+\|\Le^{(1+r)/2}u\|_2^{p-1}\Bigr]
	\|\Le^{r/2}u_t\|_2
	\|\Le^{(1+r+s)/2}u\|_2 \\
	&\leqslant
	\frac{\gamma}{8}
	\|\Le^{(1+r+s)/2}u\|_2^2+C_r.
\end{aligned}
$$
Here, the first embedding follows from the restriction
$s\leqslant1-\frac12(1-2r)(p-1)$ in \eqref{eq:step} and is automatic when $r=\frac12$.
The same argument, with Young's inequality, gives
$\big|\langle\HH u_t,\Le^{r+s}u\rangle\big|
\leqslant\frac{\alpha}{2}
\|\Le^{(r+s)/2}u_t\|_2^2
+\frac{\gamma}{8}
\|\Le^{(1+r+s)/2}u\|_2^2+C_r$.
Consequently,
\begin{equation}\label{eq:De}
E_{r+s}'(t)\leqslant-\frac{\alpha}{2}\|\Le^{(r+s)/2}u_t\|_2^2+\frac{\gamma}{4} \|\Le^{(1+r+s)/2}u\|_2^2+C_r.
\end{equation}
\noindent Furthermore,
\[
\chi_{r+s}'(t)
=
\|\Le^{(r+s)/2}u_t\|_2^2
-\|\Le^{(1+r+s)/2}u\|_2^2
+\big\langle\Le^{r+s}u,
-\alpha u_t-f(u)+\mu_1u-\HH u+b\big\rangle.
\]
The estimates above, the bootstrap hypothesis and Young's inequality
then give
\begin{equation}\label{eq:Dchi}
\chi_{r+s}'(t)\leqslant C\|\Le^{(r+s)/2}u_t\|_2^2
-\frac12\|\Le^{(1+r+s)/2}u\|_2^2+C_r.
\end{equation}
Combining \eqref{eq:De} and \eqref{eq:Dchi}, we obtain
\[
\Upsilon_{r+s}'(t)
\leqslant
-\left[\frac{\alpha}{2}-C\gamma\right]
\|\Le^{(r+s)/2}u_t\|_2^2
-\frac{\gamma}{4}
\|\Le^{(1+r+s)/2}u\|_2^2+C_r.
\]
Choosing $\gamma>0$ sufficiently small so that
$\gamma<C'_1$ and $\alpha/2-C\gamma>0$ and using
\eqref{eq:Ups-equiv}, we find constants $\beta>0$ and
$C_r>0$ such that
$\Upsilon_{r+s}'(t)+\beta\Upsilon_{r+s}(t)\leqslant C_r$.
Gronwall inequality now yields \eqref{eq:B1}.
\medskip
\underline{Case 2: $\tfrac12<r\leqslant1$.} Let $r_0=\min\{r,\tfrac12+\theta\}$. Since
$0<r_0-\tfrac12\leqslant\theta=s\bigl(\tfrac12\bigr)$, the increment $r_0-\tfrac12$ is admissible at
the level $\tfrac12$ and one further application of Case 1 gives
$\|\U(t)\|_{\mathscr H^{r_0}}\leqslant C_{r_0}$ for every $t\geqslant0$. If $r_0=r$, the proof is
complete.
Assume therefore that $r_0=\tfrac12+\theta<r\leqslant1$. Since $r_0>\tfrac12$, we have
$\tfrac32<1+r_0\leqslant2$, so that Corollary \ref{cor:fract}, which applies because
$(1+r_0)/2\leqslant1$, together with \eqref{eq:charact} and the Sobolev embedding give
$$
D\bigl(\Le^{(1+r_0)/2}\bigr)\subset[H^{1+r_0}(\Omega)]^3\subset[L^\infty(\Omega)]^3 ,
$$
whence $\|u(t)\|_{L^\infty}\leqslant C_{r_0}$ uniformly for $t\geqslant0$. This bound removes every
subcritical restriction from the estimates of Steps 1 and 2. Indeed, $\|f(u)\|_2$ and
$\|Jf(u)\|_{L^\infty}$ are now uniformly bounded by \eqref{eq:growth-g}, $\|u_t\|_2$ is controlled by
Corollary \ref{cor:bounded} and $\|\Le^{r'+s'}u\|_2\leqslant C\|\Le^{(1+r'+s')/2}u\|_2$ by
\eqref{eq:power-Le} whenever $r'+s'\leqslant1$. Hence all the terms occurring in $E_{r'+s'}$ and
$\chi_{r'+s'}$ are estimated directly and the only surviving constraint on the increment is
$s'\leqslant1-r'$. Taking $r'=r_0$ and $s'=r-r_0$ we obtain
$\|\U(t)\|_{\mathscr H^{r}}\leqslant C_{r}$ for every $t\geqslant0$, which completes the proof.
\end{proof}
\begin{Lemma}[Smoothing of the inhomogeneous component] \label{lem:comp}
Assume {\rm(H1)}-{\rm(H7)}. Let $r\in[0,1)$ and $0<s\leqslant s(r)$, where $s(r)$ is defined in \eqref{eq:step}. For every $R_r>0$, there exists a constant $R_{r+s}=R_{r+s}(R_r)>0$ such that, whenever $\|\U_0\|_{\mathscr H^r}\leqslant R_r$, the inhomogeneous component $\U^c$ satisfies
$$
 \|\U^c(t)\|_{\mathscr{H}^{r+s}}\leqslant R_{r+s},\,\,\, \forall\, t\geqslant 0 .
 $$
\end{Lemma}
\begin{proof}
The proof follows the strategy of Lemma \ref{lem:bootstrap}. Set $\F(t)=f(u(t))-\mu_1u(t)+\HH u(t)-b$, where $u$ denotes the full	solution and define
$$ E_{r+s}^c(t) =
\frac12\|\U^c(t)\|_{\mathscr H^{r+s}}^2 +\big\langle\F(t),\Le^{r+s}u^c(t)\big\rangle
\,\,\,\mbox{and}\,\,\,
\chi_{r+s}^c(t)=
\big\langle\Le^{(r+s)/2}u^c(t),\Le^{(r+s)/2}u_t^c(t)\big\rangle.
$$
We further set $\Upsilon_{r+s}^c(t)	=E_{r+s}^c(t)+\gamma\chi_{r+s}^c(t)$. Testing \eqref{eq:comp} with $\Le^{r+s}u_t^c$ and using $\langle\F,\Le^{r+s}u_t^c\rangle =\frac{d}{dt}\langle\F,\Le^{r+s}u^c\rangle	-\langle\F_t,\Le^{r+s}u^c\rangle$, where $\F_t=Jf(u)u_t-\mu_1u_t+\HH u_t$, we obtain
$$\big[E_{r+s}^c\big]'(t)
=-\alpha\big\|\Le^{(r+s)/2}u_t^c\big\|_2^2
+\big\langle Jf(u)u_t-\mu_1u_t+\HH u_t,	\Le^{r+s}u^c	\big\rangle.
$$
The terms on the right-hand side are estimated as in Lemma \ref{lem:bootstrap}. The bounds for $u$ and $u_t$ at level $r$ follow from that lemma, while the factor $\Le^{r+s}u^c$ is controlled
by \eqref{eq:key}. When $r>\tfrac12$, the same estimates are used in the simplified form provided by the uniform $L^\infty(\Omega)$ bound established in Case 2 of Lemma \ref{lem:bootstrap}. The additional term involving $\mu_1$ is estimated by \eqref{eq:power-Le} and Corollary \ref{cor:bounded} as
\[
\begin{aligned}
\mu_1\big|\langle u_t,\Le^{r+s}u^c\rangle\big|\leqslant
\mu_1\|u_t\|_2\|\Le^{r+s}u^c\|_2 \leqslant
{C_r}\|\Le^{(1+r+s)/2}u^c\|_2\leqslant
	\varepsilon\|\Le^{(1+r+s)/2}u^c\|_2^2\!+\!C_{\varepsilon,r}.
\end{aligned}
\]
The analogue of \eqref{eq:Ups-equiv} follows in the same way. Indeed,
the level-$r$ bound for $u$ controls the norm $f(u)$ in $L^{6/(5-2r-2s)}(\Omega)$, whereas the higher-order factor now involves
$u^c$. Consequently, the corresponding constants depend only on
$R_r$ and not on time. We note a simplification compared with Lemma \ref{lem:bootstrap}. Since $\F$ does not involve $u^c$, the term 	$\langle\F,\Le^{r+s}u^c\rangle$ is linear in $u^c$ and can therefore 	be absorbed with an arbitrarily small coefficient $\varepsilon>0$. Consequently, the analogue of the lower bound in 	\eqref{eq:Ups-equiv} holds with a constant arbitrarily close to 	$\tfrac12$ and the smallness condition \eqref{eq:smallness} is not
	needed at this stage.
The lower-order bound required above follows from
$\U^c=\U-\U^h$, Lemma~\ref{lem:bootstrap} and
\eqref{eq:lin-decay-r}, which give
$\sup_{t\geqslant0}\|\U^c(t)\|_{\mathscr H^r}\leqslant C_r$.
Since $\U^c(0)=0$, we also have $\Upsilon_{r+s}^c(0)=0$.
The resulting differential inequality and Gronwall's lemma therefore
yield a uniform bound for $\Upsilon_{r+s}^c(t)$. Combining this bound
with the analogue of \eqref{eq:Ups-equiv}, we conclude that
$\|\U^c(t)\|_{\mathscr H^{r+s}}\leqslant R_{r+s}(R_r)$ for every
$t\geqslant0$.
\end{proof}
\subsubsection{An exponentially attracting set in $\mathscr H^1$}
We first construct a bounded subset of $\mathscr H^1$ that
exponentially attracts every bounded subset of $\mathscr H^0$.
Consider the closed neighbourhood of the global attractor
$$\mathbf{N}_{\mathbf{A}}=\{\U\in\mathscr H^0\,|\,
\dist_{\mathscr H^0}(\U,\A)\leqslant1\}.$$
Since $\A$ is compact, $\mathbf{N}_{\mathbf{A}}$ is bounded in $\mathscr H^0$.
Moreover, since $\A$ attracts every bounded subset of
$\mathscr H^0$, the set $\mathbf{N}_{\mathbf{A}}$ is absorbing.
\begin{Theorem}\label{thm:B1}
There exist a closed bounded set $\mathbf{B}_1\subset\mathscr H^1$ and constants $N,\nu>0$ such that
\begin{equation}\label{eq:B1-attr}
\dist_{\mathscr H^0}\bigl(S(t)\mathbf{N}_{\mathbf{A}},\mathbf{B}_{1}\bigr)
\leqslant Ne^{-\nu t},\,\,\,\forall\, t\geqslant0.
\end{equation}
\end{Theorem}
\begin{proof}
Set $\theta=(3-p)/2>0$. From \eqref{eq:step}, one readily checks that $s(r)=\min\{1-r,\theta\}$ for every $r\in[0,1)$. Define $k_0=0$ and $k_{j+1}=k_j+s(k_j)$. Then $k_j=\min\{j\theta,1\}$ and hence $k_n=1$ for $n=\lceil1/\theta\rceil=\lceil2/(3-p)\rceil$. We thus obtain a finite chain $0=k_0<k_1<\cdots<k_n=1$.
Let $R_{k_0}$ be a bound for $\mathbf{N}_{\mathbf{A}}$ in $\mathscr H^0$. Applying Lemma \ref{lem:comp} recursively along this chain, we obtain constants $R_{k_j}=R_{k_j}(R_{k_{j-1}})>0$ and define $\mathbf{B}_{k_j}=\{\U\in\mathscr H^{k_j}\,|\, \|\U\|_{\mathscr H^{k_j}}\leqslant R_{k_j}\}$ for 	$j=1,\ldots,n$.	Fix $\U_0\in \mathbf{B}_{k_{j-1}}$ and write	$S(t)\U_0=\U^h(t)+\U^c(t)$. By Lemma \ref{lem:comp},	$\U^c(t)\in \mathbf{B}_{k_j}$ for every $t\geqslant0$. Therefore,
$$
\dist_{\mathscr H^0}\bigl(S(t)\U_0,\mathbf{B}_{k_j}\bigr)
	\leqslant	\|\U^h(t)\|_{\mathscr H^0}
	\leqslant	Ce^{-\eta t}\|\U_0\|_{\mathscr H^0}.
$$
Since $\mathscr H^{k_{j-1}}$ is continuously embedded in $\mathscr H^0$, taking the supremum over $\U_0\in \mathbf{B}_{k_{j-1}}$ gives
$$
\dist_{\mathscr H^0}
\bigl(S(t)\mathbf{B}_{k_{j-1}},\mathbf{B}_{k_j}\bigr)
\leqslant C_j e^{-\eta t}\, \,\, (j=1,\ldots,n).
$$
All the sets $\mathbf{B}_{k_j}$ are bounded in $\mathscr H^0$. Hence, Theorem \ref{thm:global}$(3)$ provides the Lipschitz estimate required by Lemma~\ref{lem:trans}. Applying that lemma successively along the finite chain, we find constants $N,\nu>0$ such that $\dist_{\mathscr H^0}(S(t)\mathbf{N}_{\mathbf{A}},\mathbf{B}_{k_n}) \leqslant Ne^{-\nu t}$ for every $t\geqslant0$. Since $k_n=1$, the set $\mathbf{B}_{k_n}$ is a closed bounded subset of $\mathscr H^1$. Taking $\mathbf{B}_1=\mathbf{B}_{k_n}$ proves \eqref{eq:B1-attr}.
\end{proof}
\begin{Lemma}[Uniform bound for the time derivative]\label{lem:Zt}
There exists a constant $C_{\mathbf{B}_1}>0$ such that every solution with $\U_0\in \mathbf{B}_1$ satisfies
$$\|\U_t(t)\|_{\mathscr H^0} \leqslant C_{\mathbf{B}_1},\,\,\,\forall \, t \geqslant 0.$$
\end{Lemma}
\begin{proof}
Since $\mathbf{B}_1$ is bounded in $\mathscr H^1=D(\Ee)$, Theorem \ref{thm:global}$(2)$ implies the regularity $\U\in C^1([0,\infty);\mathscr H^0)$ and $\U_t(t)=\Ee\U(t)+\Ff(\U(t))$ in $\mathscr H^0$. By \eqref{eq:EF} and Lemma \ref{lem:lipschitz}, $\|\Ee\U\|_{\mathscr H^0}\leqslant C\|\U\|_{\mathscr H^1}$ and	$\|\Ff(\U)\|_{\mathscr H^0}\leqslant	C[1+\|\U\|_{\mathscr H^0}]^3$.
Moreover, Lemma \ref{lem:bootstrap} with $r=1$ yields $\|\U(t)\|_{\mathscr H^1}\leqslant C_{\mathbf{B}_1}$. Since $\mathscr H^1\subset\mathscr H^0$, the preceding estimates give the assertion.
\end{proof}

\subsubsection{Verification of the hypotheses of Theorem \ref{thm:exp-crit}}
We first construct the bounded positively invariant set on which the
fixed-time map will be considered. After fixing this set, we verify the
Lipschitz property and the contractive-smoothing decomposition required
by Theorem \ref{thm:exp-crit}.
\begin{Lemma}[A compact exponentially attracting set]\label{lem:A}
	Define $\displaystyle \mathbf{D}= \overline{\bigcup_{\tau\geqslant0}S(\tau)\mathbf{B}_1}^{\,\mathscr H^0}$. Then the following properties hold.
	\begin{enumerate}
		\item[$1.$] The set $\mathbf{D}$ is compact in $\mathscr H^0$ and bounded in $\mathscr H^1$.

		\item[$2.$] The set $\mathbf{D}$ is positively invariant.

		\item[$3.$] There exist constants $N_1,\nu_1>0$ such that
		\begin{equation}\label{eq:A-attr}
			\dist_{\mathscr H^0}\bigl(S(t)\mathbf{N}_{\mathbf{A}},\mathbf{D}\bigr)
			\leqslant N_1e^{-\nu_1t},\,\,\, t\geqslant0.
		\end{equation}

		\item[$4.$] There exists a constant $C_{\mathbf{D}}>0$ such that every trajectory $\U(t)=S(t)\U_0$ with $\U_0\in\mathbf{D}$ satisfies
		\begin{equation}\label{eq:Zt-A}
			\|\U_t(t)\|_{\mathscr H^0}\leqslant C_{\mathbf{D}},\,\,\,\forall\,t\geqslant 0.
		\end{equation}
	\end{enumerate}
\end{Lemma}
\begin{proof}
	Since $\mathbf{B}_1$ is bounded in $\mathscr H^1$, Lemma \ref{lem:bootstrap} with $r=1$ yields a constant $R_{1}=R_{1}(\mathbf{B}_1)>0$ such that $\|S(t)\U_0\|_{\mathscr H^1} \leqslant R_{1}$, for all	$\U_0\in \mathbf{B}_1$ and $t\geqslant 0$. This bound is preserved when the closure is taken in $\mathscr H^0$.	Indeed, if $\U_m\to\U$ in $\mathscr H^0$ and	$\|\U_m\|_{\mathscr H^1}\leqslant R_{1}$, then a subsequence converges	weakly in $\mathscr H^1$. Its weak limit must coincide with $\U$ and the weak lower semicontinuity of the norm gives	$\|\U\|_{\mathscr H^1}\leqslant R_{1}$. Hence, $\mathbf{D}$ is bounded in $\mathscr H^1$. Since $\mathscr H^1$ is compactly embedded in 	$\mathscr H^0$ and $\mathbf{D}$ is closed in $\mathscr H^0$, the set $\mathbf{D}$ is compact in $\mathscr H^0$.

	For every $t\geqslant0$, the continuity of $S(t)$ in $\mathscr H^0$ and the semigroup property give
	$$
	S(t)\mathbf{D}
	\subseteq
	\overline{\bigcup_{\tau\geqslant0}S(t+\tau)\mathbf{B}_1}^{\,\mathscr H^0}
	\subseteq
	\overline{\bigcup_{\tau\geqslant0}S(\tau)\mathbf{B}_1}^{\,\mathscr H^0}
	=\mathbf{D}.
	$$
	Thus, $\mathbf{D}$ is positively invariant.

	Since $S(0)=I$, the definition of $\mathbf{D}$ gives $\mathbf{B}_1\subset \mathbf{D}$. Therefore, Theorem \ref{thm:B1} immediately yields $\dist_{\mathscr H^0}\bigl(S(t)\mathbf{N}_{\mathbf{A}},\mathbf{D}\bigr)\leqslant
	\dist_{\mathscr H^0}\bigl(S(t)\mathbf{N}_{\mathbf{A}},\mathbf{B}_1\bigr)
	\leqslant
	Ne^{-\nu t}$ for all $t\geqslant0$. Hence, \eqref{eq:A-attr} holds with $N_1=N$ and $\nu_1=\nu$.

	Finally, the boundedness of $\mathbf{D}$ in $\mathscr H^1$ and its positive
	invariance imply $\|S(t)\U_0\|_{\mathscr H^1}\leqslant R_{\mathbf{D}}$ for every
	$\U_0\in {\mathbf{D}}$ and $t\geqslant0$,where
	$R_{\mathbf D}$ denotes the corresponding radius. Repeating the argument of
	Lemma \ref{lem:Zt}, we obtain
	\[
	\|\U_t(t)\|_{\mathscr H^0}
	\leqslant
	\|\Ee\U(t)\|_{\mathscr H^0}
	+\|\Ff(\U(t))\|_{\mathscr H^0}
	\leqslant
	CR_{\mathbf{D}}+C(1+R_{\mathbf{D}})^3,
	\]
	uniformly for $\U_0\in {\mathbf{D}}$ and $t\geqslant0$. This proves
	\eqref{eq:Zt-A}.
\end{proof}
\begin{Lemma}[Spatial regularity of the nonlinear difference] \label{lem:Df}
Let $\U_0^1,\U_0^2\in \mathbf{D}$ and let $\U^j(t)=(u^j(t),u_t^j(t))=S(t)\U_0^j$, $j=1,2$, be the corresponding strong solutions of \eqref{eq:main-abs}. Define $\widetilde f_j=f_j(u^2)-f_j(u^1)$. Then, for every $t_*>0$, there exists a constant $C_{t_*,\mathbf{D}}>0$ such that
$$
\left\|	\frac{\partial\widetilde f_j}{\partial x_k}(t)	\right\|_2 \leqslant C_{t_*,\mathbf{D}}\|\U_0^1-\U_0^2\|_{\mathscr H^0}, \quad \forall \, t\in[0,t_*], \, j,k\in\{1,2,3\}.
$$
\end{Lemma}
\begin{proof}
Set $w=u^2-u^1$ and $z_\sigma=u^1+\sigma w$. By the mean value formula, $\widetilde f_j=\int_0^1\nabla f_j(z_\sigma)\cdot w\,d\sigma$. Differentiating with respect to $x_k$ and using the chain rule, we obtain
\begin{align*}
\frac{\partial\widetilde f_j}{\partial x_k}
&=\int_0^1 \bigl[D^2f_j(z_\sigma)\partial_{x_k}z_\sigma \bigr]\cdot w\,d\sigma + \int_0^1 \nabla f_j(z_\sigma)\cdot\partial_{x_k}w\,d\sigma.
\end{align*}
Hence, by \eqref{eq:growth-g} and {\rm(H7)},
\begin{align*}
\left|\frac{\partial\widetilde f_j}{\partial x_k} \right|
&\leqslant
C\left[	1+\sum_{i=1}^3	\bigl(|u_i^1|^{p-2}+|u_i^2|^{p-2}\bigr)\right]
\Bigl[\bigl(|\nabla u^1|+|\nabla u^2|\bigr)|w|+|\nabla w| \Bigl].
\end{align*}
Since $\mathbf{D}$ is bounded in $\mathscr H^1$, Lemma~\ref{lem:bootstrap} with $r=1$ gives $\sup_{t\in[0,t_*]}\|u^j(t)\|_{H^2}\leqslant C_{\mathbf{D}}$. The embeddings $H^2(\Omega)\subset L^\infty(\Omega)$ and $H^2(\Omega)\subset W^{1,3}(\Omega)$, together with H\"older's inequality, therefore yield
$$
\left\|	\frac{\partial\widetilde f_j}{\partial x_k}(t)\right\|_2
\leqslant C_{\mathbf{D}}\Bigl[	\|\nabla u^1(t)\|_3+\|\nabla u^2(t)\|_3\Bigr]\|w(t)\|_6+C_{\mathbf{D}}\|\nabla w(t)\|_2
\leqslant C_{\mathbf{D}}\|w(t)\|_{H^1}.
$$
Finally, the equivalence $D(\Le^{1/2})=[H_0^1(\Omega)]^3$ and \eqref{eq:cont-dep} imply 	$\|w(t)\|_{H^1}	\leqslant C_{t_*,\mathbf{D}}\|\U_0^1-\U_0^2\|_{\mathscr H^0}$ for $t\in[0,t_*]$. Combining the last two estimates completes the	proof.
\end{proof}
\begin{Lemma}[Verification of the exponential attractor criterion]	\label{lem:criterion} The following properties hold.
\begin{enumerate}
\item[$1.$] For every $T>0$, the map $[0,T]\times \mathbf{D}\ni(t,\U_0)\mapsto S(t)\U_0\in \mathbf{D}$ is Lipschitz continuous.
\item[$2.$] There exist $t_*>0$ and $C_*>0$ such that the decomposition $S(t_*)|_{\mathbf{D}}=S_h(t_*)|_{\mathbf{D}}+S_c(t_*)|_{\mathbf{D}}$ satisfies
$$
\|S_h(t_*)\U_0^1-S_h(t_*)\U_0^2\|_{\mathscr H^0}
\leqslant
\frac18\|\U_0^1-\U_0^2\|_{\mathscr H^0}
$$
and
$$
\|S_c(t_*)\U_0^1-S_c(t_*)\U_0^2\|_{\mathscr H^1}
\leqslant
C_*\|\U_0^1-\U_0^2\|_{\mathscr H^0},
$$
for all $\U_0^1,\U_0^2\in \mathbf{D}$.
\end{enumerate}
Consequently, the assumptions of Theorem $\ref{thm:exp-crit}$ are satisfied with $H_0=\mathscr H^0$, $H_1=\mathscr H^1$, $S_0=S_h(t_*)|_{\mathbf{D}}$ and $S_1=S_c(t_*)|_{\mathbf{D}}$.
\end{Lemma}
\begin{proof}
Let $T>0$, $\U_0^1,\U_0^2\in \mathbf{D}$ and $t_1,t_2\in[0,T]$. By the triangle inequality,
$$
\|S(t_1)\U_0^1-S(t_2)\U_0^2\|_{\mathscr H^0}
\leqslant \|S(t_1)\U_0^1-S(t_1)\U_0^2\|_{\mathscr H^0}+	\|S(t_1)\U_0^2-S(t_2)\U_0^2\|_{\mathscr H^0}.
$$
The continuous dependence estimate \eqref{eq:cont-dep} gives
$$\|S(t_1)\U_0^1-S(t_1)\U_0^2\|_{\mathscr H^0}\leqslant	C_{T,\mathbf{D}}\|\U_0^1-\U_0^2\|_{\mathscr H^0}.
$$
To estimate the second term, set $\tau_1=\min\{t_1,t_2\}$ and	$\tau_2=\max\{t_1,t_2\}$. Since the trajectories starting in ${\mathbf{D}}$ are strong, the fundamental theorem of calculus and \eqref{eq:Zt-A} yield
\begin{align*}
\|S(t_1)\U_0^2-S(t_2)\U_0^2\|_{\mathscr H^0}\leqslant\int_{\tau_1}^{\tau_2}\left\|\frac{d}{ds}S(s)\U_0^2\right\|_{\mathscr H^0}\,ds
\leqslant	C_{\mathbf{D}}|t_1-t_2|.
\end{align*}
Consequently, the map is Lipschitz with constant $L_T=\max\{C_{T,\mathbf{D}},C_{\mathbf{D}}\}$. This proves the first assertion for every $T>0$.
We now choose the time required for the decomposition. By Remark \ref{rem:homog}, there exist $C,\eta>0$ such that	$\|S_h(t)\U_0\|_{\mathscr H^0}	\leqslant Ce^{-\eta t}\|\U_0\|_{\mathscr H^0}$.	Choose $t_*>0$ so that $Ce^{-\eta t_*}\leqslant1/8$ and define $S_0\U_0:=S_h(t_*)\U_0$, $S_1\U_0:=S_c(t_*)\U_0$ with $\U_0\in {\mathbf{D}}$. Then $S(t_*)|_{\mathbf{D}}=S_0+S_1$. Since $S_h(t_*)$ is linear,
\begin{equation}\label{eq:contractive}
\|S_0\U_0^1-S_0\U_0^2\|_{\mathscr H^0}
\leqslant\frac18\|\U_0^1-\U_0^2\|_{\mathscr H^0}.
\end{equation}
It remains to establish the smoothing estimate for $S_1$. Set $(v,v_t)=S_c(t)\U_0^1-S_c(t)\U_0^2$ and $\widetilde u=u^1-u^2$. Then
\begin{equation}\label{eq:v-eq}\begin{cases}
v_{tt}+\Le v+\alpha v_t=\F,&\text{in }\Omega\times\mathbb R^+,\\
v=0,&\text{on }\partial\Omega\times\mathbb R^+,\\
v(0)=0,\,\,\, v_t(0)=0,&\text{in }\Omega,
\end{cases}
\end{equation}
where $\F=f(u^2)-f(u^1)+\mu_1\widetilde u-\HH\widetilde u$. By Lemma~\ref{lem:Df},
$f(u^2)-f(u^1)\in[H^1(\Omega)]^3$. Moreover,
$\widetilde u\in[H_0^1(\Omega)]^3$ and the regularity
$\mu\in C^3(\overline{\Omega})$ ensures that
$\HH\widetilde u\in[H^1(\Omega)]^3$. All the terms in $\F$ have zero
trace, since $u^1=u^2=0$ on $\partial\Omega$ and hence
$f(u^2)-f(u^1)=f(0)-f(0)=0$ there. Therefore,
$\F\in[H_0^1(\Omega)]^3=D(\Le^{1/2})$ and Lemma \ref{lem:Df}, \eqref{eq:cont-dep}  and the regularity of the coefficients imply that $\|\Le^{1/2}\F(t)\|_2^2\leqslant C_{t_*,\mathbf{D}}\|\U_0^1-\U_0^2\|_{\mathscr H^0}^2$ for all $t\in[0,t_*]$.
Before testing \eqref{eq:v-eq} with $\Le v_t$, we justify the
computation. Since $\U_0^j\in\mathbf D$, Lemma
\ref{lem:bootstrap} with $r=1$ gives
$u^j\in C([0,t_*];D(\Le))\cap
C^1([0,t_*];D(\Le^{1/2}))$. Assumption {\rm(H7)} and the embedding
$H^2(\Omega)\subset L^\infty(\Omega)$ imply that the Nemytskii
operator $v\mapsto f(v)$ is continuous from $D(\Le)$ into
$[H^1(\Omega)]^3$. Consequently,
$\F\in C([0,t_*];D(\Le^{1/2}))$. Choose a sequence
$\F_m\in C^1([0,t_*];D(\Le^{1/2}))$ converging to $\F$ in
$C([0,t_*];D(\Le^{1/2}))$. Since $v(0)=v_t(0)=0$, the classical
theory for the corresponding inhomogeneous linear problems provides
solutions $v_m$ in $D(\Le^{3/2})\times D(\Le)$ for which the following
computation is legitimate. The resulting estimate passes to the limit
because it involves only $\|\Le^{1/2}\F\|_2$ and the
$\mathscr H^1$-norm of $(v,v_t)$, both of which are continuous under
this approximation.
\noindent Testing \eqref{eq:v-eq} with $\Le v_t$ gives
$$\frac{d}{dt} \left[\frac12\|(v,v_t)\|_{\mathscr H^1}^2\right] +\alpha\|\Le^{1/2}v_t\|_2^2 = \langle\Le^{1/2}\F,\Le^{1/2}v_t\rangle.
$$
The Cauchy-Schwarz and Young inequalities imply
$$
\frac{d}{dt}\left[\frac12\|(v,v_t)\|_{\mathscr H^1}^2\right] +\frac{\alpha}{2}\|\Le^{1/2}v_t\|_2^2 \leqslant C\|\Le^{1/2}\F\|_2^2.
$$
Since $(v(0),v_t(0))=(0,0)$, integration over $[0,t_*]$ yields
\begin{align}\label{eq:smoothing}
\|S_1\U_0^1-S_1\U_0^2\|_{\mathscr H^1} \leqslant C_*\|\U_0^1-\U_0^2\|_{\mathscr H^0},
\end{align}
where $C_*=\sqrt{2C C_{t_*,\mathbf{D}}t_*}$.	Together with \eqref{eq:contractive}, this verifies the second assertion. Taking $T=t_*$ in the first assertion completes the verification of Theorem \ref{thm:exp-crit}.
\end{proof}
\begin{Remark}\label{rem-set} The proofs of Lemmas \ref{lem:Df} and \ref{lem:criterion} do not require the set $\mathbf{D}$ itself to be positively invariant. The only property used is the uniform boundedness in $\mathscr H^1$ of the forward trajectories issuing from $\mathbf{D}$. Consequently, both lemmas remain valid for any set whose forward	trajectories satisfy such a bound. This observation will be used in	Section \ref{sec:robust}, where the relevant set is bounded in	$\mathscr H^1$ and is invariant only for sufficiently large times.
\end{Remark}
\subsubsection{Completion of the proof of Theorem~\ref{thm:exp}}
\begin{proof} Let $\mathbf D$ be the bounded positively invariant set constructed in Lemma \ref{lem:A}. By Lemma \ref{lem:criterion}, the hypotheses of Theorem \ref{thm:exp-crit} are satisfied on $\mathbf D$. We therefore obtain a compact positively invariant set $\A^{\exp}\subset\mathbf D$ such that $\dim_F^{\mathscr H^0}(\A^{\exp})\leqslant1+\log N^*/\log2$ and
\begin{equation}\label{eq:exp-attr-A}
\dist_{\mathscr H^0}\bigl(S(t)\mathbf D,\A^{\exp}\bigr)
\leqslant Ne^{-\nu t},\,\,\, t\geqslant0.
\end{equation}
It remains to show that $\A^{\exp}$ attracts every bounded subset of $\mathscr H^0$ exponentially. Combining \eqref{eq:A-attr} and	\eqref{eq:exp-attr-A} through Lemma \ref{lem:trans}, we obtain	constants $N_2,\nu_2>0$ such that
\begin{equation}\label{eq:B0-exp}
\dist_{\mathscr H^0}\bigl(S(t)\mathbf{N}_{\mathbf{A}},\A^{\exp}\bigr)\leqslant N_2 e^{-\nu_2 t},\,\,\,  t\geqslant0,
\end{equation}
where $\nu_2$ may be chosen smaller than both $\nu$ and $\nu_1$.
Let $D\subset\mathscr H^0$ be bounded and set
$r(D)=\sup_{\U\in D}\|\U\|_{\mathscr H^0}$. Since
$\mathbf N_{\mathbf A}$ is absorbing, its entrance time can be chosen
as a nondecreasing function of $r(D)$. Denote this time by $t_D$, so
that $S(t)D\subseteq\mathbf N_{\mathbf A}$ for every $t\geqslant t_D$.
For such $t_D$, the semigroup property and \eqref{eq:B0-exp} give
\begin{align*}
\dist_{\mathscr H^0}\bigl(S(t)D,\A^{\exp}\bigr)
&=\dist_{\mathscr H^0}\bigl(S(t-t_D)S(t_D)D,\A^{\exp}\bigr)\\
&\leqslant \dist_{\mathscr H^0} \bigl(S(t-t_D)\mathbf{N}_{\mathbf{A}},\A^{\exp}\bigr)\\
&\leqslant N_2e^{-\nu_2(t-t_D)} =	N_2 e^{\nu_2 t_D}e^{-\nu_2 t}.
\end{align*}
Corollary \ref{cor:bounded} gives $\|S(t)\U\|_{\mathscr H^0}\leqslant R(r_D)$ for every $\U\in D$ and $t\geqslant0$, where $r_D=\sup_{\U\in D}\|\U\|_{\mathscr H^0}$ and $R$ is nondecreasing. Setting $R_{\exp} =\sup_{\U\in\A^{\exp}}\|\U\|_{\mathscr H^0}$, we obtain, for
$0\leqslant t\leqslant t_D$
\begin{align*}
\dist_{\mathscr H^0}\bigl(S(t)D,\A^{\exp}\bigr) \leqslant R(r(D))+R_{\exp}
\leqslant\bigl[R(r(D))+{R_{\exp}}\bigr]e^{\nu_2 {t_D}}e^{-\nu_2 t}.
\end{align*}
Consequently,
$$
\dist_{\mathscr H^0}\bigl(S(t)D,\A^{\exp}\bigr)\leqslant C(r_D)e^{-\nu_2t},\,\,\, t\geqslant0,
$$
where $C(r)= e^{\nu_2T(r)}\max\bigl\{N_2,R(r)+R_{\exp}\bigr\}$. Here $T(r)$ is the nondecreasing entrance-time function introduced above, so that $t_D=T(r_D)$. Since both $T$ and $R$ are
nondecreasing, the function $C$ is also nondecreasing. Therefore,
$\A^{\exp}$ attracts every bounded subset of $\mathscr H^0$
exponentially. Together with its compactness, positive invariance and
finite fractal dimension, this proves that $\A^{\exp}$ is an
exponential attractor for $\{S(t)\}_{t\geqslant0}$ in
$\mathscr H^0$.
Finally, since $\mathbf A$ is invariant, $S(t)\mathbf A=\mathbf A$ for every $t\geqslant0$ and
therefore
$$
\dist_{\mathscr H^0}\bigl(\mathbf A,\A^{\exp}\bigr)
=\dist_{\mathscr H^0}\bigl(S(t)\mathbf A,\A^{\exp}\bigr)
\leqslant C\bigl(r(\mathbf A)\bigr)e^{-\nu_2 t}\xrightarrow[t\to\infty]{}0 .
$$
Hence $\dist_{\mathscr H^0}(\mathbf A,\A^{\exp})=0$ and $\A^{\exp}$ being closed,
$\mathbf A\subseteq\A^{\exp}$.
\end{proof}
\begin{Remark}\label{rem:subcritical} The construction of $\A^{\exp}$ relies essentially on the subcritical assumptions {\rm(H6)}-{\rm(H7)}, including $h\equiv0$ and $p<3$. Indeed, the smoothing estimate \eqref{eq:smoothing} uses the second-order growth condition {\rm(H7)} through Lemma \ref{lem:Df} and is not available for the critical component allowed in Theorem \ref{thm:attractor}. The existence of an exponential attractor in the critical regime therefore remains open for the present system. A similar subcritical restriction appears in \cite{PS}. Finally, since $\mathbf A\subseteq\A^{\exp}$, Theorem \ref{thm:exp} also provides an alternative proof that $\mathbf A$ has finite fractal dimension in the subcritical case.
\end{Remark}
\section{Robustness under degeneration of the wave-speed gap}
\label{sec:robust}
In this section, we study the stability of the long-time dynamics under
degeneration of the elastic wave-speed gap. More precisely, the first
Lam\'e coefficient is replaced by a family
$\{\lambda_\delta\}_{\delta\in[0,1]}$ converging to a limiting
coefficient $\lambda_0$, while the shear modulus $\mu$ remains fixed.
After introducing the parameter-dependent coefficients, we recall the
abstract robustness criterion of Theorem \ref{thm:robust-abstract} and establish the uniform coercivity, elliptic regularity, dissipativity, decay and smoothing estimates needed for its
application. These results lead to the robust family of exponential attractors in Theorem \ref{thm:robust} and to the upper semicontinuity of the global attractors in Theorem
\ref{thm:usc}.

\subsection{The parameter-dependent family}
To carry out this analysis, assumptions {\rm(H1)}, {\rm(H3)}-{\rm(H5)} and {\rm(H7)} remain in force, with the sole proviso that the smallness condition \eqref{eq:smallness} in {\rm(H4)} is replaced by its uniform counterpart stated at the end of {\rm(H8)}. The conditions from {\rm(H2)} and {\rm(H6)} required in this section are stated explicitly, while the dependence on $\delta$ is specified in {\rm(H8)}. More precisely, we assume
$$
\rho\equiv1,\,\,\, \alpha>0,\,\,\, h\equiv0,\,\,\, f=g \,\,\,\mbox{and}\,\,\, 2\leqslant p<3.
$$
The shear modulus $\mu\in C^3(\overline\Omega)$ satisfies $\mu(x)\geqslant\underline\mu>0$ for all $x\in\overline\Omega$. We retain the notations $\displaystyle \overline\mu=\max_{x\in\overline\Omega}\mu(x)$ and $\displaystyle \mu_2=
\frac12 \max_{\substack{x\in\overline\Omega\\1\leqslant i,j\leqslant3}}
|\partial_{ij}^2\mu(x)|$.
The dependence on $\delta$ is described by the following hypothesis.
\begin{itemize}
\item[{\rm(H8)}] $\{\lambda_\delta\}_{\delta\in[0,1]} \subset C^{0,1}(\overline\Omega)$ is a family of first Lam\'e coefficients satisfying
\begin{align}
&\Lambda=\sup_{\delta\in[0,1]}\|\lambda_\delta\|_{C^{0,1}(\overline\Omega)}<\infty,\label{eq:H8b}\\[1mm]
&\|\lambda_\delta-\lambda_0\|_{L^\infty(\Omega)}\leqslant\delta,\,\,\,\forall\,\delta\in[0,1],\label{eq:H8c}\\[1mm]
&\lambda_\delta(x)+\mu(x)\geqslant0,\,\,\, \forall\,x\in\overline\Omega,\ \forall\,\delta\in[0,1],\label{eq:H8a}\\[1mm]
&\mu_1^{*}=\sup_{\delta\in[0,1]}\ \sup_{x\in\overline\Omega}
\frac{\bigl|\nabla\mu(x)\bigr|^{2}}{\lambda_\delta(x)+\mu(x)}\ <\ \infty .
\label{eq:H8d}
\end{align}
\end{itemize}
Whenever both the numerator and the denominator in \eqref{eq:H8d}
vanish, the quotient is defined to be zero. Moreover, we assume $2M+6\mu_2+\mu_1^{*}<\underline\mu\lambda_1$.

\subsection{Physical interpretation of the limiting regime} \label{sec:physics}
Assumption {\rm(H8)} is formulated in terms of the Lam\'e coefficients $\lambda_\delta$ and $\mu$, while its physical content is more transparent in terms of the bulk modulus
$K_\delta=\lambda_\delta+\tfrac23\mu$, the shear modulus $\mu$ and the corresponding wave speeds. This subsection interprets the parameter-dependent assumptions from this perspective and complements the discussion in Subsection \ref{sub:physical}.

We denote by
\begin{equation}\label{eq:omega0}
	\omega_0
	=
	\bigl\{
	x\in\overline\Omega
	\,|\,
	\lambda_0(x)+\mu(x)=0
	\bigr\}
\end{equation}
the \emph{transition zone} of the limit medium. By \eqref{eq:speeds}, this is precisely the set where the pressure- and shear-wave speeds coincide at $\delta=0$. Assumption {\rm(H8)} imposes no restriction on the geometry or measure of $\omega_0$.

\medskip
	
\underline{1. The hypotheses in terms of the bulk and shear moduli.}
Let
$K_\delta=\lambda_\delta+\tfrac23\mu$ denote the bulk modulus. The combinations of the Lam\'e coefficients appearing in the analysis satisfy
\begin{equation}\label{eq:moduli}
\lambda_\delta+\mu=K_\delta+\tfrac13\mu,\,\,\,
\lambda_\delta+2\mu=K_\delta+\tfrac43\mu	\,\,\, \mbox{and}\,\,\,
	\left(\frac{c_{P,\delta}}{c_S}\right)^2
	=
	\frac{K_\delta}{\mu}+\frac43.
\end{equation}
Thus, the ratio between the pressure- and shear-wave speeds is determined by the dimensionless quantity $K_\delta/\mu$. The three conditions in Remark \ref{rem:speeds} can be written as
\begin{center}
	\begin{tabular}{lll}
		\hline
		condition & equivalent form & consequence\\
		\hline
		$3\lambda_\delta+2\mu>0$
		&\,\,\, $K_\delta/\mu>0$
		& coercivity of $a$\\
		$\lambda_\delta+\mu>0$
		& \,\,\, $K_\delta/\mu>-\tfrac13$
		& coercivity of the modified form\\
		$\lambda_\delta+2\mu>0$
		& \,\,\, $K_\delta/\mu>-\tfrac43$
		& Legendre-Hadamard condition\\
		\hline
	\end{tabular}
\end{center}
The first condition gives direct coercivity of the original elastic form by Proposition \ref{prop:unperturbed}. The third concerns only the principal symbol and guarantees hyperbolicity. For a variable shear modulus, it does not by itself imply coercivity. The intermediate condition is the one used to construct the modified coercive form in Subsection \ref{sub:modified}. The associated lower-order terms are subsequently controlled through the smallness condition.

In terms of the bulk modulus, \eqref{eq:H8a} becomes
$$
K_\delta(x)\geqslant-\frac13\mu(x),
\,\,\,
x\in\overline\Omega,
\,\,\,
\delta\in[0,1].
$$
Moreover,
$$
\omega_0
=
\left\{
x\in\overline\Omega
\,|\,
K_0(x)=-\frac13\mu(x)
\right\},
$$
and \eqref{eq:H8d} reads
$$
\mu_1^*
=
\sup_{\delta\in[0,1]}
\sup_{x\in\overline\Omega}
\frac{|\nabla\mu(x)|^2}
{K_\delta(x)+\tfrac13\mu(x)}
<
\infty.
$$
Hence, {\rm(H8)} allows the bulk modulus to enter the range
$-\mu/3\leqslant K_\delta\leqslant0$, which is excluded by the classical positivity condition. In the limit medium, the endpoint
$K_0=-\mu/3$ corresponds precisely to the coincidence
$c_{P,0}=c_S$.	
	
\medskip

\underline{ 2. Why the limiting medium may remain stable.} As discussed in Subsection \ref{sub:physical}, a negative bulk modulus is incompatible with the stability of a free homogeneous isotropic body but may be stabilized by a constraint. A quantitative example is provided by the criterion of Lakes and Drugan \cite{lakesdrugan}. For a cylindrical inclusion with Lam\'e coefficients
$(\lambda^{\rm i},\mu^{\rm i})$ embedded in a matrix with shear modulus
$\mu^{\rm m}$, their condition reads
$$
\frac{\mu^{\rm m}}{\mu^{\rm i}}
>
-1-\frac{\lambda^{\rm i}}{\mu^{\rm i}},
$$
or equivalently
\begin{equation}\label{eq:LD}
\lambda^{\rm i}+\mu^{\rm i}+\mu^{\rm m}>0.
\end{equation}
At the algebraic level, \eqref{eq:H8a} is more restrictive, since it requires
$\lambda_\delta+\mu\geqslant0$ without an additional stabilizing contribution from a surrounding matrix. In the present problem, stability is instead obtained from the clamped boundary condition and the resulting coercivity of the modified form.

Strong ellipticity is preserved throughout the degeneration. Indeed,
\eqref{eq:H8a} gives
$$
K_\delta\geqslant-\frac13\mu>-\frac43\mu,
$$
and hence $c_{P,\delta}>0$ on $\overline\Omega$. The Legendre-Hadamard condition therefore holds uniformly in $\delta$, including at points of $\omega_0$, as shown in Proposition \ref{prop:elliptic}. What disappears in the limit is only the distinction between the two wave speeds. On $\omega_0$ one has
$c_{P,0}=c_S$, while the equation remains hyperbolic. Moreover, on every open subset contained in $\omega_0$, condition \eqref{eq:H8d} implies
$\nabla\mu=0$ and Remark \ref{rem:muconst} shows that the elastic operator reduces locally to a multiple of the vector Laplacian.

No smallness condition is imposed on the measure or geometry of the transition zone. Although the energy density is indefinite on $\omega_0$, where
$K_0=-\mu/3<0$, the modified form remains coercive on
$[H_0^1(\Omega)]^3$ uniformly with respect to $\delta$, as shown in Proposition \ref{prop:unifnorm}. This global coercivity follows from the integration-by-parts identity of Lemma \ref{lem:ibp} and the clamped boundary condition. Thus, stability in this setting is a property of the constrained body rather than of the energy density at each point.

	\medskip
	
	\underline{3. The physical origin of the family.}
	Three aspects of the family
	$\{\lambda_\delta\}_{\delta\in[0,1]}$ should be considered separately. They concern its shape, the direction of the degeneration and its endpoint. No single body of experimental work addresses all three.
	
	The shape of the family reflects the fact that $\mu$ remains fixed while only $\lambda_\delta$ changes. Here, fixed means independent of $\delta$ and not spatially constant. Thus, $\mu=\mu(x)$ may remain genuinely heterogeneous, as required in {\rm(H8)} and illustrated by Example \ref{ex:H8example}. This behavior is consistent with Gassmann relations of poroelasticity \cite{gassmann}, which relate the local moduli of a heterogeneous medium. At each point, the shear modulus of an isotropic porous rock is unaffected by the pore fluid because a shear deformation produces no change in pore volume, whereas the bulk modulus changes. We refer to \cite{smith} for the standard formulation of the fluid-substitution relations. Fluid substitution therefore preserves the profile $x\mapsto\mu(x)$ while changing $x\mapsto K(x)$. Departures from the predicted shear invariance have been observed in carbonates at seismic and ultrasonic frequencies \cite{adam}, so {\rm(H8)} should be regarded as an idealization rather than an exact description of a saturated rock. Poroelasticity, however, does not provide the required direction of degeneration, since saturation increases the bulk modulus.
	
	The direction is suggested by displacive phase transitions, particularly the $\alpha\!\to\!\beta$ transition of quartz discussed in Subsection \ref{sub:physical}. Resonance measurements show that the linear incompressibilities
	$K_1=(C_{11}+C_{12}+C_{13})/3$ and
	$K_3=(C_{33}+2C_{13})/3$ decrease sharply near the transition temperature of approximately
	$573^\circ\mathrm{C}$, whereas the shear moduli
	$C_{44}$,
	$(C_{11}+C_{33}-2C_{13})/4$ and
	$C_{66}=(C_{11}-C_{12})/2$
	decrease only slightly \cite{ohno}. The reduction in Young modulus exceeds $30\%$, while the shear moduli change by only a few percent \cite{peng}. Quartz is trigonal, so these elastic constants cannot be identified directly with the Lam\'e coefficients of an isotropic model. The measurements nevertheless support the relevant mechanism, namely the softening of the volumetric response with comparatively little change in the shear response. This disparity may also be reflected in the corresponding effective isotropic moduli.
	
	The endpoint is motivated by constrained media with negative bulk modulus, as discussed in part $2$. Thus, the fixed shear profile, the decrease of the bulk modulus and the limiting negative value are supported by distinct physical considerations. Assumption {\rm(H8)} combines these features into a single parameter-dependent model and Example \ref{ex:H8example} provides an explicit realization.
	
	\medskip
	\underline{4. The coupling condition and a scaling law.}
	Condition \eqref{eq:H8d} is the only hypothesis in {\rm(H8)} that couples the two moduli. It requires the shear modulus to become locally homogeneous wherever the bulk modulus approaches the critical value $-\mu/3$. In the composite interpretation, the inclusion is shear-uniform, while the surrounding region may remain heterogeneous. This behavior is expressed by the finiteness of $\mu_1^*$ and motivates the choice of cut-offs with disjoint transition regions in Example \ref{ex:H8example}.
	
	The condition also admits a dimensionless interpretation. Suppose that the shear modulus changes by a relative amount
	$\varepsilon_\mu=\Delta\mu/\underline\mu$ over a length $\ell$ and let $L$ denote the diameter of $\Omega$, so that
	$\lambda_1\simeq L^{-2}$. Define
	$$
	\kappa
	:=
	\inf_{\delta\in[0,1]}
	\inf_{\{x\,|\,\nabla\mu(x)\neq0\}}
	\bigl(\lambda_\delta+\mu\bigr)
	=
	\inf_{\delta\in[0,1]}
	\inf_{\{\nabla\mu\neq0\}}
	\rho\bigl(c_P^2-c_S^2\bigr).
	$$
	Assume that $\kappa>0$, as in Example \ref{ex:H8example}. Since
	$|\nabla\mu|\simeq\Delta\mu/\ell$, we have
	$\mu_1^*\simeq(\Delta\mu)^2/(\ell^2\kappa)$. Hence the requirement
	$\mu_1^*\lesssim\underline\mu\lambda_1$ contained in the smallness condition of {\rm(H8)} becomes
	\begin{equation}\label{eq:scaling}
		\varepsilon_\mu
		\lesssim
		\frac{\ell}{L}
		\sqrt{\frac{\kappa}{\underline\mu}}.
	\end{equation}
	Thus, the admissible shear contrast depends on the ratio between the heterogeneity length and the size of the body, together with the square root of the residual wave-speed gap relative to the shear modulus.
	
	This estimate also explains the role of the disjoint transition regions. If $\nabla\mu$ were nonzero where the wave-speed gap vanishes, then $\kappa$ would approach zero and \eqref{eq:scaling} would force
	$\varepsilon_\mu\to0$ as $\delta\to0$. This is the collapse to a constant shear modulus described in Remark \ref{rem:H8forced}. Separating the transition regions keeps the gap uniformly positive where $\mu$ varies and allows the shear profile to remain nonconstant in the limit. The term $6\mu_2$ imposes a similar condition on the second derivatives of $\mu$, formally of the order
	$\varepsilon_\mu\lesssim(\ell/L)^2$. Since this condition does not involve $\kappa$, it does not deteriorate as the wave-speed gap closes.
	
	\medskip
\underline{5. Interpretation of the robustness theorem.} Theorem \ref{thm:robust} can now be interpreted in terms of elastic bodies. Consider a limit medium containing a transition zone where the bulk modulus reaches the critical value $-\mu/3$, together with nearby media for which this value has not yet been attained. Assumption \eqref{eq:H8c} measures their discrepancy by $\delta$ in the uniform norm. The theorem provides exponential attractors for the corresponding dynamics satisfying
$$
\dist_{\rm sym}
\bigl(
\mathbf A_\delta^{\exp},
\mathbf A_0^{\exp}
\bigr)
\leqslant
C\delta^\tau,
$$
where $C>0$ and $\tau\in(0,1)$ are independent of $\delta$. These finite-dimensional sets attract the long-time response of the body to the load $b$ at a uniform exponential rate.

Thus, the long-time dynamics is stable as the bulk modulus approaches its critical value along the family prescribed in {\rm(H8)}. The limit medium may therefore be used as an idealization of the nearby media, with an error controlled quantitatively in the symmetric Hausdorff distance. This conclusion is stronger than the upper semicontinuity of global attractors obtained in \cite{bocanegra} for constant coefficients. It is two-sided and provides a H\"older rate, together with uniform bounds on the fractal dimension and the rate of exponential attraction.

Exponential attractors are positively invariant rather than fully invariant and are not unique. The family obtained in Theorem \ref{thm:robust} is therefore a robust selection of exponential attractors rather than a canonical family. Remark \ref{rem:usc-rate} explains the relation between this result and the upper semicontinuity of the global attractors.
\begin{Remark}[On the r\^ole of a constant shear modulus]\label{rem:H8forced}
		Condition \eqref{eq:H8d} is the uniform counterpart of the pointwise window
		\eqref{eq:window-sharp} and the choice matters. Had the heterogeneity of $\mu$ been measured
		globally instead, through \eqref{eq:window-global}, a family whose wave-speed gap closes
		\emph{uniformly}, $\epsilon_0^\delta=\min_{\overline\Omega}(\lambda_\delta+\mu)\to0$, would be
		forced to satisfy $|\nabla\mu|=O\bigl(\sqrt{\epsilon_0^\delta}\bigr)\to0$, so that the limit medium
		would necessarily have a constant shear modulus. The setting would then collapse to
		$\mu\equiv\mu_0$, $\HH=0$, $\mu_1^{*}=\mu_2=0$ and $\Le^\delta=-\Delta_e^\delta$ by
		Remark \ref{rem:muconst} and Section \ref{sec:framework} would be needed nowhere. The medium of
		Example \ref{rem:example} exhibits the collapse at the level of the constants: taking there
		$\epsilon_0=\mu_0\eta=\delta$ gives $\eta=\delta/\mu_0\to0$ and hence $\mu\to\mu_0$. That medium is
		of course not a family in the sense of {\rm(H8)}, where $\mu$ is fixed and only $\lambda_\delta$
		varies and no conflict with Example \ref{ex:H8example} should be read into it. It is precisely by
		localizing the requirement on the transition zone that \eqref{eq:H8d} allows the shear modulus to
		remain genuinely variable in the limit.
\end{Remark}

\begin{Example}[An explicit degenerating family with variable shear modulus]\label{ex:H8example}
	We construct a family satisfying {\rm(H8)} whose limit medium has a genuinely variable shear modulus. Let
	$\Omega=B_1(0)\subset\mathbb R^3$, so that $\lambda_1=\pi^2$ and choose radial cut-offs satisfying
	\begin{align*}
		&\varphi_1\in C^\infty([0,1]),
		\,\,\,
		\varphi_1=0\ \mbox{on}\ [0,\tfrac12],
		\,\,\,
		0\leqslant\varphi_1\leqslant1,
		\,\,\,
		\varphi_1'>0\ \mbox{on}\ (\tfrac12,1),
		\\
		&\varphi_2\in C^\infty(\overline\Omega),
		\,\,\,
		\varphi_2=0\ \mbox{on}\ \{|x|\leqslant\tfrac14\},
		\,\,\,
		\varphi_2=1\ \mbox{on}\ \{|x|\geqslant\tfrac12\},
	\end{align*}
	with $0\leqslant\varphi_2\leqslant1$ and
	$\varphi_2>0$ on $\{|x|>\tfrac14\}$. Fix
	$\mu_0>0$, $c\geqslant1$ and $\eta\in(0,1)$ and define
	\begin{equation}\label{eq:H8example}
		\mu(x)
		=
		\mu_0\bigl[1+\eta\varphi_1(|x|)\bigr],
		\,\,\,
		\lambda_\delta(x)
		=
		-\mu(x)+\delta+(c-\delta)\varphi_2(x),
		\,\,\,
		\delta\in[0,1].
	\end{equation}
	Since $\varphi_1$ vanishes near the origin,
	$\mu\in C^\infty(\overline\Omega)$ and
	$\underline\mu=\mu_0>0$.
	
	We have
	$$
	\lambda_\delta+\mu
	=
	\delta+(c-\delta)\varphi_2
	\in[\delta,c],
	$$
	so \eqref{eq:H8a} holds. Moreover,
	$\lambda_\delta-\lambda_0=\delta(1-\varphi_2)$ and hence
	$$
	\|\lambda_\delta-\lambda_0\|_{L^\infty(\Omega)}
	=
	\delta.
	$$
	Thus, \eqref{eq:H8c} holds with equality. The functions
	$\lambda_\delta$ are smooth and their $C^{0,1}(\overline\Omega)$ norms are bounded uniformly in $\delta$, which verifies \eqref{eq:H8b}.
	
	To verify \eqref{eq:H8d}, observe that
	$\nabla\mu=0$ on $\{|x|\leqslant\tfrac12\}$, whereas
	$\lambda_\delta+\mu=c$ on
	$\{|x|\geqslant\tfrac12\}$. Therefore,
	$$
	\mu_1^*
	=
	\sup_{\delta\in[0,1]}
	\sup_{x\in\overline\Omega}
	\frac{|\nabla\mu(x)|^2}
	{\lambda_\delta(x)+\mu(x)}
	\leqslant
	\frac{\mu_0^2\eta^2\|\varphi_1'\|_\infty^2}{c}
	<
	\infty.
	$$
	At points where $\delta=0$ and $\lambda_0+\mu=0$, the quotient is understood according to the convention in \eqref{eq:H8d}. Since
	$\mu_2=O(\eta)$ and $\mu_1^*=O(\eta^2)$, the uniform smallness condition in {\rm(H8)} holds for all sufficiently small $\eta$ whenever
	$2M<\mu_0\pi^2$.
	
	For the limit medium,
	$$
	\lambda_0+\mu
	=
	c\varphi_2,
	$$
	and therefore
	$$
	\omega_0
	=
	\{|x|\leqslant\tfrac14\}.
	$$
	This set has nonempty interior, while $\mu$ is not constant and is strictly increasing in $|x|$ on
	$\{\tfrac12<|x|<1\}$. The separation of the transition regions of the two cut-offs keeps $\mu_1^*$ finite while allowing $\mu$ to remain variable.
	
	Finally, on $\omega_0$ one has
	$$
	3\lambda_\delta+2\mu
	=
	3\delta-\mu_0
	\,\,\,\mbox{and}\,\,\,
	\frac{K_\delta}{\mu}
	=
	-\frac13+\frac{\delta}{\mu_0}.
	$$
	Hence $3\lambda_\delta+2\mu<0$ on $\omega_0$ whenever
	$\delta<\mu_0/3$. At $\delta=0$, the bulk modulus reaches the endpoint
	$K_0/\mu=-1/3$ and
	$\lambda_0=-\underline\mu$ on $\omega_0$. Thus, the limit medium is covered by neither Proposition \ref{prop:unperturbed} nor the strict sufficient condition in Remark \ref{rem:korn}, while the modified framework remains applicable.
\end{Example}
	\medskip
	\begin{Remark}[Fixed phase spaces for the robustness analysis]
		The comparison of attractors corresponding to different values of $\delta$ requires a common phase space equipped with a fixed metric. For this reason, we use the $\delta$-independent spaces and norms in \eqref{eq:fixedspaces}. The natural energy norms involving
		$(\Le^\delta)^{1/2}$ and $\Le^\delta$ depend on $\delta$ and therefore do not provide a fixed metric in which the symmetric Hausdorff distance can be evaluated directly. Propositions \ref{prop:unifnorm} and \ref{prop:unifH2} show that these parameter-dependent norms are uniformly equivalent to the fixed norms in \eqref{eq:fixedspaces}. Consequently, all members of the family can be compared in the same phase spaces with constants independent of $\delta$.
	\end{Remark}
	
\subsection{The parameter-dependent dynamics}
We now place the family of problems associated with $\{\lambda_\delta\}_{\delta\in[0,1]}$ in a common functional setting. The condition \eqref{eq:H8d} is the uniform counterpart of the sharp
condition \eqref{eq:window-sharp} discussed in Remark \ref{rem:window}. We therefore use the same shift $\mu_1=\mu_1^{*}$ for every $\delta\in[0,1]$ and define
\begin{equation}\label{eq:formdelta}
	\Le^\delta=-\Delta_e^\delta-\HH+\mu_1^{*} I
\end{equation}
through the symmetric form
\begin{align*}
	\langle u,v\rangle_{e,\delta}
	&= \! \int_\Omega
	\Bigl[
	(\lambda_\delta\!+\! \mu)\dive(u)\dive(v)
	+\mu \nabla u \!\cdot\! \nabla v +
	\dive(v) u \!\cdot\! \nabla\mu
	+\dive(u) v \!\cdot\! \nabla\mu
	+\mu_1^{*} u \!\cdot\! v
	\Bigr] dx .
\end{align*}
Since $\mu$, $\HH$ and $\mu_1^{*}$ are fixed, the family of operator $\{\Le^\delta\}_{\delta\in[0,1]}$ depends on $\delta$ only through $\lambda_\delta$.
To compare the resulting dynamics, we use the fixed phase spaces
\begin{equation}\label{eq:fixedspaces}
\mathscr H^0= [H_0^1(\Omega)]^3\times[L^2(\Omega)]^3\,\,\,\mbox{and}\,\,\,
\mathscr H^1= [H^2(\Omega)\cap H_0^1(\Omega)]^3	\times[H_0^1(\Omega)]^3,
\end{equation}
endowed with the norms
$$
\|(u,v)\|_{\mathscr H^0}^2
=
\|\nabla u\|_2^2+\|v\|_2^2,
\,\,\,\mbox{and}\,\,\,
\|(u,v)\|_{\mathscr H^1}^2
=
\|u\|_{H^2}^2+\|\nabla v\|_2^2.
$$
The next two propositions show that these fixed norms are uniformly
equivalent to the natural norms induced by $\Le^\delta$.
\begin{Proposition}[Uniform equivalence of the energy norms] \label{prop:unifnorm} Under {\rm(H8)}, there exists ${C_{\rm eq}}>0$, independent of $\delta$, such that
\begin{equation}\label{eq:unifcoerc}
\underline\mu\|\nabla u\|_2^2 \leqslant \|u\|_{e,\delta}^2	\leqslant {C_{\rm eq}}\|\nabla u\|_2^2,\,\,\,\forall\, u\in[H_0^1(\Omega)]^3,\, \delta\in[0,1].
\end{equation}
Consequently, $D\bigl((\Le^\delta)^{1/2}\bigr)=[H_0^1(\Omega)]^3$ and the norms $\|(\Le^\delta)^{1/2}u\|_2$ and $\|\nabla u\|_2$ are uniformly equivalent with respect to $\delta\in[0,1]$.
\end{Proposition}
\begin{proof}
The lower bound follows from the proof of Proposition \ref{prop:coercive}, with $\mu_1$ replaced by $\mu_1^{*}$. Indeed, \eqref{eq:H8a} ensures that $\lambda_\delta+\mu\geqslant0$, while \eqref{eq:H8d} provides the	pointwise estimate needed in \eqref{eq:young-form}. For the upper bound, \eqref{eq:H8b} gives $\max_{\overline\Omega}(\lambda_\delta+\mu)\leqslant\Lambda+\overline\mu$ uniformly in $\delta$. The remaining terms are controlled by $\|\nabla\mu\|_\infty$, $\mu_1^{*}$  and the Poincar\'e inequality. Therefore, ${C_{\rm eq}}$ is independent of $\delta$. The characterization of the square-root domain follows from	Proposition \ref{prop:selfadj}.
\end{proof}
\begin{Proposition}[Uniform elliptic regularity] \label{prop:unifH2} Under {\rm(H8)}, one has
$$D(\Le^\delta)=[H^2(\Omega)\cap H_0^1(\Omega)]^3,\,\,\,\forall\, \delta\in[0,1].$$
Moreover, there exists	$C_{\rm ell}>0$, independent of $\delta$, such that
\begin{equation}\label{eq:uniform-H2} C_{\rm ell}^{-1}\|u\|_{H^2}
\leqslant \|\Le^\delta u\|_2 \leqslant C_{\rm ell}\|u\|_{H^2},\,\,\,u\in D(\Le^\delta).
\end{equation}
\end{Proposition}
\begin{proof}
The domain characterization follows from Theorem \ref{thm:domain}. It remains only to verify that the constants can be chosen independently of $\delta$. The principal part of $\Le^\delta$ satisfies the Legendre condition with ellipticity constant $\underline\mu$, independently of $\delta$,	by Proposition \ref{prop:elliptic} and \eqref{eq:H8a}. Moreover, \eqref{eq:H8b} yields a uniform bound for the Lipschitz norms of the principal coefficients. The $[H^2(\Omega)]^3$-regularity estimate for strongly elliptic Dirichlet systems therefore gives $\|u\|_{H^2} \leqslant C\bigl[\|\Le^\delta u\|_2+\|u\|_2\bigr]$ with $C$ independent of $\delta$. Proposition \ref{prop:unifnorm}, together with the Poincar\'e inequality, allows the $L^2$-term to be absorbed, yielding the first inequality in \eqref{eq:uniform-H2}. The reverse inequality follows directly from definition of $\Le^\delta$, the uniform coefficient bound \eqref{eq:H8b} and the boundedness of the lower-order terms.
\end{proof}
The preceding results allow us to associate a dynamical system with each value of the parameter. We record this fact and introduce the notation used in the sequel.
\begin{Corollary}[Dynamics for each parameter] \label{cor:dynamics-delta} Assume {\rm(H8)}. For every $\delta\in[0,1]$, the problem associated with $\lambda_\delta$ generates a continuous nonlinear semigroup $\{S_\delta(t)\}_{t\geqslant0}$ on $\mathscr H^0$. The corresponding homogeneous problem generates an exponentially stable linear semigroup, denoted by $\{S_h^\delta(t)\}_{t\geqslant0}$. Moreover, $\{S_\delta(t)\}_{t\geqslant0}$ possesses a unique compact global attractor, denoted by $\mathbf A_\delta$ and at least one exponential attractor in $\mathscr H^0$.
\end{Corollary}
\begin{proof}
Propositions \ref{prop:unifnorm} and \ref{prop:unifH2} show that the energy and strong phase spaces associated with $\Le^\delta$ coincide with the fixed spaces \eqref{eq:fixedspaces}. Moreover, the coercivity and ellipticity arguments of Sections~\ref{sec:framework}-\ref{sec:exp} remain valid when the strict condition $\lambda+\mu>0$ in {\rm(H2)} is replaced by \eqref{eq:H8a}-\eqref{eq:H8d}. Consequently, the well-posedness, linear stability and attractor results established in those sections apply to each fixed $\delta\in[0,1]$.
\end{proof}
The preceding corollary concerns each fixed value of $\delta$. The
energy bounds and the decay of the homogeneous dynamics can, in fact,
be made uniform over the entire parameter interval.
\begin{Lemma}[Uniform energy and linear decay] \label{lem:unifdyn1} Assume {\rm(H8)}. The constants in Proposition $\ref{prop:energy}$ and Corollary $\ref{cor:bounded}$ can be chosen independently of $\delta\in[0,1]$. In particular, the energy equivalence constants $C_0,C_1,C_2$ and the boundedness function $R$ are uniform with respect to $\delta$.
Moreover, there exist constants $C,\eta>0$, independent of $\delta$, such that
$$	\|S_h^\delta(t)\U_0\|_{\mathscr H^0} \leqslant 	Ce^{-\eta t}\|\U_0\|_{\mathscr H^0},\,\,\,	t\geqslant0,\,\,\, \U_0\in\mathscr H^0,\,\,\,\delta\in[0,1].$$
\end{Lemma}
\begin{proof} The proof of Proposition \ref{prop:energy} shows that the energy constants depend only on the structural quantities $M$, $m_f$, $|\Omega|$, $\|b\|_2$, $\lambda_1$, on $\mu_2$ and $\mu_1^{*}$. All these quantities are independent of  $\delta$. Indeed, $\mu$ is fixed, the bound $\mu_1^{*}$ is uniform by \eqref{eq:H8d} and Proposition \ref{prop:unifnorm} gives \eqref{eq:unifcoerc} with $\underline\mu$ and ${C_{\rm eq}}$ independent of $\delta$. Hence the constants $C_0,C_1,C_2$ in Proposition \ref{prop:energy} can be chosen	uniformly. The proof of Corollary \ref{cor:bounded} uses only these energy estimates and therefore its boundedness function $R$ is also independent of $\delta$.

For the homogeneous dynamics, the argument of Remark \ref{rem:homog} depends on the damping coefficient $\alpha$, the Poincar\'e constant and the coercivity and continuity bounds for the elastic form. These quantities are uniform by Proposition~\ref{prop:unifnorm}. Consequently, the constants in the decay estimate can be chosen independently of $\delta$.
\end{proof}
\subsection{Main result}
We first recall the discrete counterpart of the notion of an
exponential attractor.
\begin{Definition}\label{discrete-exp}
	Let $X$ be a metric space and let $L:X\to X$ be a map. The iterates of $L$ define the discrete semigroup $S^n x=L^n x$, $x\in X$ and $n\in\mathbb N$. A set $\mathbf M\subset X$ is called an \underline{exponential attractor} for the map $L$ if it is compact in $X$, has finite fractal dimension, is positively invariant under $L$ and attracts bounded subsets of $X$
	exponentially with respect to the number of iterations. More precisely,
	$L(\mathbf M)\subseteq\mathbf M$ and there exists $\gamma>0$ such
	that, for every bounded set $B\subset X$, there is a constant $C(B)>0$
	satisfying
	$$
	\dist_X\bigl(L^nB,\mathbf M\bigr)
	\leqslant
	C(B)e^{-\gamma n},
	\,\,\, n\in\mathbb N.
	$$
\end{Definition}
We now consider families of discrete exponential attractors depending
on a perturbation parameter. Robustness means both upper and lower
semicontinuity with respect to this parameter, equivalently, continuity
in the symmetric Hausdorff distance. The existence of such a family,
together with uniform attraction and finite-dimensionality estimates,
follows from the fixed-space version of
\cite[Theorem~3.1]{FGMZ} stated below.
\begin{Theorem}\label{thm:robust-abstract}
Let $H_1$ and $H_0$ be Banach spaces, with $H_1$ compactly embedded into $H_0$. Suppose that $\mathbf B\subset H_0$ is closed and bounded and that $\{L_\delta:\mathbf B\to\mathbf B\}_{\delta\in[0,1]}$ is a family of maps satisfying the following conditions with constants independent of $\delta$.
\begin{itemize}
\item[$1.$] $L_\delta=L_\delta^{(0)}+L_\delta^{(1)}$, where $L_\delta^{(0)},L_\delta^{(1)}$ map $\mathbf{B}$ into $H_0$ and, for all $\U^1,\U^2\in\mathbf{B}$,$$
\bigl\|L_\delta^{(0)}\U^1-L_\delta^{(0)}\U^2\bigr\|_{H_0}\!\leqslant \varkappa\|\U^1-\U^2\|_{H_0} \,\,\,\mbox{and}\, \,\,\,
\bigl\|L_\delta^{(1)}\U^1-L_\delta^{(1)}\U^2\bigr\|_{H_1}\! \leqslant c_*\|\U^1-\U^2\|_{H_0},$$
with $\varkappa<\tfrac12$.
\item[$2.$] $\|L_\delta^{n}\U-L_0^{n}\U\|_{H_0}\leqslant C \delta L^{n}$, for every $n\in\mathbb N$ and every $\U\in\mathbf{B}$.
\end{itemize}
Then there exists a family $\{\mathbf M_\delta\}_{\delta\in[0,1]}$ such that $\mathbf M_\delta\subset\mathbf B$ is an exponential attractor for $L_\delta$. Moreover, there exist constants $C_3,C_4,C_5,\nu>0$ and $\tau\in(0,1)$, independent of $\delta$, such that
	\begin{align}
		\dist_{H_0}
		\bigl(L_\delta^n\mathbf B,\mathbf M_\delta\bigr)
		&\leqslant C_3e^{-\nu n},\,\,\,\, n\in\mathbb N,\, \delta\in[0,1],
		\label{robusta}\\
		\dim_F^{H_0}(\mathbf M_\delta)
		&\leqslant C_4,\,\,\,\,\delta\in[0,1],
		\label{robustb}\\
		\dist_{\rm sym}
		\bigl(\mathbf M_\delta,\mathbf M_0\bigr)
		&\leqslant C_5\delta^\tau,\,\,\,\, \delta\in[0,1].
		\label{robustc}
	\end{align}
The constants $C_3,C_4,C_5$ and $\tau$ can be explicitly determined from the constants appearing in the assumptions. Here, $\dist_{\rm sym}(X,Y)=\max\{\dist_{H_0}(X,Y),\dist_{H_0}(Y,X)\}$, $X,Y\subset H_0$.
\end{Theorem}
\begin{proof}
We apply \cite[Theorem 3.1]{FGMZ} with $\mathscr{E}(\delta)=H_0$, $\mathscr{E}_1(\delta)=H_1$, $\mathbf B_\delta=\mathbf B$ and $\Pi_\delta=I_{\mathbf B}$	for every $\delta\in[0,1]$. Since the embedding $H_1\subset H_0$ is compact, the unit ball of $H_1$ is totally bounded in $H_0$. Hence its Kolmogorov entropy $\mathbb M(r)=\log N_{r}\bigl[B_{H_1}(0,1),H_0\bigr]$ is finite for every $r>0$, where $N_{r}\bigl[B_{H_1}(0,1),H_0\bigr]$ is the minimal number of balls of radius $r$ of $H_0$ needed to cover $B_{H_1}(0,1)$. As the spaces and their norms are independent of $\delta$, the same function $\mathbb M$ applies to the entire family. This verifies \cite[(3.1)]{FGMZ}. The compatibility condition \cite[(3.2)]{FGMZ} holds with the two constants appearing there equal to $1$ and $0$, since $\mathscr E(\delta)=\mathscr E(0)=H_0$ and $\mathbf B_\delta=\mathbf B_0=\mathbf B$. Assumption $1$ is precisely the uniform contractive-smoothing decomposition required in 	\cite[(3.3)]{FGMZ}.	Finally, $\Pi_\delta=I_{\mathbf B}$ maps $\mathbf B_\delta$ onto	$\mathbf B_0$  and \cite[(3.4)]{FGMZ} reduces to $\|L_\delta^n\U-L_0^n\U\|_{H_0} \leqslant C\delta L^n$, which is assumption $2$. Thus, all the hypotheses of \cite[Theorem 3.1]{FGMZ} are satisfied and its conclusions yield	\eqref{robusta}-\eqref{robustc}.
\end{proof}

The decomposition allowed by assumption $1$ is essential here. The criterion of
\cite[Theorem 4.4]{EMZ2} requires instead the full map $L_\delta$ to be smoothing from $H_0$ into
$H_1$, a property which is not available for \eqref{eq:main-abs}: for $L_\delta=S_\delta(t_*)$ the
homogeneous component $S_h^\delta(t_*)$ does not regularize and only the inhomogeneous component
$S_c^\delta(t_*)$ produces the required gain of regularity. This is also what distinguishes the
present hyperbolic problem from the strongly damped wave equation of Remark \ref{rem:PS}, whose
linear dynamics is smoothing.

The main result of this section is given by the following theorem.
\begin{Theorem}[Robust family of exponential attractors] \label{thm:robust}
	Assume {\rm(H8)}. Then there exists a family $\{\mathbf A_\delta^{\exp}\}_{\delta\in[0,1]}$ such that $\mathbf A_\delta^{\exp}$ is an exponential attractor for $(\mathscr H^0,S_\delta(t))$ for every $\delta\in[0,1]$. Moreover:
	\begin{enumerate}
		\item[$1.$] The fractal dimensions are uniformly bounded, that is, $\sup_{\delta\in[0,1]} \dim_F^{\mathscr H^0} \bigl(\mathbf A_\delta^{\exp}\bigr)<\infty$.

		\item[$2.$] Each exponential attractor contains the corresponding global attractor, that is, $	\mathbf A_\delta\subseteq\mathbf A_\delta^{\exp}$.

		\item[$3.$] The family is H\"older continuous at $\delta=0$ in the symmetric
		Hausdorff distance. More precisely, there exist $C>0$ and $\tau\in(0,1)$ such that $
		\dist_{\rm sym}
		\bigl(\mathbf A_\delta^{\exp},\mathbf A_0^{\exp}\bigr)
		\leqslant C\delta^\tau$, $\delta\in[0,1]$.
	\end{enumerate}
\end{Theorem}
\subsection{Proof of Theorem \ref{thm:robust}}
\subsubsection{Verification of the hypotheses of Theorem \ref{thm:robust-abstract}}
The first step in verifying the hypotheses of Theorem \ref{thm:robust-abstract} is to construct the closed bounded set $\mathbf B\subset\mathscr H^0$, independent of $\delta$.
\begin{Lemma}[Construction of a common invariant set]\label{lem:unifdyn3}
Assume {\rm(H8)}. There exist constants $R_\star,{t_{\mathbf B}}>0$, independent of $\delta$, such that the closed ball $$\mathbf B = \bigl\{\U\in\mathscr H^1 \,|\, \|\U\|_{\mathscr H^1}\leqslant R_\star\bigr\}$$ satisfies $S_\delta(t)\mathbf B\subseteq\mathbf B$ for every $t\geqslant t_{\mathbf B}$ and $\mathbf A_\delta\subseteq\mathbf B$  for every $\delta\in[0,1]$. Moreover, $\mathbf B$ is compact in $\mathscr H^0$.
\end{Lemma}
\begin{proof}
We apply the bootstrap procedure of Lemma \ref{lem:bootstrap} up to the level $r=1$. This level is not reached directly from $r=0$. Indeed, by \eqref{eq:step}, the first admissible increment is $s(0)=\theta=\frac{3-p}{2}\leqslant\frac12$. Accordingly, as in the proof of Theorem \ref{thm:B1}, we consider a finite chain $0=k_0<k_1<\cdots<k_N=1$ with $N\leqslant\left\lceil\frac{2}{3-p}\right\rceil$ in which every increment is admissible for the bootstrap estimate.	At each level $k_j$, the proof of Lemma \ref{lem:bootstrap} yields a perturbed functional $\Upsilon_{k_j}^\delta$ and constants 	$\beta_j,C_j,c_j,c_j'>0$ such that
\begin{equation}\label{eq:unifgronwall}
\frac{d}{dt}\Upsilon_{k_j}^\delta(t)+\beta_j\Upsilon_{k_j}^\delta(t)\leqslant C_j\,\,\,\,
\mbox{and}
\,\,\,\, c_j\|\U\|_{\mathscr H^{k_j}}^2-C_j
\leqslant \Upsilon_{k_j}^\delta(\U)	\leqslant c_j'\|\U\|_{\mathscr H^{k_j}}^2+C_j.
\end{equation}
Once the bound at level $k_{j-1}$ has been fixed, these constants depend only on $\alpha$, $\underline\mu$, $\lambda_1$, $\|b\|_2$, $\mu_2$, $\mu_1^{*}$, the growth constants in {\rm(H5)} and {\rm(H7)}, the uniform 	elliptic constant $C_{\rm ell}$ and the previously obtained bound. The Propositions \ref{prop:unifnorm} and \ref{prop:unifH2}, together with interpolation, give the corresponding uniform equivalences at all intermediate levels $k_j$. Hence, by induction on $j$, the constants in \eqref{eq:unifgronwall} can be chosen independently of $\delta$.

Integrating \eqref{eq:unifgronwall}, we obtain
\begin{equation}\label{derivada}
\Upsilon_{k_j}^\delta(t)	\leqslant	e^{-\beta_jt}\Upsilon_{k_j}^\delta(\U_0)	+\frac{C_j}{\beta_j},\,\,\, \forall \, t \geqslant 0.
\end{equation}
If $\|\U_0\|_{\mathscr H^{k_j}}\leqslant R$, then $\Upsilon_{k_j}^\delta(\U_0)\leqslant c_j'R^2+C_j$.	Therefore, for
$$ t_j(R) =
\frac1{\beta_j}
\left[\log\frac{\beta_j(c_j'R^2+C_j)}{C_j}\right]_+,\,\,\, [s]_+:=\max\{s,0\},
$$
we have $\Upsilon_{k_j}^\delta(t)\leqslant	\frac{2C_j}{\beta_j},\,\,\,t\geqslant t_j(R)$. It follows from \eqref{eq:unifgronwall} that $\|S_\delta(t)\U_0\|_{\mathscr H^{k_j}}\leqslant R_j$ for  every $t\geqslant t_j(R)$, where 	$R_j^2	= 	\frac1{c_j}	\left[\frac{2C_j}{\beta_j}+C_j\right]$. Both $R_j$ and $t_j(R)$ are independent of $\delta$.

We now iterate this argument along the chain $k_0,\ldots,k_N$. After reaching the uniform bound at level $k_j$, we apply the estimate at level $k_{j+1}$ to the time-translated trajectory $t\longmapsto	S_\delta(t+T_j)\U_0$, where $T_j:=t_j(R_{j-1})$ is the entrance time produced at level $k_j$. Since the chain has finitely many levels, induction gives a radius $R_\star>0$ and an entrance time ${t_{\mathbf B}}>0$, both independent of	$\delta$, such that  $S_\delta(t)\mathbf B\subseteq\mathbf B$ for every $t\geqslant {t_{\mathbf B}}$ and $\delta\in[0,1]$.

It remains to prove that $\mathbf A_\delta\subseteq\mathbf B$. For each	fixed $\delta$, Theorem \ref{thm:attractor} shows that $\mathbf A_\delta$ is bounded in $\mathscr H^1$, write $	\sup_{\U\in\mathbf A_\delta} \|\U\|_{\mathscr H^1}\leqslant\varrho_\delta$.	Applying the preceding bootstrap argument with initial radius $\varrho_\delta$, we find a time $T_\delta>0$ such that $S_\delta(T_\delta)\mathbf A_\delta\subseteq\mathbf B$. Although $T_\delta$ may depend on $\delta$, the radius of $\mathbf B$ does not. Since the global attractor is strictly invariant, we obtain $\mathbf A_\delta\subseteq\mathbf B$.
Finally, since $\mathscr H^1$ is reflexive and compactly embedded into
$\mathscr H^0$, the closed bounded set $\mathbf B$ is compact in
$\mathscr H^0$.
\end{proof}
\begin{Remark}[Uniform boundedness of the orbits issuing from $\mathbf B$] \label{orbit-B} Although $\mathbf B$ is only known to satisfy $S_\delta(t)\mathbf B\subseteq\mathbf B$ for every $t\geqslant t_{\mathbf B}$, the forward trajectories issuing from $\mathbf B$ are uniformly bounded in $\mathscr H^1$ for all $t\geqslant0$. Indeed, at the final bootstrap level $k_N=1$, if $\U_0\in\mathbf B$, then from \eqref{eq:unifgronwall}, we obtain $\Upsilon_1^\delta(\U_0)\leqslant c_N'R_\star^2+C_N$. Therefore, using \eqref{derivada}, we find
$$\Upsilon_1^\delta(t)\leqslant e^{-\beta_Nt}\bigl(c_N'R_\star^2+C_N\bigr)
+\frac{C_N}{\beta_N}, \,\,\, \forall\,t\geqslant0.
$$
Combining this estimate with $c_N\|S_\delta(t)\U_0\|_{\mathscr H^1}^2-C_N \leqslant \Upsilon_1^\delta(t)$, we obtain
$$\|S_\delta(t)\U_0\|_{\mathscr H^1}^2 \leqslant \frac1{c_N}
\left[ c_N'R_\star^2+2C_N+\frac{C_N}{\beta_N}\right]=: C_\star^{2}.$$
Hence, $\sup_{\delta\in[0,1]} \sup_{\U_0\in\mathbf B}\sup_{t\geqslant0}\|S_\delta(t)\U_0\|_{\mathscr H^1}	\leqslant C_\star$.
Notice that this does not imply that $S_\delta(t)\mathbf B\subseteq\mathbf B$ for $0\leqslant t<{t_{\mathbf B}}$, since $C_\star$ may be larger than $R_\star$. It does, however, provide the uniform bound on forward orbits required in Remark \ref{rem-set} and hence allows Lemmas \ref{lem:Df} and \ref{lem:criterion} to be applied to $\mathbf B$ with constants independent of $\delta$.
\end{Remark}
We now define the discrete maps to which Theorem \ref{thm:robust-abstract} will be applied.
\begin{Lemma}\label{lem:unifdyn4}
Let $\mathbf B$ be the common ball given by Lemma $\ref{lem:unifdyn3}$. The constants in Lemmas $\ref{lem:Df}$ and $\ref{lem:criterion}$, when restricted to trajectories issuing from $\mathbf B$, can be chosen independently of $\delta\in[0,1]$. Moreover, there exists $t_*\geqslant {t_{\mathbf B}}$, independent of $\delta$, such that the	fixed-time maps
$$ L_\delta=S_\delta(t_*)\big|_{\mathbf B}:\mathbf B\to\mathbf B$$
satisfy the contractive-smoothing estimates of Theorem $\ref{thm:robust-abstract}$ with constants independent of $\delta$.
\end{Lemma}
\begin{proof} By Remark \ref{orbit-B}, the forward trajectories issuing from $\mathbf B$ are uniformly bounded. Hence, in view of Remark \ref{rem-set}, Lemmas \ref{lem:Df} and \ref{lem:criterion} apply to $\mathbf B$. Their constants depend only on $C_\star$, the growth constants in {\rm(H5)} and {\rm(H7)}, the fixed Sobolev embeddings and the uniform	energy, continuous-dependence and elliptic regularity estimates. Consequently, these constants can be chosen independently of	$\delta$.

Let $C,\eta>0$ be the uniform constants in the linear decay estimate of Lemma \ref{lem:unifdyn1}  and set	$t_*=\max\left\{{t_{\mathbf B}},\frac1\eta[\log(8C)]_+\right\}$. Then $t_*$ is independent of $\delta$ and $Ce^{-\eta t_*}\leqslant\frac18$. Moreover, since $t_*\geqslant {t_{\mathbf B}}$, Lemma \ref{lem:unifdyn3} gives $S_\delta(t_*)\mathbf B\subseteq\mathbf B$. Thus, $L_\delta:\mathbf B\to\mathbf B$ is well defined.	As in Lemma \ref{lem:criterion}, decompose $L_\delta	=L_\delta^{(0)}+L_\delta^{(1)}$, where $L_\delta^{(0)}\U=S_h^\delta(t_*)\U$ and $L_\delta^{(1)}\U	= S_c^\delta(t_*)\U$.	The uniform decay estimate gives
$$
\|L_\delta^{(0)}\U^1-L_\delta^{(0)}\U^2\|_{\mathscr H^0}
\leqslant
\frac18\|\U^1-\U^2\|_{\mathscr H^0}
$$
and
$$
\|L_\delta^{(1)}\U^1-L_\delta^{(1)}\U^2\|_{\mathscr H^1}
\leqslant
c_*\|\U^1-\U^2\|_{\mathscr H^0},
$$
where $c_*>0$ is independent of $\delta$. Therefore, the maps $L_\delta$ satisfy the contractive-smoothing hypothesis of	Theorem \ref{thm:robust-abstract} uniformly with respect to $\delta$.
\end{proof}
We now verify the parameter-continuity hypothesis of Theorem \ref{thm:robust-abstract}. The key point is to compare $S_\delta(t)$ and $S_0(t)$ for initial data in $\mathbf B\subset\mathscr H^1$. The strong regularity is essential because the perturbation of the elastic operator produces a term in $[H^{-1}(\Omega)]^3$, which cannot be paired with a velocity belonging
only to $[L^2(\Omega)]^3$.
\begin{Lemma}\label{thm:close}
Assume {\rm(H8)} and let $T>0$. There exists a constant $C(T)>0$, depending on $T$, $R_\star$ and the structural constants, but	independent of $\delta$, such that
\begin{equation*}
\sup_{t\in[0,T]}\bigl\|S_\delta(t)\U_0-S_0(t)\U_0\bigr\|_{\mathscr H^0}
\leqslant {C(T)} \delta,\,\,\, \forall\,\U_0\in\mathbf{B},\ \forall\,\delta\in[0,1].
\end{equation*}
\end{Lemma}
\begin{proof}
Fix $\U_0=(u_0,u_1)\in\mathbf{B}$ and let $u^\delta$ and $u^0$ be the corresponding solutions of \eqref{eq:main-abs} with Lam\'e coefficients $\lambda_\delta$ and $\lambda_0$, respectively. By Remark \ref{orbit-B}, there exists $R>0$, independent of
 $\delta$, $t$ and $\U_0\in\mathbf B$, such that
\begin{equation}\label{eq:unifH1bound}
\|u^\delta(t)\|_{H^2}+\|u^\delta_t(t)\|_{H^1}+\|u^0(t)\|_{H^2}+\|u^0_t(t)\|_{H^1}\leqslant R,\,\,\,\forall\, t\geqslant0.
\end{equation}

\underline{Step 1: the equation for $w$.}
Since $\mu$, $\rho$, $\HH$ and $\mu_1^{*}$ are the same for both
problems, the dependence of $\Le^\delta$ on $\delta$ occurs only
through $\lambda_\delta$. Writing $\zeta_\delta=\lambda_\delta-\lambda_0$, so that
$\|\zeta_\delta\|_{L^\infty(\Omega)}\leqslant\delta$ by \eqref{eq:H8c}, the contributions of $\mu$,
$\HH$ and $\mu_1^{*}$ cancel in the difference and \eqref{eq:lame-expanded} gives
$(\Le^\delta-\Le^0)u=-\nabla\bigl[\zeta_\delta\,\dive u\bigr]$. Subtracting the equations satisfied by $u^\delta$ and $u^0$ and writing $\Le^\delta u^\delta-\Le^0u^0=\Le^\delta w+(\Le^\delta-\Le^0)u^0$, we obtain
\begin{equation}\label{eq:weq}
w_{tt}+\Le^\delta w+\alpha w_t+\bigl[f(u^\delta)-f(u^0)\bigr]
=	\mu_1^{*} w-\HH w	+\nabla\bigl[\zeta_\delta\,\dive u^0\bigr]\,\,\,\text{in }\Omega\times\mathbb R^+,
\end{equation}
in the sense of $[H^{-1}(\Omega)]^3$,
together with $w=0$ on $\partial\Omega\times\mathbb R^+$, $w(0)=0$ and $w_t(0)=0$.

\underline{Step 2: the energy identity.}
Multiplying \eqref{eq:weq} by $w_t$ and integrating over $\Omega$, we
note that $\U_0\in\mathbf B\subset\mathscr H^1$ and Theorem \ref{thm:global}$(2)$ yield
$w_t(t)\in[H_0^1(\Omega)]^3$. Hence, the perturbation term is well
defined in the duality between $[H^{-1}(\Omega)]^3$ and
$[H_0^1(\Omega)]^3$  and integration by parts gives
$$\langle\nabla[\zeta_\delta\dive u^0],w_t\rangle
=-\int_\Omega\zeta_\delta\dive u^0\dive w_t\,dx.$$
Moreover,
$\langle w_{tt},w_t\rangle
=\frac12\frac{d}{dt}\|w_t\|_2^2$ and
$\langle\Le^\delta w,w_t\rangle
=\langle w,w_t\rangle_{e,\delta}
=\frac12\frac{d}{dt}\|w\|_{e,\delta}^2$. Therefore,
\begin{equation}\label{eq:wenergy}
\frac12 \frac{d}{dt}
\Bigl[\|w_t\|_2^2 \!+\! \|w\|_{e,\delta}^2\Bigr]\!+\! \alpha\|w_t\|_2^2
=\! -\bigl\langle f(u^\delta)\!-\! f(u^0),w_t\bigr\rangle+\bigl\langle\mu_1^{*} w \!-\! \HH w,w_t\bigr\rangle\!-\!  \int_\Omega\zeta_\delta\dive u^0\dive w_t\, dx .
\end{equation}
\underline{Step 3: integration by parts in time.}
The last term in \eqref{eq:wenergy} cannot be estimated directly because
$\dive w_t$ is not controlled by $\|w_t\|_2$. We therefore use
$$-\int_\Omega\zeta_\delta\dive u^0\dive w_t\,dx
=-\frac{d}{dt}\int_\Omega\zeta_\delta\dive u^0\dive w\,dx
+\int_\Omega\zeta_\delta\dive u_t^0\dive w\,dx$$
and introduce the modified energy
\begin{equation}\label{eq:Gfunctional}
E(t)=\tfrac12\|w_t(t)\|_2^2+\tfrac12\|w(t)\|_{e,\delta}^2+\int_\Omega\zeta_\delta\dive u^0(t)\dive w(t)\,dx .
\end{equation}
Since $w(0)=w_t(0)=0$, we have $E(0)=0$ and
\begin{equation}\label{eq:Gprime}
	E'(t)+\alpha\|w_t\|_2^2
	=-\langle f(u^\delta)-f(u^0),w_t\rangle
	+\langle\mu_1^{*} w-\HH w,w_t\rangle
	+\int_\Omega\zeta_\delta\dive u_t^0\dive w\,dx .
\end{equation}

\underline{Step 4: $E$ is equivalent to the energy up to $O(\delta^2)$.}
Set $\Phi(t)=\|w_t(t)\|_2^2+\|\nabla w(t)\|_2^2$. By
\eqref{eq:div-grad}, \eqref{eq:unifH1bound}  and Young  inequality
\begin{equation}\label{eq:crossterm}
	\Bigl|\int_\Omega\zeta_\delta\dive u^0\dive w\,dx\Bigr|
	\leqslant\delta R\|\nabla w\|_2
	\leqslant\frac{\underline\mu}{4}\|\nabla w\|_2^2
	+\frac{R^2}{\underline\mu}\delta^2 .
\end{equation}
Combining \eqref{eq:crossterm} with \eqref{eq:unifcoerc}, we obtain
\begin{equation}\label{eq:Gequiv}
c_1\Phi(t)-c_2\delta^2	\leqslant E(t)	\leqslant c_3\Phi(t)+c_2\delta^2,
\end{equation}
where $c_1=\min\{\tfrac12,\tfrac{\underline\mu}{4}\}$,
$c_3=\max\{\tfrac12,\tfrac{{C_{\rm eq}}}{2}
+\tfrac{\underline\mu}{4}\}$ and
$c_2=R^2/\underline\mu$ are independent of $\delta$.

\underline{Step 5: estimate of the right-hand side of \eqref{eq:Gprime}.}
Using \eqref{eq:div-grad}, \eqref{eq:unifH1bound} and Young inequality, we obtain
\begin{equation}\label{eq:term2}
\Bigl|\int_\Omega\zeta_\delta\dive u_t^0\dive w\,dx\Bigr|
\leqslant\delta R\|\nabla w\|_2
\leqslant\frac12\|\nabla w\|_2^2+\frac{R^2}{2}\delta^2.
\end{equation}
This is where the uniform boundedness of the trajectories in $\mathscr H^1$ is essential and explains the restriction $\U_0\in\mathbf B$. Moreover, Lemma \ref{lem:lipschitz} applied on a fixed ball of $\mathscr H^0$ of radius $R$, gives $\|f(u^\delta)-f(u^0)\|_2\leqslant C_R\|\nabla w\|_2$. Hence,
\begin{equation}\label{eq:term1}
\bigl|\langle f(u^\delta)-f(u^0),w_t\rangle\bigr|
\leqslant C_R\|\nabla w\|_2\|w_t\|_2
\leqslant\frac{C_R}{2}\Phi(t).
\end{equation}
Finally, $\|\HH w\|_2\leqslant6\mu_2\|w\|_2$. Thus, Poincar\'e  and
Young  inequalities yield
\begin{equation}\label{eq:term3}
	\bigl|\langle\mu_1^{*} w-\HH w,w_t\rangle\bigr|
	\leqslant\bigl[\mu_1^{*}+6\mu_2\bigl]\|w\|_2\|w_t\|_2
	\leqslant
	\frac{\mu_1^{*}+6\mu_2}{2\sqrt{\lambda_1}}\Phi(t).
\end{equation}

\underline{Step 6: conclusion of the energy estimate.}
Inserting \eqref{eq:term2}, \eqref{eq:term1} and \eqref{eq:term3} into
\eqref{eq:Gprime} and dropping the nonnegative damping term, we obtain
$$E'(t)\leqslant c_4\Phi(t)+c_5\delta^2,$$
where
$c_4=\frac{C_R}{2}+\frac12+\frac{\mu_1^{*}+6\mu_2}{2\sqrt{\lambda_1}}$ and $c_5=\frac{R^2}{2}$.
The lower bound in \eqref{eq:Gequiv} gives $\Phi(t)\leqslant c_1^{-1}\bigl[E(t)+c_2\delta^2\bigl]$ and therefore
$$
E'(t)
\leqslant
\frac{c_4}{c_1}E(t)
+\left[\frac{c_4c_2}{c_1}+c_5\right]\delta^2.
$$
Since $E(0)=0$, Gronwall inequality yields, for $0\leqslant t\leqslant T$,
$$
E(t)
\leqslant
\frac{c_1}{c_4}
\left[\frac{c_4c_2}{c_1}+c_5\right]
\left(e^{c_4T/c_1}-1\right)\delta^2
=:c_6(T)\delta^2.
$$
Using again \eqref{eq:Gequiv}, we conclude that
$\Phi(t)\leqslant c_1^{-1}[c_6(T)+c_2]\delta^2$ for
$0\leqslant t\leqslant T$. Since
$$\|S_\delta(t)\U_0-S_0(t)\U_0\|_{\mathscr H^0}^2=\Phi(t)$$ by
\eqref{eq:fixedspaces}, the assertion follows with
$C(T)=[(c_6(T)+c_2)/c_1]^{1/2}$.
\end{proof}
The finite-time estimate above controls the difference between the
flows over one time step. Combined with the uniform Lipschitz continuity
of the fixed-time maps, it yields the parameter-continuity estimate for
their iterates required by Theorem \ref{thm:robust-abstract}.
\begin{Lemma}\label{lem:iterclose}	Assume {\rm(H8)} and let $L_\delta$ be the fixed-time maps introduced	in Lemma $\ref{lem:unifdyn4}$. Then there exist constants $C>0$ and 	$L>1$, independent of $\delta$, such that
\begin{equation*}
\bigl\|L_\delta^{\,k}\U_0-L_0^{\,k}\U_0\bigr\|_{\mathscr H^0} \leqslant C\delta L^k,\,\,\,\forall \,k\in\mathbb N,\,\, \U_0\in\mathbf B,\,\,\delta\in[0,1].
\end{equation*}
\end{Lemma}
\begin{proof}
By \eqref{eq:cont-dep} and Lemma \ref{lem:unifdyn1}, the maps $L_\delta$ are Lipschitz continuous on $\mathbf B$, uniformly in	$\delta$. Thus, after increasing the constant if necessary, there exists	$L=L(t_*)>1$, independent of $\delta$, such that
$$\|L_\delta\U^1-L_\delta\U^2\|_{\mathscr H^0}\leqslant L\|\U^1-\U^2\|_{\mathscr H^0}, \,\,\, \U^1,\U^2\in\mathbf B.$$
Let $C_0=C(t_*)$ be the constant in Lemma \ref{thm:close}. We prove inductively that
\begin{equation}\label{induction}
\|L_\delta^{\,k}\U_0-L_0^{\,k}\U_0\|_{\mathscr H^0}\leqslant C_0\delta\sum_{j=0}^{k-1}L^j, \,\,\, \U_0\in\mathbf B.
\end{equation}
The estimate is immediate for $k=1$ from Lemma \ref{thm:close}. Assume that it holds for some $k\geqslant1$. Since $\mathbf B$ is invariant under $L_\delta$ and $L_0$, both $L_\delta^{\,k}\U_0$ and $L_0^{\,k}\U_0$ belong to $\mathbf B$. Hence
\begin{align*}
\|L_\delta^{\,k+1}\U_0-L_0^{\,k+1}\U_0\|_{\mathscr H^0}
&\leqslant	\|L_\delta(L_\delta^{\,k}\U_0) -L_\delta(L_0^{\,k}\U_0)\|_{\mathscr H^0}+
\|L_\delta(L_0^{\,k}\U_0)-L_0(L_0^{\,k}\U_0)\|_{\mathscr H^0}\\
&\leqslant L\|L_\delta^{\,k}\U_0-L_0^{\,k}\U_0\|_{\mathscr H^0} +C_0\delta\\
&\leqslant C_0\delta\sum_{j=0}^{k}L^j.
\end{align*}
Thus \eqref{induction} holds for every $k\in\mathbb N$.	Since $L>1$, $\sum_{j=0}^{k-1}L^j=\frac{L^k-1}{L-1}\leqslant\frac{L^k}{L-1}$. The conclusion follows with $C=C_0/(L-1)$.
\end{proof}
\begin{Lemma}[Uniform exponential attraction to $\mathbf B$]\label{lem:unifattr} Assume {\rm(H8)}. The absorbing ball $\mathbf B_0\subset\mathscr H^0$ of the Proposition $\ref{prop:absorbing}$ can be chosen independently of $\delta$. Moreover, for every bounded set $D\subset\mathscr H^0$, there exists an entrance time ${t_D}\geqslant0$, independent of $\delta$, such that $S_\delta(t){D}\subseteq\mathbf B_0$ for all $t\geqslant t_D$ and	$\delta\in[0,1]$. In addition, there exist constants $N_1,\nu_1>0$, independent of	$\delta$, such that
	$$
	\dist_{\mathscr H^0} \bigl(S_\delta(t)\mathbf B_0,\mathbf B\bigr)
	\leqslant N_1e^{-\nu_1t}, \,\,\,\forall \, t\geqslant0,\, \forall\, \delta\in[0,1].
	$$
\end{Lemma}
\begin{proof}
	The radius of the absorbing ball in Proposition \ref{prop:absorbing} and the corresponding entrance times depend only on the constants in Proposition \ref{prop:energy}, on $\mu_2$ and $\mu_1^{*}$ and on $\underline\mu$. These quantities are independent of $\delta$ by  Lemma \ref{lem:unifdyn1} and \eqref{eq:unifcoerc}. Hence, $\mathbf B_0$ can be chosen independently of $\delta$ and for every	bounded set $D\subset\mathscr H^0$, its entrance time $t_D$ can also be	chosen uniformly with respect to $\delta$.

	For the exponential attraction, we follow the bootstrap construction in the proof of Theorem \ref{thm:B1}. The finite chain	$0=k_0<k_1<\cdots<k_n=1$ depends only on $p$. At each level, the radii provided by Lemma \ref{lem:comp} depend only on the structural constants $\mu_2$, $\mu_1^{*}$ and $C_{\rm ell}$ of Proposition \ref{prop:unifH2}. Hence they are	independent of $\delta$. By the choice of $R_\star$ in	Lemma \ref{lem:unifdyn3}, the set obtained at the final level is contained in $\mathbf B$. Moreover, the linear decay estimate used at every step is uniform by Lemma \ref{lem:unifdyn1}, while the Lipschitz estimate \eqref{eq:cont-dep} also holds with a constant independent of $\delta$. Applying Lemma \ref{lem:trans} along this finite chain and using
	$\mathbf{B}_{k_n}\subseteq\mathbf B$ yields the conclusion, with constants independent of $\delta$.
\end{proof}
\subsubsection{Completion of the Proof of Theorem \ref{thm:robust}}
\begin{proof}
\underline{Step 1. The discrete family.}
We apply Theorem \ref{thm:robust-abstract} with $H_0=\mathscr H^0$, $H_1=\mathscr H^1$, the compact set $\mathbf B$ given by Lemma \ref{lem:unifdyn3} and the fixed-time maps $L_\delta$
introduced in Lemma \ref{lem:unifdyn4}. The embedding $\mathscr H^1\subset\mathscr H^0$ is compact, while hypotheses $1$ and $2$ follow from Lemmas \ref{lem:unifdyn4} and
\ref{lem:iterclose}, respectively. Hence there exist exponential attractors $\mathbf M_\delta\subset\mathbf B$ for $L_\delta$, with uniformly bounded fractal dimension and H\"older continuous dependence on $\delta$, such that
\begin{equation}\label{eq:attrB}
\dist_{\mathscr H^0}
\bigl(L_\delta^{\,n}\mathbf B,\mathbf M_\delta\bigr)
\leqslant C_0e^{-an},\,\,\, n\in\mathbb N,
\end{equation}
where $C_0,a>0$ are independent of $\delta$.

\underline{Step 2. Passage to continuous time.} Following the classical construction of \cite{EFNT}, used in the parameter-dependent setting in \cite[(4.53)]{FGMZ}, define $
\mathscr{F}_\delta:[0,t_*]\times\mathbf M_\delta\longrightarrow\mathscr H^0$ by  ${\mathscr F_\delta}(t,\U)=S_\delta(t)\U$ and set
$$
\mathbf A_\delta^{\exp}
:=
{\mathscr F_\delta}([0,t_*]\times\mathbf M_\delta)
=
\bigcup_{t\in[0,t_*]}S_\delta(t)\mathbf M_\delta.
$$
Since ${\mathscr F_\delta}$ is continuous (Lemma \ref{lem:criterion}$(1)$) and $[0,t_*]\times\mathbf M_\delta$ is compact, $\mathbf A_\delta^{\exp}$ is compact in $\mathscr H^0$.
The set $\mathbf A_\delta^{\exp}$ is positively invariant. Indeed, let
$\U=S_\delta(t)\V$, where $t\in[0,t_*]$ and
$\V\in\mathbf M_\delta$. For $s\geqslant0$, write
$s+t=nt_*+\tau$, where $n$ is a nonnegative integer and
$\tau\in[0,t_*)$. Since $L_\delta^n\mathbf M_\delta\subseteq\mathbf M_\delta$, the semigroup
property gives $S_\delta(s)\U
=
S_\delta(\tau)L_\delta^n\V
\in
S_\delta(\tau)\mathbf M_\delta
\subseteq
\mathbf A_\delta^{\exp}$. Thus
$$S_\delta(s)\mathbf A_\delta^{\exp}
\subseteq\mathbf A_\delta^{\exp}\,\,\,\forall\, s\geqslant0.$$
By Lemmas \ref{lem:criterion}$(1)$ and \ref{lem:unifdyn4}, ${\mathscr F_\delta}$ is Lipschitz continuous, with a constant independent of $\delta$. Since Lipschitz maps do not increase fractal dimension, we obtain $\dim_F^{\mathscr H^0}\bigl(\mathbf A_\delta^{\exp}\bigr)
\leqslant 1+\dim_F^{\mathscr H^0}\bigl(\mathbf M_\delta\bigr)<\infty$. Consequently,
$$
\sup_{\delta\in[0,1]}\dim_F^{\mathscr H^0}\bigl(\mathbf A_\delta^{\exp}\bigr)<\infty.
$$
Finally, for $t\geqslant0$, write $t=nt_*+\tau$, where $n$ is a
nonnegative integer and $\tau\in[0,t_*)$. Since
$S_\delta(\tau)\mathbf M_\delta\subseteq\mathbf A_\delta^{\exp}$, the
uniform Lipschitz continuity of $S_\delta(\tau)$ with respect to the
initial data and \eqref{eq:attrB} yield
$$
\dist_{\mathscr H^0}\bigl(S_\delta(t)\mathbf B,\mathbf A_\delta^{\exp}\bigr)
\leqslant
\dist_{\mathscr H^0} \bigl(S_\delta(\tau)L_{\delta}^{n}\mathbf B,S_\delta(\tau)\mathbf M_\delta\bigr)
\leqslant C\dist_{\mathscr H^0}\bigl(L_\delta^n\mathbf B,\mathbf M_\delta\bigr)\leqslant CC_0e^{-an}.
$$
Since $n\geqslant t/t_*-1$, we conclude that
\begin{equation}\label{eq:attrBcont}
\dist_{\mathscr H^0} \bigl(S_\delta(t)\mathbf B,\mathbf A_\delta^{\exp}\bigr)\leqslant C_0'e^{-\nu t},\,\,\, \forall\, t\geqslant0,
\end{equation}
where $C_0'>0$ is independent of $\delta$ and $\nu=a/t_*$.

\underline{Step 3. The H\"older estimate.}
Let $\U\in\mathbf M_\delta$ and $t\in[0,t_*]$. By the H\"older
estimate obtained in Step 1, there exists $\U^0\in\mathbf M_0$ such that $\|\U-\U^0\|_{\mathscr H^0}\leqslant C_5\delta^\tau$. Since $S_0(t)\U^0\in\mathbf A_0^{\exp}$, Lemma \ref{thm:close} and \eqref{eq:cont-dep} yield
\begin{align*}
	\|S_\delta(t)\U-S_0(t)\U^0\|_{\mathscr H^0}
	&\leqslant
	\|S_\delta(t)\U-S_0(t)\U\|_{\mathscr H^0}
	+
	\|S_0(t)\U-S_0(t)\U^0\|_{\mathscr H^0}\\
	&\leqslant
	C(t_*)\delta+K(t_*)C_5\delta^\tau
	\leqslant C\delta^\tau,
\end{align*}
where $0<\tau<1$ and $\delta\in[0,1]$. Taking the supremum over
$t\in[0,t_*]$ and $\U\in\mathbf M_\delta$ gives
$$
\dist_{\mathscr H^0}
\bigl(\mathbf A_\delta^{\exp},\mathbf A_0^{\exp}\bigr)
\leqslant C\delta^\tau.
$$
The same argument, using the reverse Hausdorff estimate between
$\mathbf M_0$ and $\mathbf M_\delta$, gives
$$
\dist_{\mathscr H^0}
\bigl(\mathbf A_0^{\exp},\mathbf A_\delta^{\exp}\bigr)
\leqslant C\delta^\tau.
$$
Therefore,
$$
\dist_{\rm sym}
\bigl(\mathbf A_\delta^{\exp},\mathbf A_0^{\exp}\bigr)
\leqslant C\delta^\tau.
$$

\underline{Step 4. Exponential attraction and inclusion of the global
	attractor.}
Since $\mathbf B$ need not be absorbing in $\mathscr H^0$, the
exponential attraction of arbitrary bounded sets does not follow
directly from \eqref{eq:attrBcont}. We first combine it with Lemma
\ref{lem:unifattr}.
Apply Lemma \ref{lem:trans} with
$\mathbf{K}_1=\mathbf B_0$, $\mathbf{K}_2=\mathbf B$ and
$\mathbf{K}_3=\mathbf A_\delta^{\exp}$. Its hypotheses follow from Lemma
\ref{lem:unifattr}, \eqref{eq:attrBcont} and the uniform
continuous-dependence estimate \eqref{eq:cont-dep}. Hence there exist
$N_2,\nu_2>0$, independent of $\delta$, such that
\begin{equation}\label{eq:B0exp}
\dist_{\mathscr H^0}
\bigl(S_\delta(t)\mathbf B_0,\mathbf A_\delta^{\exp}\bigr)\leqslant N_2e^{-\nu_2t},
	\,\,\,t\geqslant0,\,\,\delta\in[0,1].
\end{equation}
In particular, one may take
$\nu_2=\nu_1\nu/(\vartheta_0+\nu_1+\nu)$.
Let $D\subset\mathscr H^0$ be bounded and set
$r_D=\sup_{\U\in D}\|\U\|_{\mathscr H^0}$. By Lemma
\ref{lem:unifattr}, there exists an entrance time $t_D\geqslant0$,
depending only on $r_D$, such that
$S_\delta(t_D)D\subseteq\mathbf B_0$ for every $\delta$. Thus, for
$t\geqslant t_D$, \eqref{eq:B0exp} gives
\begin{align*}
	\dist_{\mathscr H^0}\bigl(S_\delta(t)D,\mathbf A^{\exp}_\delta\bigr)
	&=\dist_{\mathscr H^0}\bigl(S_\delta(t-t_D)S_\delta(t_D)D,\mathbf A^{\exp}_\delta\bigr)\\
	&\leqslant\dist_{\mathscr H^0}\bigl(S_\delta(t-t_D)\mathbf{B}_0,\mathbf A^{\exp}_\delta\bigr)
	\ \leqslant\ N_2e^{\nu_2 t_D}e^{-\nu_2 t}.
\end{align*}
It remains to consider $0\leqslant t\leqslant t_D$. By Corollary
\ref{cor:bounded} and Lemma \ref{lem:unifdyn1},
$\|S_\delta(t)\U\|_{\mathscr H^0}\leqslant R({r_D})$ for every
$\U\in D$. Moreover, since $\mathbf M_\delta\subseteq\mathbf B$, the
same results and the definition of $\mathbf A_\delta^{\exp}$ give
$\|\V\|_{\mathscr H^0}\leqslant R({r_{\mathbf B}})$ for every
$\V\in\mathbf A_\delta^{\exp}$, where
$r_{\mathbf B}=\sup_{\U\in\mathbf B}\|\U\|_{\mathscr H^0}$. Hence, by the triangle inequality,
$$
\|S_\delta(t)\U-\V\|_{\mathscr H^0}
\leqslant R(r_D)+R(r_{\mathbf B}),\,\,\,
\U\in D,\,\,\V\in\mathbf A_\delta^{\exp}.
$$
Taking the infimum over $\V\in\mathbf A_\delta^{\exp}$, then the supremum over
$\U\in D$ and using that $e^{\nu_2(t_D-t)}\geqslant1$ for $0\leqslant t\leqslant t_D$, we obtain
$$
\dist_{\mathscr H^0}
\bigl(S_\delta(t)D,\mathbf A_\delta^{\exp}\bigr)
\leqslant
\bigl[R(r_D)+R(r_{\mathbf B})\bigr]
e^{\nu_2t_D}e^{-\nu_2t}.
$$
Combining both estimates, we obtain
$$
\dist_{\mathscr H^0}
\bigl(S_\delta(t)D,\mathbf A_\delta^{\exp}\bigr)
\leqslant C_De^{-\nu_2t},
\,\,\, t\geqslant0,
$$
where $C_D$ is independent of $\delta$. Therefore, $\mathbf A_\delta^{\exp}$ is an exponential attractor for $(\mathscr H^0,S_\delta(t))$.
Finally, Lemma \ref{lem:unifdyn3} gives
$\mathbf A_\delta\subseteq\mathbf B$. Using the invariance of
$\mathbf A_\delta$ and \eqref{eq:attrBcont}, we obtain
$$
\dist_{\mathscr H^0}
\bigl(\mathbf A_\delta,\mathbf A_\delta^{\exp}\bigr)
=
\dist_{\mathscr H^0}
\bigl(S_\delta(t)\mathbf A_\delta,\mathbf A_\delta^{\exp}\bigr)
\leqslant
\dist_{\mathscr H^0}
\bigl(S_\delta(t)\mathbf B,\mathbf A_\delta^{\exp}\bigr)
\xrightarrow[t\to\infty]{}0.
$$
Since $\mathbf A_\delta^{\exp}$ is closed,
$\mathbf A_\delta\subseteq\mathbf A_\delta^{\exp}$.
\end{proof}
The finite-time comparison of the semigroups and the common bound for
the global attractors also yield the following stability property.
\begin{Theorem}[Upper semicontinuity of the global attractors] \label{thm:usc}
Assume {\rm(H8)}. Then
	$$
	\lim_{\delta\to0^+}
	\dist_{\mathscr H^0}
	\bigl(\mathbf A_\delta,\mathbf A_0\bigr)=0.
	$$
\end{Theorem}
\begin{proof} Let $\epsilon>0$. Since $\mathbf A_0$ attracts the bounded set $\mathbf B$ under $S_0$, there exists $T>0$ such that
$$\dist_{\mathscr H^0}\bigl(S_0(T)\mathbf B,\mathbf A_0\bigr)<\frac{\epsilon}{2}.$$
By Lemma \ref{lem:unifdyn3}, $\mathbf A_\delta\subseteq\mathbf B$ and the invariance of $\mathbf A_\delta$ gives $\mathbf A_\delta=S_\delta(T)\mathbf A_\delta$. Hence, by Lemma \ref{thm:close},
\begin{align*}
\dist_{\mathscr H^0}\bigl(\mathbf A_\delta,\mathbf A_0\bigr)&=\dist_{\mathscr H^0}\bigl(S_\delta(T)\mathbf A_\delta,\mathbf A_0\bigr)\\
&\leqslant\sup_{\U_0\in\mathbf A_\delta}\bigl\|S_\delta(T)\U_0-S_0(T)\U_0\bigr\|_{\mathscr H^0}
+\dist_{\mathscr H^0}\bigl(S_0(T)\mathbf A_\delta,\mathbf A_0\bigr)\\
&\leqslant C(T)\delta+\dist_{\mathscr H^0}\bigl(S_0(T)\mathbf B,\mathbf A_0\bigr)<		C(T)\delta+\frac{\epsilon}{2}.
\end{align*}
Choosing $\delta<\epsilon/[2C(T)]$ completes the proof.
\end{proof}
\begin{Remark}\label{rem:usc-rate}
The H\"older continuity established in Theorem \ref{thm:robust} is quantitative and two-sided, whereas Theorem \ref{thm:usc} provides only qualitative upper semicontinuity. This difference reflects the distinct nature of the two objects. Global attractors are unique, but quantitative continuity generally requires additional structural	assumptions, such as the hyperbolicity of the equilibria in $\mathscr N$. Exponential attractors are not unique and this 	flexibility permits a selection that depends H\"older continuously on the parameter. Moreover, the inclusions 	$\mathbf A_\delta\subseteq\mathbf A_\delta^{\exp}$ and the H\"older		estimate of Theorem \ref{thm:robust} imply
$$
\dist_{\mathscr H^0}\bigl(\mathbf A_\delta,\mathbf A_0^{\exp}\bigr)\leqslant C\delta^\tau.
$$
This does not provide a convergence rate toward $\mathbf A_0$, since $\mathbf A_0^{\exp}$ may be strictly larger than $\mathbf A_0$.
\end{Remark}

\subsection*{Data availability}
\noindent  Not applicable. No datasets were generated or analyzed during this study.
\subsection*{Competing interests}
\noindent The authors declare that they have no competing interests.
\subsection*{Acknowledgements}
\noindent The first author gratefully acknowledges the Departamento de Matem\'atica Aplicada a las TIC (DMATIC) at the Universidad Polit\'ecnica de Madrid for its hospitality and financial support during his research stay.
\subsection*{Funding}
\noindent L. E. Bocanegra-Rodr\'iguez and P. N. Seminario-Huertas were supported
by the Madrid Government (Comunidad de Madrid, Spain) under the
Multiannual Agreement 2023-2026 with the Universidad Polit\'ecnica de
Madrid, Line A, Emerging PhD Researchers
(DOCTORES-EMERGENTES-24-UMDLRU-108-XVLWPC). All authors were also
supported by the Spanish Ministerio de Ciencia, Innovaci\'on y Universidades under the project PID2025-174633NA-I00.


\begin{thebibliography}{99}
\bibitem{adam} L. Adam, M. Batzle and I. Brevik.
\newblock {Gassmann's fluid substitution and shear modulus variability in carbonates at laboratory seismic and ultrasonic frequencies}.
\newblock {\em Geophysics}, Vol. 71(6), F173-F183, 2006.
\bibitem{ADN} S. Agmon, A. Douglis and L. Nirenberg.
\newblock {Estimates near the boundary for solutions of elliptic partial differential equations satisfying general boundary conditions II}.
\newblock {\em Comm. Pure Appl. Math.}, Vol. 17, 35-92, 1964.
\bibitem{arrieta} J. Arrieta, A. N. Carvalho and J. K. Hale.
\newblock {A damped hyperbolic equation with critical exponent}.
\newblock {\em Comm. Partial Differential Equations}, Vol. 17(5-6), 841-866, 1992.
\bibitem{BLR} C. Bardos, G. Lebeau and J. Rauch.
\newblock {Sharp sufficient conditions for the observation, control, and stabilization of waves from the boundary}.
\newblock {\em SIAM J. Control Optim.}, Vol. 30(5), 1024-1065, 1992.
\bibitem{bellassoued} M. Bellassoued.
\newblock {Energy decay for the elastic wave equation with a local time-dependent nonlinear damping}.
\newblock {\em Acta Math. Sin. (Engl. Ser.)}, Vol. 24(7), 1175-1192, 2008.
\bibitem{bocanegra} L. E. Bocanegra-Rodr\'iguez, M. A. Jorge Silva, T. F. Ma and P. N. Seminario-Huertas.
\newblock {Longtime dynamics of a semilinear Lam\'e system}.
\newblock {\em J. Dynam. Differential Equations}, Vol. 35, 1435-1456, 2023.
\bibitem{Ma} M. M. Cavalcanti, L. H. Fatori and T. F. Ma.
\newblock {Attractors for wave equations with degenerate memory}.
\newblock {\em J. Differential Equations}, Vol. 260(1), 56-83, 2016.
\bibitem{Von} I. Chueshov and I. Lasiecka.
\newblock {\em Von Karman Evolution Equations: Well-posedness and Long Time Dynamics}.
\newblock Springer, New York, 2010.
\bibitem{CL-memoirs} I. Chueshov and I. Lasiecka.
\newblock {Long-time behavior of second order evolution equations with nonlinear damping}.
\newblock {\em Mem. Amer. Math. Soc.}, Vol. 195, no. 912, 2008.
\bibitem{EFNT} A. Eden, C. Foias, B. Nicolaenko and R. Temam.
\newblock {\em Exponential Attractors for Dissipative Evolution Equations}.
\newblock RAM 37, Masson, Paris, 1994.
\bibitem{EMZ} M. Efendiev, A. Miranville and S. Zelik.
\newblock {Exponential attractors for a nonlinear reaction-diffusion system in $\mathbb{R}^3$}.
\newblock {\em C. R. Acad. Sci. Paris S\'er. I Math.}, Vol. 330, 713-718, 2000.
\bibitem{EMZ2} M. Efendiev, A. Miranville and S. Zelik.
\newblock {Exponential attractors for a singularly perturbed Cahn-Hilliard system}.
\newblock {\em Math. Nachr.}, Vol. 272, 11-31, 2004.
\bibitem{FGMZ} P. Fabrie, C. Galusinski, A. Miranville and S. Zelik.
\newblock {Uniform exponential attractors for a singularly perturbed damped wave equation}.
\newblock {\em Discrete Contin. Dyn. Syst.}, Vol. 10(1-2), 211-238, 2004.
\bibitem{fujiwara} D. Fujiwara.
\newblock {Concrete characterization of the domains of fractional powers of some elliptic differential operators of the second order}.
\newblock {\em Proc. Japan Acad.}, Vol. 43, 82-86, 1967.
\bibitem{gassmann} F. Gassmann.
\newblock {\"Uber die Elastizit\"at por\"oser Medien}.
\newblock {\em Vierteljahrsschr. Naturforsch. Ges. Z\"urich}, Vol. 96, 1-23, 1951.
\bibitem{HR} K. W. Hajduk and J. C. Robinson.
\newblock {Explicit characterisation of the fractional power spaces of the Dirichlet Laplacian and Stokes operators}.
\newblock arXiv:2108.03139, 2021.
\bibitem{hale} J. K. Hale.
\newblock {\em Asymptotic Behavior of Dissipative Systems}.
\newblock Math. Surveys Monogr. 25, Amer. Math. Soc., Providence, 1988.
\bibitem{heinz} E. Heinz.
\newblock {Beitr\"age zur St\"orungstheorie der Spektralzerlegung}.
\newblock {\em Math. Ann.}, Vol. 123, 415-438, 1951.
\bibitem{john} F. John.
\newblock {Almost global existence of elastic waves of finite amplitude arising from small initial disturbances}.
\newblock {\em Comm. Pure Appl. Math.}, Vol. 41(5), 615-666, 1988.
\bibitem{kato61} T. Kato.
\newblock {A generalization of the Heinz inequality}.
\newblock {\em Proc. Japan Acad.}, Vol. 37, 305-308, 1961.
\bibitem{kato-book} T. Kato.
\newblock {\em Perturbation Theory for Linear Operators}, 2nd ed.
\newblock Grundlehren Math. Wiss. 132, Springer, Berlin, 1980.
\bibitem{Kolsky} H. Kolsky.
\newblock {\em Stress Waves in Solids}.
\newblock Dover, New York, 1963.
\bibitem{ladyzhenskaya} O. Ladyzhenskaya.
\newblock {\em Attractors for Semigroups and Evolution Equations}.
\newblock Cambridge University Press, Cambridge, 1991.
\bibitem{lakes} R. S. Lakes, T. Lee, A. Bersie and Y. C. Wang.
\newblock {Extreme damping in composite materials with negative-stiffness inclusions}.
\newblock {\em Nature}, Vol. 410, 565-567, 2001.
\bibitem{lakesdrugan} R. S. Lakes and W. J. Drugan.
\newblock {Dramatically stiffer elastic composite materials due to a negative stiffness phase?}
\newblock {\em J. Mech. Phys. Solids}, Vol. 50, 979-1009, 2002.
\bibitem{Lame} G. Lam\'e.
\newblock {\em Le\c{c}ons sur la th\'eorie math\'ematique de l'\'elasticit\'e des corps solides}.
\newblock Gauthier-Villars, Paris, 1866.
\bibitem{LM} J.-L. Lions and E. Magenes.
\newblock {\em Non-Homogeneous Boundary Value Problems and Applications, Vol. I}.
\newblock Springer, Berlin, 1972.
\bibitem{Love} A. E. H. Love.
\newblock {\em A Treatise on the Mathematical Theory of Elasticity}, 4th ed.
\newblock Dover, New York, 1944.
\bibitem{singular} T. F. Ma and R. N. Monteiro.
\newblock {Singular limit and long-time dynamics of Bresse systems}.
\newblock {\em SIAM J. Math. Anal.}, Vol. 49(4), 2468-2495, 2017.
\bibitem{mclean} W. McLean.
\newblock {\em Strongly Elliptic Systems and Boundary Integral Equations}.
\newblock Cambridge University Press, Cambridge, 2000.
\bibitem{necas} J. Ne\v{c}as.
\newblock {\em Direct Methods in the Theory of Elliptic Equations}.
\newblock Springer Monographs in Mathematics. Springer, Berlin, 2012.
\newblock Translated from the 1967 French original \emph{Les m\'ethodes directes en th\'eorie des
\'equations elliptiques} by G. Tronel and A. Kufner.
\bibitem{ohno} I. Ohno, K. Harada and C. Yoshitomi.
\newblock {Temperature variation of elastic constants of quartz across the $\alpha$-$\beta$ transition}.
\newblock {\em Phys. Chem. Minerals}, Vol. 33, 1-9, 2006.
\bibitem{PS} V. Pata and M. Squassina.
\newblock {On the strongly damped wave equation}.
\newblock {\em Comm. Math. Phys.}, Vol. 253(3), 511-533, 2005.
\bibitem{pazy} A. Pazy.
\newblock {\em Semigroups of Linear Operators and Applications to Partial Differential Equations}.
\newblock Appl. Math. Sci. 44, Springer, New York, 1983.
\bibitem{peng} Z. Peng and S. A. T. Redfern.
\newblock {Mechanical properties of quartz at the $\alpha$-$\beta$ phase transition: implications for tectonic and seismic anomalies}.
\newblock {\em Geochem. Geophys. Geosyst.}, Vol. 14(1), 18-28, 2013.
\bibitem{Pujol} J. Pujol.
\newblock {\em Elastic Wave Propagation and Generation in Seismology}.
\newblock Cambridge University Press, Cambridge, 2003.
\bibitem{Sideris} T. C. Sideris.
\newblock {The null condition and global existence of nonlinear elastic waves}.
\newblock {\em Invent. Math.}, Vol. 123(2), 323-342, 1996.
\bibitem{smith} T. M. Smith, C. H. Sondergeld and C. S. Rai.
\newblock {Gassmann fluid substitutions: a tutorial}.
\newblock {\em Geophysics}, Vol. 68(2), 430-440, 2003.
\bibitem{temam} R. Temam.
\newblock {\em Infinite-Dimensional Dynamical Systems in Mechanics and Physics}, 2nd ed.
\newblock Appl. Math. Sci. 68, Springer, New York, 1997.
\bibitem{ya1988} K. Yamamoto.
\newblock {Singularities of solutions to the boundary value problems for elastic and Maxwell's equations}.
\newblock {\em Japan. J. Math. (N.S.)}, Vol. 14(1), 119-163, 1988.
\bibitem{yam1989} K. Yamamoto.
\newblock {Exponential energy decay of solutions of elastic wave equations with the Dirichlet condition}.
\newblock {\em Math. Scand.}, Vol. 65(2), 206-220, 1989.
\end{thebibliography}
\end{document}